\documentclass[twoside,11pt]{article}

\usepackage{amsmath, mathtools, amssymb, amsthm}
\usepackage{natbib}
\usepackage{microtype}
\usepackage{enumerate, multicol}
\usepackage{bbm}
\usepackage[linesnumbered,ruled,vlined,algo2e]{algorithm2e}
\usepackage[font=small,labelfont=bf]{caption}
\usepackage{array} 
\usepackage[table]{xcolor}
\usepackage{tikz}
\usetikzlibrary{positioning, arrows, shapes.geometric, calc}
\usepackage{fullpage}
\usepackage[colorlinks=true, linkcolor=blue, urlcolor=blue, citecolor=blue]{hyperref}
\usepackage{algorithm,algorithmic}
\usepackage{breakcites}
\makeatletter
\newenvironment{breakablealgorithm}
  {
    \begin{center}
      \refstepcounter{algorithm}
      \hrule height.8pt depth0pt \kern2pt
      \parskip 0pt
      \renewcommand{\caption}[2][\relax]{
        {\raggedright\textbf{\fname@algorithm~\thealgorithm} ##2\par}%
        \ifx\relax##1\relax 
          \addcontentsline{loa}{algorithm}{\protect\numberline{\thealgorithm}##2}%
        \else 
          \addcontentsline{loa}{algorithm}{\protect\numberline{\thealgorithm}##1}%
        \fi
        \kern2pt\hrule\kern2pt
     }
  }
  {
     \kern2pt\hrule\relax
   \end{center}
  }
\makeatother

\definecolor{gaussianfill}{RGB}{220,230,250}
\definecolor{kssgfill}{RGB}{220,250,230}
\definecolor{usgfill}{RGB}{245,230,250}

\SetKwInput{KwInput}{Input}
\SetKwInput{KwOutput}{Output}
\SetKwRepeat{KwRepeat}{repeat}{until}

\theoremstyle{plain}
\newtheorem{theorem}{Theorem}[section]
\newtheorem{lemma}[theorem]{Lemma}
\newtheorem{corollary}[theorem]{Corollary}
\newtheorem{proposition}[theorem]{Proposition}
\newtheorem{remark}[theorem]{Remark}
\newtheorem{definition}[theorem]{Definition}
\numberwithin{equation}{section}

\newcommand{\EE}{\mathbb{E}}

\newcommand{\NN}{\mathbb{N}}
\newcommand{\PP}{\mathbb{P}}
\newcommand{\QQ}{\mathbb{Q}}
\newcommand{\RR}{\mathbb{R}}

\newcommand{\Smat}{\mathbb{S}}

\newcommand{\conv}{\operatorname{conv}}

\newcommand{\E}{\mathbb{E}}
\newcommand{\Pp}{\mathbb{P}}

\newcommand{\cC}{\mathcal{C}}

\newcommand{\cL}{\mathcal{L}}

\newcommand{\cO}{\mathcal{O}}
\newcommand{\cP}{\mathcal{P}}

\newcommand{\cU}{\mathcal{U}}

\newcommand{\diag}{\operatorname{diag}}
\newcommand{\diam}{\mathrm{diam}}

\newcommand{\tr}{\operatorname{tr}}

\newcommand{\argmax}{\mathop{\mathrm{argmax}}}
\newcommand{\argmin}{\mathop{\mathrm{argmin}}}

\def\T{^\top}

\allowdisplaybreaks

\begin{document}

\begin{center}
{\LARGE Polynomial-Time Near-Optimal Estimation over Certain Type-2 Convex Bodies}

{\large
\begin{center}
Matey Neykov
\end{center}}

{Department of Statistics and Data Science, Northwestern University\\
\texttt{mneykov@northwestern.edu}}
\end{center}

\begin{abstract}
We develop polynomial-time algorithms for near-optimal minimax mean estimation under $\ell_2$-squared loss in a Gaussian sequence model under convex constraints. The parameter space is an origin-symmetric, type-2 convex body $K \subset \mathbb{R}^n$, and we assume additional regularity conditions: specifically, we assume $K$ is well-balanced, i.e., there exist known radii $r, R > 0$ such that $r B_2 \subseteq K \subseteq R B_2$, as well as oracle access to the Minkowski gauge of $K$. Under additional conditions guaranteeing an efficient approximate quadratic-form-maximization oracle on $K$, our procedures achieve the minimax rate up to factors that depend polylogarithmically on the dimension, while remaining computationally efficient.

We further extend our methodology to the linear regression and robust heavy-tailed settings, establishing polynomial-time near-optimal estimators when the constraint set satisfies the regularity conditions above. To the best of our knowledge, these results provide the first general framework for attaining statistically near-optimal performance under such broad geometric constraints while preserving computational tractability.
\end{abstract}

\tableofcontents

\section{Introduction}

Many statistical estimation problems are naturally formulated as
\emph{constrained mean estimation} or \emph{constrained regression} problems.
Given observations corrupted by stochastic noise, one seeks to estimate a
signal known to belong to a prescribed constraint set
$K\subset\mathbb{R}^n$. Of particular interest in this work are
origin-symmetric convex constraints possessing favorable geometric structure,
such as ellipsoids, hyperrectangles, $\ell_p$-balls for $p\geq 2$, and,
more generally, quadratically convex orthosymmetric sets. Such constraint
classes arise naturally in nonparametric smoothing, signal recovery, and
regression. For example, classical Sobolev smoothness classes become
ellipsoidal constraints after expansion in an appropriate orthonormal basis,
leading to the well-known Gaussian sequence formulation of nonparametric
estimation
\cite{pinsker1980optimal, tsybakov2009introduction}. More general
quadratically convex orthosymmetric constraints, including hyperrectangles and
$\ell_p$-balls for $p>2$, were studied in the seminal work of
\cite{donoho1990minimax}. Ellipsoidal constraints have also played an important
role in regression and prediction problems
\cite{goldenshluger2001adaptive,goldenshluger2003optimal,pathak2024noisy}.

From a decision-theoretic perspective, a fundamental goal is to characterize
the \emph{minimax risk} over the constraint class $K$. There is a long
history of such results for structured type-2 classes. Pinsker's classical
theorem characterizes the asymptotically minimax risk over ellipsoids in the
Gaussian sequence model \cite{pinsker1980optimal}, while
\cite{donoho1990minimax} showed that, for quadratically convex orthosymmetric
sets, the minimax risk is closely related to geometric quantities such as
Kolmogorov widths. More recently, a substantially more general theory has
emerged. In particular, \cite{neykov2022minimax} derived sharp minimax rates
for constrained mean estimation in the Gaussian sequence model for essentially
\emph{any} bounded convex constraint, expressing the optimal rate in terms of
the local geometry of $K$. Subsequently,
\cite{prasadan2024information} relaxed convexity to the substantially weaker
assumption of star-shapedness, thereby covering an even broader range of
constraint sets. In the regression setting,
\cite{prasadan2025characterizingminimaxratenonparametric} obtained analogous
minimax characterizations.

Despite this information-theoretic progress, a major obstacle remains.
The general minimax-optimal procedures developed in these works are not
computationally efficient. Their constructions rely on carefully designed
\emph{local packing sets} adapted to the geometry of $K$. Such packings can
contain exponentially many points and are typically difficult to construct or
even approximate. Consequently, although these estimators attain the optimal
statistical rate for very general constraint sets, they need not admit
polynomial-time implementations.

The goal of the present work is therefore different: rather than treating
arbitrary convex or star-shaped constraints, we focus on the important subclass
of origin-symmetric \emph{type-2} convex bodies and ask whether their additional
geometric structure can be exploited to obtain computationally efficient,
near-minimax estimators.

\subsection*{Our contribution}

In this paper, we develop a new framework that yields \emph{polynomial-time, 
near-minimax} estimators for constrained mean estimation, robust mean estimation with heavy-tailed noise, and constrained linear regression. Our approach bypasses the packing-set bottleneck entirely. Instead, we build on a classical inequality of B.~Carl, which relates the \emph{entropy numbers} of a convex body to its \emph{Kolmogorov $n$-widths}. This geometric connection allows us to translate statistical rate calculations into width-approximation problems, avoiding the need to construct packing sets.

The central technical ingredient of our work is a new algorithm for 
efficiently approximating optimal projections of the Kolmogorov widths of the constraint set. These approximations are sufficiently accurate to be used within an iterative estimation procedure whose risk matches the minimax rate up to factors that are polylogarithmic in the dimension (see Remark~\ref{tracking:Upsilon:K:rem}).

There is, however, a price to pay: our methods require nontrivial geometric 
regularity of the constraint set. First, Carl's inequality yields meaningful 
constants only for \emph{type-2} (see \eqref{type2:constant} and the discussion surrounding it for the definition of a type-2 set), origin-symmetric convex bodies.\footnote{In fact, Carl's inequality holds for all origin-symmetric convex bodies, but the resulting constants are too coarse unless the body is type-2.} Second, the efficient approximation of the optimal projections of the Kolmogorov widths requires the constraint set to be \emph{type-2}\footnote{Here and throughout the paper, when we say type-2 convex bodies we mean convex bodies in $\RR^n$ with sufficiently small type-2 constants (see  \eqref{type2:constant}).}, \emph{well-balanced}, and to satisfy several additional structural conditions. While these assumptions cover a wide range of constraint classes used in modern statistics, they do exclude convex classes that are meaningful and are statistically minimax-tractable (at least in theory and some in practice). Examples of sets over which our algorithm cannot work well include rigid bodies such as the $\ell_p$ balls for $1 \leq p < 2$.

On the other hand, we note that the ``dream'' of obtaining a polynomial-time minimax algorithm (even with a polylogarthmic slack) whenever one has polynomial-time access to the Minkowski gauge of the origin-symmetric convex set $K$, without imposing any additional structural assumptions on $K$, is likely unattainable — at least assuming the standard low-degree conjecture for planted submatrix recovery. We demonstrate this in appendix \ref{computational:appendix}.

\subsection{Related Work}

We now review several related works. In the Gaussian sequence model setting, when the constraint set is an ellipsoid, a classical result due to Pinsker \cite{pinsker1980optimal} shows that linear estimators are optimal. One can implement Pinsker's estimator in polynomial time. A similar result holds for linear regression over ellipsoidal constraints \cite{goldenshluger2001adaptive,goldenshluger2003optimal}. Since ellipsoids are type-2 sets, our framework provides alternative polynomial-time, near-minimax-rate-optimal estimators. Notably, there exist ellipsoids over which the \textbf{least squares estimator is grossly suboptimal} (see, e.g., \cite{zhang2013nearly}). We emphasize that the least squares estimator is perhaps the most natural estimator in this problem: it simply projects the data $Y$ onto the set $K$. This fact highlights the need to develop computationally efficient algorithms beyond least squares. Another computationally natural procedure is the Bayes estimator corresponding to the uniform prior on $K$, namely the posterior mean. For a convex set $K$, the posterior density is proportional to
\begin{align*}
\exp\left(-\frac{\|Y-\theta\|_2^2}{2\sigma^2}\right)\mathbf{1}(\theta\in K),
\end{align*}
and is therefore log-concave; consequently, its mean can be approximated in polynomial time using algorithms for sampling truncated Gaussian distributions over convex bodies \cite{cousins2018gaussian}. Nevertheless, on suitably stretched ellipsoids this estimator can exhibit the same unfavorable behavior as the least squares estimator and can therefore also be grossly suboptimal relative to the minimax risk. We demonstrate this in appendix \ref{posterior:mean}. 

A related geometric perspective appears in \cite{mourtada2023universal}, who studies Gaussian location models from the viewpoint of universal coding and sequential prediction. He shows that minimax redundancy is characterized by global covering numbers, and contrasts this with non-sequential Gaussian mean estimation, whose minimax risk is governed by local covering numbers; our work takes the latter local-entropy characterization as its statistical starting point and exploits type-2 geometry and Carl's inequality to obtain computationally efficient estimators via Kolmogorov-width approximations.

Chandrasekaran and Jordan \cite{chandrasekaran2013computational} study computational--statistical tradeoffs in the Gaussian sequence model by projecting onto tractable convex outer relaxations of the parameter space, with weaker relaxations reducing computation at the cost of increased statistical error or sample size. In contrast, our goal is to obtain polynomial-time, near-minimax estimation over the original constraint set for a broad class of type-2 convex bodies, using Carl's inequality and efficiently approximated Kolmogorov-width projections rather than replacing the constraint by an outer relaxation.

\cite{zhang2013nearly} proposes a near minimax optimal estimator for a linear regression problem with an $\ell_q$ constraints on the signal for $q \leq 1$. The author's algorithm is implementable only when $q = 1$. Our work is significantly different, as we require the constraint set to be of type-2 (which the $\ell_1$ ball is not). Despite this major difference, there are also some similarities: both our work and \cite{zhang2013nearly} use approximations of the Kolmogorov widths of the constraint set, as well as Euclidean projections. \cite{javanmard2012minimax} argue that for bounded symmetric convex polytopes, $K = \{x \in \RR^n : \|Ax \|_{\infty} \leq 1\}$, defined by some matrix $A \in \RR^{m \times n}$, the truncated series estimator\footnote{A truncated series estimator uses the optimal Euclidean projection from a carefully selected Kolmogorov width, to project the data on a lower dimensional subspace.} is nearly optimal for the Gaussian sequence model problem, provided that $m$ is not too large. However, the authors mention that they do not know how to efficiently implement their estimator. They also formulate a class of constraints of the type $\{x \in \RR^n : \|A x\|_p \leq 1\}$ for $p > 0$, although they only show the near optimality result for the truncated series when $p = \infty$. \cite{javanmard2012minimax} do conjecture however (see Conjecture 7.1 therein), that the truncated series estimator remains near optimal (within $\log m$ of optimal) for sets of the type $\{x \in \RR^n : \|A x\|_p \leq 1\}$ for $p \geq 2$. While we cannot verify their conjecture, we do show that our estimator provably works near optimally for \textit{any} well-balanced constraint of the type $K = \{x \in \RR^n : \|A x\|_p \leq 1\}$ for any $p \geq 2$ (we lose logarithmic factors in $m$).

Finally, we note a recent related work by \cite{kur2022efficient}, which develops computationally efficient estimators for multidimensional convex and Lipschitz nonparametric regression. The estimators proposed there are fundamentally different from those studied in this paper, as they are specifically tailored to the structure of convex and Lipschitz constraints, whereas our approach is designed to address constrained estimation problems in a considerably more general setting.

\subsection{Definitions and Notation}

We now outline some definitions and commonly used notation. 

We use a lot of absolute constants which we denote with $c, C, L, T$ and variations, whose values can change from line to line. We use $\rho_K$ to denote the Minkowski gauge (also called Minkowski functional) of a convex body $K$. The gauge $\rho_K$ is a norm, whenever $K$ is bounded, origin-symmetric (which we will assume throughout the paper) and absorbing. We use $\|\cdot\|_p$ to denote the standard $\ell_p$ norm for $p \geq 1$. We use $\|\cdot\|$ to denote an arbitrary norm on $\RR^n$. We use $\|\cdot\|_{\operatorname{op}}$, $\|\cdot\|_F$ and $\|\cdot\|_*$ to denote the operator, Frobenius and nuclear norms of a matrix respectively, and we use $\lambda_{\min}(A)$ and $\lambda_{\max}(A)$ to denote the minimum and maximum eigenvalues for symmetric positive semi-definite (psd) matrix $A$. The trace of a matrix is denoted with $\operatorname{tr}$. For two matrices $A,B \in \RR^{m \times n}$, we put $\langle A, B\rangle = \operatorname{tr}(A\T B)$ for the standard Frobenius (Hilbert-Schmidt) dot product. We write $\leq$ and $\geq$ for symmetric psd matrices to mean inequalities in the psd sense. We typically use the letter $\Pi$ to denote projections -- often when indexed by a closed convex body $\Pi_{K}$ means (weak) Euclidean projection on $K$. We write $\operatorname{diag}(A)$ for a square matrix $A$ to extract its diagonal values, while $\operatorname{diag}(v)$ for a vector $v$ means a matrix with diagonal entries equal to $v$ and $0$ entries otherwise. We denote the $n$-dimensional identity matrix with $I_n$, occasionally omitting the index when the dimension is clear from the context. We use $S^{n-1}$ for the unit Euclidean sphere in $\RR^n$, while $B_2(x, r)$ means the closed Euclidean ball $\{y \in \RR^n: \|x-y\|_2 \leq r\}$ centered at $x$, with the shorthand $B_2 = B_2(0,1)$ for the unit Euclidean ball. We will write $a \wedge b$ and $a\vee b$ to denote the minimum and maximum of two real numbers $a, b \in \RR$. Occasionally we will use $O(1)$ and $\Omega(1)$ to denote that a quantity is upper or lower bounded by an absolute constant, respectively. The inequalities $\lesssim$ and $\gtrsim$ suppress absolute constant factors. We write $\asymp$ if both $\lesssim$ and $\gtrsim$ hold.

We continue with the definition of Kolmogorov width which is central to this work.

\begin{definition}[Kolmogorov width]\label{def_Kolmogorov_width}
    Let $K\subset \RR^n$ be a compact set. The $k$-th Kolmogorov width of $K$ with respect to the Euclidean norm is defined as 
\begin{equation}\label{eq: definition of Kolmogorov width}
    d_{k}(K) = \inf_{P \in \cP_k} \sup_{\theta \in K} \left\lVert \theta -P \theta \right\rVert_{2}^{},
\end{equation}
where $\mathcal{P}_{k}$ is the set of all projection operators that project a vector onto a $k$ dimensional subspace of $\RR^n$.
\end{definition}

Next we give a definition of the entropy numbers of a convex body.

\begin{definition}[Entropy numbers]\label{def_entropy_numbers}
    Let $K\subset \RR^n$ be a compact set. The $k$-th entropy number of $K$ with respect to the Euclidean norm is defined as 
\begin{equation}\label{eq: definition of entroy numbers}
    e_{k}(K) = \inf \bigl\{ \epsilon > 0 : \exists x_1, \dots, x_{k} \in \mathbb{R}^n \text{ such that } 
K \subset \bigcup_{i=1}^{k} \bigl( x_i + \epsilon B_2 \bigr) \bigr\}.
\end{equation}
\end{definition}

The local entropy of a body is defined as follows.

\begin{definition}[Local Metric Entropy] \label{definition:local:metric:entropy} Given a set $K$, define the $\eta$-packing number of $K$ as the maximum cardinality $M=M(\eta,K)$ of a set $\{\nu_1,\dots,\nu_{M}\}\subset K$ such that $\|\nu_i-\nu_j\|>\eta$ for all $i\ne j$. Given a fixed sufficiently large constant $c > 0$, define the local metric entropy of $K$ as $\log M_K^{\operatorname{loc}}(\eta)$ where $M_K^{\operatorname{loc}}(\eta) = \sup_{\nu\in K} M(\eta/c, B(\nu,\eta)\cap K)$. When clear from the context we drop the index $K$ from $M_K^{\operatorname{loc}}(\eta)$.
\end{definition}

The local entropy plays a fundamental role in what follows as \cite{neykov2022minimax, prasadan2024information} show it determines the minimax rate for the constraint mean estimation problem.

\subsection{Geometric background: type-2 constants and Carl's inequality}\label{type2:carl:background}

We briefly recall the Banach-space geometry underlying the main argument.  If
$K\subset\RR^n$ is an origin-symmetric convex body, its Minkowski gauge
$\rho_K$ is a norm, and it is useful to regard $K$ as the unit ball of the
normed space
\begin{align*}
    E_K := (\RR^n,\rho_K).
\end{align*}
The relevant geometric parameter is the Rademacher type-2 constant.  We use
the following second-moment convention: $T_2(K)$ is the smallest $T\geq 1$
such that, for every finite collection $x_1,\ldots,x_m\in\RR^n$,
\begin{align}\label{type2:constant}
    \EE_{\epsilon}\rho_K^2\!\left(\sum_{i=1}^m \epsilon_i x_i\right)
    \leq T^2\sum_{i=1}^m \rho_K^2(x_i),
\end{align}
where $\epsilon_1,\ldots,\epsilon_m$ are independent Rademacher random
variables.  Thus a small type-2 constant means that random signed sums behave,
in the norm $\rho_K$, approximately as they do in a Hilbert space.  Carl
\cite{carl1985inequalities} uses the equivalent first-moment definition of
Rademacher type; by the Kahane--Khintchine inequality the two conventions agree
up to an absolute constant, which is immaterial for our purposes.

A few examples help calibrate this assumption.  Type is invariant under
invertible linear changes of coordinates, and hence every ellipsoid has
$T_2(K)=1$.  For $2\leq p<\infty$, let $B_p^n$ be the unit ball of the $\ell_p$ norm. Then
$T_2(B_p^n)\lesssim \sqrt p \wedge \sqrt{\log(1+n)}$ \cite[see equation (8) e.g.]{bhattiprolu2021framework}, which implies
$T_2(B_\infty^n)\lesssim \sqrt{\log(1+n)}$; the same latter bound therefore
holds for hyperrectangles.  In the opposite direction, for $1\leq p<2$,
taking $x_i=e_i$ in \eqref{type2:constant} gives
\begin{align*}
    T_2(B_p^n)\geq n^{1/p-1/2}.
\end{align*}
This illustrates why $\ell_p$ balls with $p<2$ fall outside the regime of this
paper.  More generally, every $n$-dimensional normed space has a finite type-2
constant (indeed $T_2\lesssim\sqrt n$ by John's theorem), so in finite
dimensions our phrase ``type-2 body'' should be understood as a body whose
type-2 constant is sufficiently small, typically constant or
polylogarithmic in $n$.  In Section~\ref{examples:section} we show, in
particular, that the QCO bodies considered there satisfy
$T_2(K)\lesssim\sqrt{\log n}$.

We next recall the form of Carl's inequality that is used in the paper.  For a
bounded linear operator $S:E\to F$, let $\varepsilon_k(S)$ denote the smallest
radius such that $S(B_E)$ can be covered by $k$ translates of that radius
multiple of $B_F$, and define Carl's \emph{entropy modulus}
\begin{align}\label{carl:entropy:modulus}
    g_m(S):=\inf_{k\geq1} k^{1/m}\varepsilon_k(S).
\end{align}
Theorem~5 of \cite{carl1985inequalities} (see also Appendix \ref{appendix:carl:result} where we state the result in our notation), in the real case, states the following.
If $E$ has Rademacher type $p$, $F^*$ (the dual space of $F$) has Rademacher type $q$, and $s_j(S)$
denotes either the $j$-th Kolmogorov number or the $j$-th Gelfand number of
$S$, then
\begin{align}\label{carl:original:form}
    \left(\prod_{j=1}^m s_j(S)\right)^{1/m}
    \lesssim
    \tau_p(E)\tau_q(F^*)\,
    m^{-1+1/p+1/q}\, g_m(S).
\end{align}
Here $\tau_p$ denotes Carl's first-moment Rademacher type-$p$ constant; when
$p=2$ it is equivalent, up to an absolute factor, to the second-moment
constant $T_2$ in \eqref{type2:constant}. The case relevant here is $p=q=2$,
for which the power of $m$ disappears.
By contrast, without any type assumption Carl's universal bound has an
additional factor $m$:
\begin{align}\label{carl:universal:form}
    \left(\prod_{j=1}^m s_j(S)\right)^{1/m}
    \lesssim m\,g_m(S).
\end{align}
This is the precise sense in which the type-2 hypothesis makes Carl's
inequality quantitatively useful for us.

To specialize \eqref{carl:original:form} to a convex body, take
$E=E_K$, $F=(\RR^n,\|\cdot\|_2)$ (thus $F^* = F$), and let $S$ be the identity map.  Then
$S(B_E)=K$, the dual space $F^*$ is again Euclidean and has type-2 constant
one, and $\varepsilon_k(S)=e_k(K)$ with the convention of Definition
\ref{def_entropy_numbers}.  Carl's operator-theoretic Kolmogorov numbers are
indexed so that $d_j(S)$ corresponds to approximation by a subspace of
dimension strictly smaller than $j$; this differs by one from Definition
\ref{def_Kolmogorov_width}.  Using monotonicity to absorb this harmless index
shift gives
\begin{align}\label{carls:ineq}
    \left(\prod_{j=1}^m d_j(K)\right)^{1/m}
    \lesssim T_2(K)\inf_{k\geq1} k^{1/m}e_k(K).
\end{align}
We emphasize that our $e_k(K)$ is the \emph{non-dyadic} covering convention:
it uses $k$ balls.  Many modern texts instead reserve ``entropy number'' for
the dyadic quantity obtained using $2^{k-1}$ balls; Carl's entropy modulus
\eqref{carl:entropy:modulus}, and all entropy numbers in this paper, use the
former convention.

\subsection{Organization}

The paper is structured as follows. In Section \ref{background:cond:algo:section} we give background results and we state and prove our main result conditionally on the existence of a quadratic form maximization (QFM) oracle for the convex body $K$. In Section \ref{QFM:section} we argue that indeed such QFM oracles exist for a variety of type-2 convex bodies. In Section \ref{robust:section} we extend the results of Section \ref{background:cond:algo:section} to the problem of robust heavy-tailed mean estimation. In Section \ref{regression:section} we discuss the extension of our main result to a linear regression setting. In Section \ref{examples:section} we make some remarks about sets on which our algorithm can be implemented. We close out the paper with a brief discussion of selected open problems in Section \ref{discussion:section}.

\section{Conditional Algorithm}\label{background:cond:algo:section}

Let us observe a single data point $Y = \mu + \xi$ where $\mu \in K \subset \RR^n$ for some known compact convex set $K$, and $\xi=(\xi_1,\ldots,\xi_n)$ is a mean-zero sub-Gaussian random vector with i.i.d. components. We assume that the sub-Gaussian parameter of $\xi_i$, $\sigma^2$ (or an accurate upper bound) is known to the statistician. Here by sub-Gaussian parameter we understand the smallest constant $\sigma^2 \in \RR_+$ such that for all $\lambda \in \RR$ we have $\EE \exp(\lambda \cdot \xi_i) \leq \exp(\lambda^2\cdot \sigma^2/2)$. We remark that this assumption implies that $\operatorname{Var}(\xi_i) \leq \sigma^2$. Denote the variance $\operatorname{Var}(\xi_i) = \underline \sigma^2$, noting that $\underline \sigma^2 \leq \sigma^2$.

When $\xi \sim N(0,\sigma^2 I_n)$, which is part of our model class, \cite{neykov2022minimax} and later \cite{prasadan2024information} showed that when $K$ is a star-shaped set the minimax rate is given by $(\eta^{\star 2} \wedge d^2)$\footnote{In fact, it is easy to see that $\eta^\star \lesssim d$ and hence the minimax rate can be written as simply $\eta^{\star2}$.}, where $d = \diam(K)$ and $\eta^\star := \eta^\star(K) = \sup \{\eta \geq 0 : \eta^{ 2}/\sigma^2 \leq \log M^{\operatorname{loc}}_K(\eta)\}$ with $M^{\operatorname{loc}}_K(\eta)$ being the local entropy of the set $K$. The goal of the present paper is to construct a polynomial-time algorithm which nearly achieves the minimax rate over a range of convex bodies $K$. Observe that we define the minimax rate with $\sigma$ rather than $\underline \sigma$ (in the Gaussian case which gives the lower bound, we have $\underline \sigma \equiv \sigma$). We will later see (see Section \ref{robust:section}) that this can be relaxed to using $\underline \sigma$ in the definition of $\eta^\star$, but it requires having more than a single sample. 

Throughout the paper, we impose several requirements on the set $K$.  First,
$K$ is symmetric about $0$, so that its Minkowski gauge $\rho_K$ defines a
norm on $\RR^n$.  Second, its type-2 constant $T_2(K)$, defined in
\eqref{type2:constant}, is ``small'' in the sense discussed in Section
\ref{type2:carl:background}; for example, a $\sqrt{\log n}$ factor is small
for our purposes, whereas a $\sqrt n$ factor is not.  Third, we assume an
oracle for computing $\rho_K$ (which implies a weak membership oracle for
$K$; see \cite{grotschel2012geometric}).  Finally, there are known constants
$r,R\in\RR_+$ such that $rB_2\subseteq K\subseteq RB_2$.  The geometric dependence of the running time on these radii is polynomial in $\log(R/r)$; see Remark~\ref{complexity:of:r:and:R:rem} so $r$ and
$R$ may in principle be exponentially small or large in $n$.  We refer to
such a set as well-balanced.  Also, $R=d/2$ always works when the diameter
$d=\diam(K)$ is known, since $K$ is origin symmetric; otherwise necessarily
$R\geq d/2$.

The bridge between entropy and widths is the specialization of Carl's
inequality recorded in \eqref{carls:ineq}.  We now use it to establish a
result guaranteeing the existence of a sufficiently small Kolmogorov width.

\begin{lemma}\label{Kwidth:bound:lemma}
Let $K$ be a type-2 convex body with type-2 constant $T_2(K)$, such that $rB_2 \subseteq K \subseteq RB_2$ for some known $r,R$. Let $\kappa:= \kappa(K) \geq 1$ be a given scalar which can depend on $K$. Let $\eta^\star := \eta^\star(K)$ be the minimax rate defined through the equation $\eta^\star = \sup_\eta\{\eta \geq 0: \eta^2/\sigma^2 \leq \log M_K^{\operatorname{loc}}(\eta)\}$. Assume that $\kappa^{-1/2} \tilde d/T_2(K) \geq 2\eta^\star$ and $\lceil \lceil\log_c(\kappa^{1/2}\cdot T_2(K))\rceil \cdot 4 \eta^{\star 2}/\sigma^2\rceil \leq \tilde d^2/(C^2\sigma^2)$ for $\tilde d = 2R \geq d = \diam(K)$ and where $C > 0$ is some absolute constant. Then for $m = \lceil \tilde d^2/(C^2\sigma^2) \rceil \wedge n$, we have
\begin{align}\label{Kwidth:smaller:than:diameter:up:to:constant}
    d_{m}(K) \lesssim \kappa^{-1/2}\cdot \tilde d/c, 
\end{align}
with $c$ being the constant from the definition of local entropy, which is sufficiently large.
\end{lemma}

\begin{remark}
    Let us remark here that the constant $c$ can be taken to be \textbf{sufficiently large}. Since $\kappa \geq 1$ the statement \eqref{Kwidth:smaller:than:diameter:up:to:constant} is not vacuous and does not follow from the trivial bound $d_m(K) \leq \tilde d$.
\end{remark}

\begin{proof}[Proof of Lemma \ref{Kwidth:bound:lemma}]
    First consider the case when $K$ is a singleton, i.e., $K = \{0\}$ since we are assuming it is origin-symmetric. Then the statement is clearly true since all $d_m(K) \equiv 0$ for all $m \in [n]$. Furthermore if $\sigma = 0$, $m = n$ and thus \eqref{Kwidth:smaller:than:diameter:up:to:constant} continues to hold. Next, consider the non-trivial case when $K$ contains at least two points and $\sigma > 0$. It is easy to see that the former implies that for small $\delta$, $\log M^{\operatorname{loc}}(\delta) \geq \log 2 c > 1$. This in turn implies, $\eta^\star > 0$.
    
    Observe now that for any $\delta \geq 2\eta^*$ we have $\log M^{\operatorname{loc}}(\delta) \leq \log M^{\operatorname{loc}}(2\eta^{\star}) \leq 4\eta^{\star 2}/\sigma^2$ (see \cite{neykov2022minimax} Lemma II.8 where it is shown that for convex sets the map $\delta \mapsto \log M^{\operatorname{loc}}(\delta)$ is non-increasing). Set $\kappa^{-1/2}\cdot \tilde d/T_2(K) = \Upsilon$. Under our assumptions we have that $\Upsilon \geq 2\eta^*$. Next, by Lemma 3 of \cite{yang1999information} (see also Appendix \ref{appendix:yang:barron} where we restate this lemma in our notation) we know that 
\begin{align}
    \log M(\Upsilon c^{j-1}) - \log M(\Upsilon c^{j}) \leq \log M^{\operatorname{loc}}(\Upsilon c^j) \leq 4\eta^{\star2}/\sigma^2,
\end{align}
for $j = 0, 1, \ldots$. Here recall from the definition of local entropy that we use $M(\delta)$ to denote the $\delta$ packing of $K$. Note that for $j \geq \lceil \log_c (\kappa^{1/2}\cdot T_2(K)) \rceil $ we have $\log M(\Upsilon c^{j}) = 0$. Thus upon summing all the inequalities for $j = 0, \ldots, \lceil \log_c (\kappa^{1/2}\cdot T_2(K)) \rceil$ we have
\begin{align*}
        \log M(\Upsilon c^{-1}) \leq \lceil \log_c (\kappa^{1/2}\cdot T_2(K)) \rceil 4 \eta^{\star2}/\sigma^2,
\end{align*}

This implies that for $k = M(\Upsilon c^{-1}) \leq \exp(\lceil \log_c (\kappa^{1/2}\cdot T_2(K)) \rceil 4 \eta^{\star 2}/\sigma^2)$, we have $e_k(K) \leq \kappa^{-1/2} \cdot \tilde d/(T_2(K)c)$. Hence by \eqref{carls:ineq} and since the Kolmogorov widths are decreasing in $m$, we have
\begin{align}\label{super:important:bound}
    d_m(K) \leq (\prod_{k = 1}^m d_k(K))^{1/m}  \lesssim T_2(K) \exp(\lceil \log_c (\kappa^{1/2}\cdot T_2(K)) \rceil 4 \eta^{\star 2}/\sigma^2/m) \kappa^{-1/2} \cdot \tilde d/(T_2(K)c).
\end{align}
Next select $m = \lceil \tilde d^2/(C^2\sigma^2)\rceil \wedge n$, where $\tilde d = 2R$ is an upper bound on $d$. Since both $R$ and $\sigma$ are assumed to be known, this is a legitimate choice of $m$. By assumption $\lceil \lceil \log_c (\kappa^{1/2}\cdot T_2(K)) \rceil \cdot 4 \eta^{\star 2}/\sigma^2\rceil \leq \tilde d^2/(C^2\sigma^2)$, it is clear by \eqref{super:important:bound} that:
\begin{align*}
    d_{m}(K) \lesssim \kappa^{-1/2}\cdot \tilde d/c, 
    \end{align*}
where we note that in case $m = n$ the above inequality is trivial as $d_n(K) = 0$.
\end{proof}

\begin{remark}[Near-optimality of linear projection estimators]\label{remark:linear:estimators:over:type-2:sets} We now digress to provide a useful observation.
Lemma~\ref{Kwidth:bound:lemma}, together with Carl's inequality, also suggests a connection between type-$2$ geometry and the performance of linear projection estimators in the Gaussian sequence model. Indeed, if $P_m$ is an orthogonal projection onto an $m$-dimensional subspace, then the estimator $\widehat\mu=P_mY$ satisfies
\begin{align*}
    \sup_{\mu\in K}
    \EE_\mu\|\widehat\mu-\mu\|_2^2
    =
    \sigma^2 m+
    \sup_{\mu\in K}\|(I-P_m)\mu\|_2^2.
\end{align*}
Consequently, optimizing over all rank-$m$ projections gives
\begin{align*}
    \inf_{P_m}
    \sup_{\mu\in K}
    \EE_\mu\|P_mY-\mu\|_2^2
    =
    \sigma^2m+d_m(K)^2.
\end{align*}
Thus, if one can find
\begin{align*}
    m\lesssim \frac{\eta^{\star2}}{\sigma^2}
    \qquad\text{such that}\qquad
    d_m(K)\lesssim \eta^\star,
\end{align*}
then a linear projection estimator is minimax optimal up to an absolute constant.

Carl's inequality and the argument of Lemma~\ref{Kwidth:bound:lemma} almost yield such a conclusion. In particular, by telescoping the local entropy bounds across scales, one obtains, up to absolute constants,
\begin{align*}
    \log M(K,C\eta^\star)
    \lesssim
    \frac{\eta^{\star2}}{\sigma^2}
    \left(
        1+\log\frac{\diam(K)}{\eta^\star}
    \right).
\end{align*}
Combining this global covering bound with Carl's inequality therefore suggests the existence of an $m$ satisfying
\begin{align*}
    m
    \lesssim
    \frac{\eta^{\star2}}{\sigma^2}
    \left(
        1+\log\frac{\diam(K)}{\eta^\star}
    \right),
    \qquad
    d_m(K)\lesssim T_2(K)\eta^\star.
\end{align*}
Hence the best linear projection estimator should satisfy a risk bound of the form
\begin{align*}
    \inf_m\left\{
        \sigma^2m+d_m(K)^2
    \right\}
    \lesssim
    \left(
        T_2(K)^2+
        1+\log\frac{\diam(K)}{\eta^\star}
    \right)\eta^{\star2}.
\end{align*}
We note that the factor $\log(\diam(K)/\eta^\star)$ need not, in general, be logarithmic in the ambient dimension, since the diameter $d=\operatorname{diam}(K)$ may be arbitrarily large relative to the noise level. In the Gaussian sequence model, however, this dependence on the global diameter can be removed by a sample-splitting localization argument: using one independent Gaussian copy to construct a pilot estimator (namely the Euclidean projection of the copy onto $K$), the remaining estimation problem can be localized to the symmetric body $2K\cap C\sqrt{n}\sigma B_2$,
so that the corresponding logarithmic factor is replaced by
$\log\left(\frac{C\sqrt{n}\sigma}{\eta^\star}\right) \lesssim \log n.$
We emphasize that this produces a sample-split, adaptive estimator rather than a single fixed linear estimator. Details of this construction are given in Appendix \ref{appendix:linear:type2}.
It is natural to ask whether the logarithmic loss above is necessary. More specifically, we conjecture that for type-$2$ bodies one may have
\begin{align*}
    \inf_m\left\{
        \sigma^2m+d_m(K)^2
    \right\}
    \lesssim
    \operatorname{poly}(T_2(K))\,\eta^{\star2}.
\end{align*}
Such a result would imply that a suitably chosen linear projection estimator is minimax optimal, up to a factor depending only on the type-$2$ constant. The main difficulty is that the minimax rate is governed by \emph{local} metric entropy, whereas Carl's inequality relates Kolmogorov widths to \emph{global} metric entropy; it is presently unclear whether the type-$2$ assumption alone is sufficient to bridge this gap without an additional multiscale logarithmic loss. In addition, it is worth emphasizing that this observation does not by itself yield a polynomial-time algorithm for computing the best linear projection estimator, since one must still determine the appropriate dimension $m$ and compute an optimal Kolmogorov-width projection $P_m$, neither of which is obviously computationally tractable.
\end{remark}

Unfortunately, we do not know the optimal projection that yields a small Kolmogorov width. In addition, even calculating the Kolmogorov width $d_m(K)$ precisely in general is likely NP hard -- see, for instance, \cite{brieden2002geometric} where the author shows that calculating $d_0(K)$ for (even symmetric) convex polytopes cannot be computed in polynomial-time within a factor of $1.090$, unless P=NP. Thus, the fact \eqref{Kwidth:smaller:than:diameter:up:to:constant} appears useless at first glance. However, since we are interested in having a minimax optimal estimator up to ``small'' factors, we do not need to find the precise optimal Kolmogorov projection. Taking inspiration from \cite{varadarajan2007approximating}, we propose to solve the following semi-definite program (SDP) relaxation of the Kolmogorov width problem
\begin{align}
    &\argmin_{X} \max_{p \in K}  p\T X p ,\label{quadratic:maximization}\\
    &\mbox{subject to }\operatorname{tr}(X) = n-m,\nonumber\\
    &~~~~~~~~~~~~~~~0 \leq X \leq I, X\T = X \nonumber
\end{align}

As observed in \cite{varadarajan2007approximating}, the optimal value $X^{\star}$ that this SDP returns, satisfies $p\T X^{\star} p \leq d^2_{m}(K)$ for all $p \in K$. The proof of this fact is elementary, and we include a proof below for convenience of the reader. The proof also helps motivate program \eqref{quadratic:maximization}.

\begin{lemma}\label{SDP:relaxation:bound}
    The solution to the above SDP problem, $X^{\star}$, satisfies: $p\T X^{\star} p \leq d^2_{m}(K)$ for all $p \in K$.
\end{lemma}
\begin{proof}[Proof of Lemma \ref{SDP:relaxation:bound}] Let $P_{n-m} = I - P_m$, where $P_m$ is the optimal Kolmogorov width projection, i.e., 
\begin{align*}
    \sup_{p \in K} \|(I - P_m)p\|_2^2 = d^2_m(K).
\end{align*} 
This is equivalent to $\sup_{p \in K} p\T P_{n-m}p \le d^2_m(K)$. Note that $P_{n-m}$ is a feasible point to the SDP hence the result follows.
\end{proof}

Notably, solving the SDP requires maximizing a quadratic form over a convex body (cf. \eqref{quadratic:maximization}), a formidable task on its own. The authors of \cite{varadarajan2007approximating} do not encounter this issue, since their body 
$K$ is a convex polytope with finitely many vertices, making the quadratic maximization tractable (the maximum is attained at a vertex and can be found by enumerating all vertices). In our setting, no such guarantee holds. We will now assume that there exists an oracle $\cO_K$ which takes a symmetric psd matrix $X$ and outputs a vector $\cO_K(X) \in K$ such that 
\begin{align}\label{oracle:equation}
    \cO_K(X)\T X \cO_K(X) \geq \kappa^{-1}\max_{p \in K} p\T X p,
\end{align}
for some ``small''\footnote{Here by small we mean that $\kappa(K)$ is a constant slowly growing with $n$.} constant $\kappa = \kappa(K) \geq 1$ with probability at least $1-q$ for some $q \in (0,1)$ with $\log q^{-1}$ being a parameter involved in the number of iterations needed to compute the oracle, whose value will be specified later (see Section \ref{QFM:section}). In other words, we suppose there exists an approximate solution to the quadratic maximization problem. Later we will see that there is a rich class of type-2 bodies for which such an oracle indeed exists. In order to state our algorithm for approximately solving the SDP, we consider the projection of a vector $v \in \mathbb{R}^n$ onto the convex set
\begin{align*}
S_k = \{\, w \in \mathbb{R}^n : 0 \le w_i \le 1,\ \sum_{i=1}^n w_i = k \,\},
\end{align*}
where $k \in \{0,\ldots, n\}$ is given. The projection is defined as
\begin{align*}
\Pi_{S_k}(v) = \arg\min_{w \in S_k} \frac{1}{2}\|w - v\|_2^2.
\end{align*}
$\Pi_{S_k}(v)$ can be implemented efficiently (see Algorithm \ref{alg:capped_simplex}). Given the oracle $\cO_K$ and the projection $\Pi_{S_k}(\cdot)$, we may approximately solve the SDP with projected sub-gradient descent. We formalize the algorithm below:

\begin{breakablealgorithm}
    \renewcommand{\algorithmicrequire}{\textbf{Input:}}
    \renewcommand{\algorithmicensure}{\textbf{Output:}}
    \caption{SDP with projected sub-gradient descent}
    \label{algorithm:SDP:Kolmogorov:width}
   \begin{algorithmic}[1]
        \REQUIRE Convex body $K$, value $m \in \{0,\ldots, n\}$, oracle $\cO_K$ for solving quadratic maximization over $K$, a tuning parameter $C \in \RR_+$, $M = 1 \vee \lceil 4(n-m)C^4\kappa^2\rceil$ maximal number of iterations, $\gamma = 1/(\kappa C^2) \in \RR$ is the learning rate.
        \ENSURE A symmetric psd matrix $X^{\star\star}$

        \STATE Set $X^{(1)} = \diag(((n-m)/n, (n-m)/n, \ldots, (n-m)/n))$
        \IF{$m == n$}
            \RETURN $X^{(1)}$
        \ENDIF
        \FOR{ k= 1, \ldots, M}
            \STATE Evaluate $p_k = \cO_K(X^{(k)})$
            \STATE Set $\gamma_k = \gamma/\|p_k\|_2^2$
            \STATE $X^{(k+1)} = X^{(k)} - \gamma_k p_k p_k\T$
            \STATE Eigendecompose $X^{(k + 1)} = U \Lambda U\T$
            \STATE $v = \Pi_{S_{n-m}}(\operatorname{diag}(\Lambda))$
            \STATE $ X^{(k + 1)} = U \operatorname{diag}(v) U\T$
        \ENDFOR
        \RETURN $X^{\star\star} := \frac{\sum_{k = 1}^{M }\gamma_k X^{(k)}}{\sum_{k = 1}^{M}\gamma_k }$
    \end{algorithmic} 
\end{breakablealgorithm}

Next we write down the sub-routine for $\Pi_{S_{n-m}}$ used in Algorithm \ref{algorithm:SDP:Kolmogorov:width} as Algorithm \ref{alg:capped_simplex}. We give the justification of Algorithms \ref{algorithm:SDP:Kolmogorov:width} and \ref{alg:capped_simplex} below.

\begin{breakablealgorithm}
    \renewcommand{\algorithmicrequire}{\textbf{Input:}}
    \renewcommand{\algorithmicensure}{\textbf{Output:}}
\caption{Projection onto the capped simplex 
$S_k = \{ w  : 0 \le w_i \le 1, \sum_i w_i = k \}$}
\label{alg:capped_simplex}
\begin{algorithmic}[1]
\REQUIRE Vector $v \in \mathbb{R}^n$, integer $k \in \{0,\ldots, n\}$
\ENSURE Projected vector $w^\star = \arg\min_{w \in S_k} \|w - v\|_2^2$
\STATE Define breakpoints $\mathcal{B} = \{v_i - 1, v_i : i=1,\ldots,n\} \cup \{+\infty\}$
\STATE Sort the breakpoints in increasing order
\FOR{each interval $[b_j, b_{j+1})$ between consecutive breakpoints}
    \STATE Let $I_0 = \{i : v_i - b_j \ge 1\}$, $I_1 = \{i : 0 < v_i - b_j < 1\}$
        \IF{$|I_1| == 0$}
        \IF {$|I_0| == k$}
            \STATE Set $\theta = b_j$
            \STATE \textbf{break}
        \ENDIF
        \STATE \textbf{next}
    \ENDIF
    \STATE Compute
    \begin{align*}
        \theta = \frac{\sum_{i \in I_1} v_i + |I_0| - k}{|I_1|}
    \end{align*}
    \IF{ $\theta \in [b_j, b_{j+1})$}
        \STATE \textbf{break}
    \ENDIF
\ENDFOR
\STATE Return $w_i^\star = \min\{1, \max\{0, v_i - \theta\}\}$ for all $i$
\end{algorithmic}
\end{breakablealgorithm}

We start by justifying Algorithm \ref{alg:capped_simplex} which finds $\Pi_{S_k} v$. Specifically, we have:

\begin{lemma}[Projection onto the capped simplex]
\label{thm:capped-simplex}
Let $v \in \mathbb{R}^n$ and $k \in \{0,\ldots, n\}$.  
Then the unique minimizer of
\begin{align*}
\min_{w \in S_k} \frac{1}{2}\|w - v\|_2^2
\end{align*}
is given componentwise by
\begin{align*}
w_i^\star = \min\{1, \max\{ 0, v_i - \theta^\star \}\},
\end{align*}
where the scalar $\theta^\star$ is a solution of
\begin{align*}
\sum_{i=1}^n \min\{1, \max\{ 0, v_i - \theta^\star \}\} = k.
\end{align*}
In addition, observe that $w_i$ are ordered in the same way, in other words $(v_i - v_j) (w_i - w_j) \geq 0$ for any $i,j$.
\end{lemma}

The proof of Lemma \ref{thm:capped-simplex} is deferred to the appendix. It is clear that this justifies Algorithm \ref{alg:capped_simplex}. We now argue that the projection of a matrix, on the set of symmetric psd matrices which are smaller than the identity matrix in the psd sense and have fixed trace, is given by a simple spectral decomposition plus a projection on the capped simplex.

\begin{lemma} Let $X$ be a symmetric matrix. Then the projection of $X$ onto the set $\{W \text{ symmetric psd} : 0 \leq W \leq I, \operatorname{tr}(W) = k\}$ in the Frobenius norm (i.e., $W^\star = \argmin_W \|X - W\|_F^2$) is given by
\begin{align*}
    W^\star = U \operatorname{diag}(\Pi_{S_k}(\operatorname{diag}(\Lambda))) U\T, \text{ where } X = U \Lambda U\T,
\end{align*}
is the spectral decomposition of $X$. 
\end{lemma}

\begin{proof} Let $W^{\star}$ be the minimizer of $\|X - W^{\star}\|_F^2$ for $W^\star \in \{W \text{ symmetric psd} : 0 \leq W \leq I, \operatorname{tr}(W) = k\}$. Suppose $W^\star = U^\star \Lambda^\star U^{\star\top}$, where $\Lambda^\star$ is a diagonal matrix sorted in descending way (same as sorting as the diagonal entries of $\Lambda$). Let $\lambda = \operatorname{diag}(\Lambda)$ and $\lambda^\star = \operatorname{diag}(\Lambda^\star)$. We have

\begin{align*}
    \|X - W^{\star}\|_F^2 = \|X\|_F^2 - 2 \operatorname{tr} (X\T W) + \|W\|_F^2 = \|\lambda\|_2^2 - 2 \operatorname{tr} (X\T W)  +  \|\lambda^\star\|_2^2.
\end{align*}
Note that 
\begin{align*}
    \operatorname{tr} (X\T W) = \sum_{i} \sum_{j} (u_i\T u_j^\star)^2 \lambda_i \lambda^\star_j = \lambda\T A \lambda^\star,
\end{align*}
where $A = \{A_{ij}\}_{i,j \in [n]} = (u_i\T u_j^\star)^2$ is a doubly stochastic matrix (by the orthonormality of $U$ and $U^\star$). By Birkhoff–von Neumann's theorem $A = \sum_{i \in [L]}\alpha_i P_i$ for some permutation matrices $P_i$, $L \in \NN$ and $\sum_{i \in [L]} \alpha_i = 1$. Thus

\begin{align*}
    \operatorname{tr} (X\T W) = \sum_{i \in [L]}\alpha_i \lambda\T P_i \lambda^\star \leq \lambda\T \lambda^\star,
\end{align*}
where the last inequality comes from the rearrangement inequality since both $\lambda$ and $\lambda^\star$ are vectors sorted in descending order. 

We conclude that 
\begin{align*}
    \|X - W^{\star}\|_F^2 \geq \|\lambda\|_2^2 - 2 \lambda\T \lambda^\star + \|\lambda^\star\|_2^2 = \|\lambda - \lambda^\star\|_2^2.
\end{align*}

The above inequality is also clearly an equality whenever $U^\star = U$. Since the projection of $\lambda$ onto the set $S_k$ is $\Pi_{S_k}(\lambda)$ as we showed in Lemma \ref{thm:capped-simplex} the result follows. 

\end{proof}

The following theorem establishes the validity of Algorithm~\ref{algorithm:SDP:Kolmogorov:width}.

\begin{theorem} \label{main:SDP:theorem} Given as input a convex body $K$, a value $m$ for which $d_m(K) \leq L \kappa^{-1/2} \tilde d/c$ for a constant $L > 0$, and an oracle $\cO_K$ satisfying \eqref{oracle:equation} with probability at least $1 - q$ with $q = \frac{M^{-1}}{1 + \tilde d^3/\sigma^3}$, Algorithm \ref{algorithm:SDP:Kolmogorov:width} outputs a matrix $X^{\star\star}$ which, up to an absolute constant factor, satisfies
\begin{align*}
    \sup_{p \in K} p\T X^{\star\star} p \lesssim \tilde d^2/c^2, 
\end{align*}
with probability at least $1 - (1 + \tilde d^3/\sigma^3)^{-1}$.
\end{theorem}

\begin{remark}
    Let us remark that similarly to Lemma \ref{Kwidth:bound:lemma} the above statement is not vacuous, and contracts the diameter $\tilde d$, provided that the constant $c$ is taken to be sufficiently large and $C$ from the definition of $\gamma$ in Algorithm \ref{algorithm:SDP:Kolmogorov:width} is also sufficiently large. 
\end{remark}

\begin{remark}\label{poly:time:oracle:remark}
    Observe that since $M = 1 + \lceil 4(n-m)C^4\kappa^2\rceil$ and the oracle works with polynomially many iterations in $\log q^{-1}$ (see Section \ref{QFM:section}), if $\kappa$ is a ``small'' constant (i.e., $\kappa \ll n$), and in addition $\log \tilde d/\sigma$ is polynomial in $n$, we are guaranteed to have $\log [M (1 + \tilde d^3/\sigma^3)]$ is polynomial in $n$, so that our oracle algorithm will run in polynomial time.
\end{remark}

\begin{proof}[Proof of Theorem \ref{main:SDP:theorem}] We first establish the high probability bound. Since we have $M$ iterations, and the oracle succeeds with probability $1-M^{-1}/(1 + \tilde d^3/\sigma^3)$ every time, by the union bound it will fail with probability smaller than $1/(1 + \tilde d^3/\sigma^3)$, which establishes the claim. Let us now show the guarantee assuming that the oracle never fails.

Let us begin by observing that the function $f(X) := \sup_{p \in K} p\T X p = \sup_{p \in K} \langle X, pp\T \rangle$ is convex as a maximum of linear functions. It is also well known that $p^\star p^{\star \top}$ is a sub-gradient of this function, where $p^\star = \argmax_{p \in K} \langle X, pp\T \rangle$. Next, note that by our definition of the oracle, the vector $\cO_K(X)$ is near sub-gradient of the map in the following sense:
\begin{align}\label{near:subgradient:property}
    f(X) = \sup_{p \in K} \langle X, pp\T \rangle \leq \langle X, \cO_K(X)\cO_K(X)\T \rangle + (1-\kappa^{-1})\cdot f(X)
\end{align}

Put $g^{(k)} = p_k p_k\T$. Recall that $X^{\star}$ is the true optimum, which satisfies by Lemma \ref{SDP:relaxation:bound} that $\sup_{p \in K}p\T X^{\star} p \leq d^2_{m}(K)$ for all $p \in K$. Observe that we have 
\begin{align*}
    f(X^{(k)}) - f(X^{\star}) \leq \langle g^{(k)}, X^{(k)} - X^\star\rangle + (1-\kappa^{-1})\cdot f(X^{(k)}),
\end{align*}
where we used \eqref{near:subgradient:property} and the fact that $\langle g^{(k)}, X^\star \rangle \leq f(X^\star)$. We will now repeat the original proof of the sub-gradient descent due to Shor \cite{shor2012minimization, boyd2004convex}, using the near sub-gradient property above. Since Euclidean projection on a convex set shrinks the distance to any point in the set, we have
\begin{align*}
    \|X^{(k + 1)} - X^\star\|_F^2 & \leq \|X^{(k)} - X^\star\|_F^2 + 2\gamma_k \langle g^{(k)}, X^\star - X^{(k)} \rangle + \gamma_k^2\|g^{(k)}\|_F^2\\
    &\leq \|X^{(k)} - X^\star\|_F^2 + 2\gamma_k(f(X^\star) - \kappa^{-1}\cdot f(X^{(k)})) + \gamma_k^2\|g^{(k)}\|_F^2.
\end{align*}

Thus after $M$ iterations we have
\begin{align}\label{eq:weighted-subgradient-sum}
    2\sum_{k= 1}^M \gamma_k(\kappa^{-1}\cdot f(X^{(k)}) - f(X^\star))  & \leq \|X^{(1)} - X^\star\|_F^2 + \sum_{k = 1}^M\gamma_k^2\|g^{(k)}\|_F^2.
\end{align}
Since $f$ is convex
\begin{align}
f(X^{\star\star})
&=
f\left(
\frac{\sum_{k=1}^M\gamma_kX^{(k)}}
{\sum_{k=1}^M\gamma_k}
\right)\nonumber\\
&\leq
\frac{
\sum_{k=1}^M\gamma_k f(X^{(k)})
}{
\sum_{k=1}^M\gamma_k
}.
\label{eq:convexity-f-average}
\end{align}
Combining \eqref{eq:weighted-subgradient-sum} and
\eqref{eq:convexity-f-average}, we obtain
\begin{align}
\kappa^{-1}f(X^{\star\star})-f(X^\star)
\leq
\frac{
\|X^{(1)}-X^\star\|_F^2
+
\sum_{k=1}^M\gamma_k^2\|g^{(k)}\|_F^2
}{
2\sum_{k=1}^M\gamma_k
}.
\label{eq:average-SDP-bound}
\end{align}

Next we note that $\|g^{(k)}\|_F \leq d_0(K)^2 \leq d^2$ for any $k$ (since $\cO_K(X^{(k)}) \in K$). Furthermore, $\|X^{(1)}\|_F^2 = \sum_{i \in [n]}\lambda^2_i(X^{(1)}) \leq \sum_{i \in [n]}\lambda_i(X^{(1)}) \leq n-m$, and a similar logic shows that $\|X^\star\|_F^2 \leq (n-m)$. Hence $\|X^{(1)} - X^\star\|_F^2 \leq 4 (n-m)$. Hence upon selecting, $\gamma_k = \gamma/\|g^{(k)}\|_F$ (noting that $\|g^{(k)}\|_F = \|p_k p_k\T\|_F = \|p_k\|_2^2$) we obtain
\begin{align*}
    \kappa^{-1}\cdot f(X^{\star\star}) - f(X^\star) \leq \frac{4(n-m) + M\gamma^2}{2M\gamma/d^2 }.
\end{align*}
Therefore taking into account that $\gamma = 1/(\kappa\cdot C^2)$, it suffices to take $M \geq 4(n-m)C^4\kappa^2$. We thus conclude that  $f(X^{\star \star})$ satisfies
\begin{align*}
    f(X^{\star \star}) &\leq \kappa\cdot(\kappa^{-1}\cdot\tilde d^2/C^2 + f(X^\star)) \leq \kappa\cdot(\kappa^{-1} \cdot \tilde d^2/C^2 + d_m(K)^2) \leq \kappa\cdot(\kappa^{-1} \cdot \tilde d^2/C^2 + L^2\cdot \kappa^{-1}\cdot \tilde d^2/c^2)\\
    &\lesssim \tilde d^2/c^2,
\end{align*}
where we used $\tilde d\geq d$, which is what we aimed to show.

\end{proof}

Once we have obtained the matrix $X^{\star \star} = V\Lambda V\T = \sum_{i \in [n]} \lambda_i v_i v_i\T$, our algorithm proceeds with constructing the matrix $A^{\star\star} = (I - X^{\star\star})^{1/2} = \sum_{i \in [n]}\sqrt{1-\lambda_i}v_iv_i\T$. Consider now the vector $A^{\star\star}Y$. It is equal to
\begin{align*}
    A^{\star\star}Y =  A^{\star\star}\mu + A^{\star\star}\xi.
\end{align*}
In order to estimate $A^{\star\star}\mu$, we will use $A^{\star\star}Y$. We have the following lemma justifying that this is not a bad choice.

\begin{lemma}\label{sub:gaussian:vector:estimator:1}
     Suppose $\tilde d \gtrsim \sigma$. Then, the estimator $A^{\star\star}Y$ of $A^{\star\star}\mu$ satisfies
    \begin{align*}
        \|A^{\star\star}Y - A^{\star\star}\mu \|_2^2 \leq \tilde d^2/T^2,
    \end{align*}
    with probability at least $1 - \exp(-\tilde d^2/(\sigma^2C^{'2}))$.
    \end{lemma}

\begin{proof}

Applying the Hanson-Wright inequality \cite{rudelson2013hanson}, we obtain:
\begin{align}\label{HW:concentration:bound}
    \PP(\|A^{\star\star}\xi\|_2^2 - \underline \sigma^2 \operatorname{tr}(A^{\star\star 2})  \geq t) \leq \exp(-c t^2/(\sigma^4\|A^{\star\star2}\|^2_F) \wedge t/(\sigma^2\|A^{\star\star2}\|_{\operatorname{op}}))
\end{align}
     Note that $\|A^{\star\star2}\|_{\operatorname{op}} \leq 1$.  In addition $\|A^{\star\star2}\|^2_F = \operatorname{tr}((I- X^{\star\star})^2)\leq \operatorname{tr}(I- X^{\star\star}) = m$.

    Hence we can select $t = \tilde d^2/c^2 \gtrsim m \sigma^2$ to claim that \eqref{HW:concentration:bound} holds with high probability, i.e., probability at least $1 - \exp(-c' \tilde d^4/(\sigma^4 m)\wedge \tilde d^2/(c^2\sigma^2))$. Comparing $\tilde d^4/(\sigma^4 m) \gtrsim \tilde d^2/\sigma^2$ yields $\tilde d^2 \gtrsim \sigma^2 m$ which is true as $m = \lceil \tilde d^2/(C^2\sigma^2) \rceil \wedge n$ and $\tilde d \gtrsim \sigma$. This implies $c'\tilde d^4/(\sigma^4 m)\wedge \tilde d^2/(c^2\sigma^2) \asymp \tilde d^2/(C^{'2}\sigma^2)$ for some $C' > 0$. So that the bound is always $1 - \exp(-\tilde d^2/(C^{'2}\sigma^2))$. Since the mean of $\|A^{\star\star}\xi\|_2^2$ is $ \underline \sigma^2 \operatorname{tr}(A^{\star\star 2}) = \underline \sigma^2 m \lesssim \sigma^2 \tilde d^2/(C^2\sigma^2)$ we complete the proof.
\end{proof}

\begin{remark} Lemma \ref{sub:gaussian:vector:estimator:1} is the reason why we assumed that the noise $\xi$ has independent components, so that we can apply the Hanson-Wright inequality. We note that the formula for the minimax rate remains unchanged even if the noise has dependent entries \cite{prasadan2024information}. In Section \ref{robust:section} we will see that if we have more samples the noise can be arbitrary provided that its covariance matrix has a bounded operator norm.
\end{remark}

Our main algorithm starts at $0$ and works by iteratively refining the estimate arguing that at each step we have better and better approximation of $\mu$. These approximations will be shown to decay exponentially as the iteration proceeds, and in that sense the following lemma is very useful. Assuming that we have ''trapped'' $\mu$ in an Euclidean ball of certain radius we demonstrate how we can find an even better estimator of $\mu$ provided that the conditions of Lemma \ref{Kwidth:bound:lemma} are satisfied.

\begin{lemma}\label{sub:gaussian:vector:estimator} Suppose we have an estimator $\hat \mu_j \in K$ of $\mu$ such that $\|\hat \mu_j - \mu\|_2 \leq \tilde d_j$ happens with high probability, for some $\tilde d_j \in \RR_+$ such that $\tilde d_j \geq 2r \vee C\sigma$. Let $K_{(j)} = K \cap B(0, \tilde d_j/2)$. Obtain $m$ and the matrix $A_j^{\star\star}$ from Lemma \ref{Kwidth:bound:lemma} and Algorithm \ref{algorithm:SDP:Kolmogorov:width} for the set $K_{(j)}$\footnote{We explain later that the oracle $\cO_{K_{(j)}}$ can also be run in polynomial-time, provided that the oracle $\cO_K$ can be run in polynomial-time.}, assuming that the conditions of Lemma \ref{Kwidth:bound:lemma} hold. Then the estimator $\tilde \mu_{j + 1} = (A^{\star\star}_j Y - A_j^{\star\star} \hat \mu_j)/2$ of $(\mu - \hat \mu_j)/2$ 
used in Lemma \ref{sub:gaussian:vector:estimator:1} satisfies:
\begin{align*}
    \PP(\|\tilde \mu_{j + 1} -( \mu - \hat \mu_j)/2\|_2 \geq \tilde d_j / L, \|\hat \mu_j - \mu\|_2 \leq \tilde d_j) \leq \exp(-\tilde d_j^2/(C^{'2} \sigma^2)) + \frac{1}{1 + \tilde d_j^3/\sigma^3}.
\end{align*}
for some absolute constant $L$.
    
\end{lemma}

\begin{proof}
     Let us first condition on the event that our oracle did not fail (which happens with probability at least $1 - (1 + \tilde d_j^3/\sigma^3)^{-1}$). On the event $\|\hat \mu_j - \mu\|_2 \leq \tilde d_j$ we have that, $\nu = (\mu - \hat \mu_j)/2 \in K_{(j)} = K \cap B(0, \tilde d_j/2)$. 
Recall that $A_j^{\star\star} = (I - X_j^{\star\star})^{1/2}$, and let $X_j^{\star\star} = V \Lambda V\T $ be the eigendecomposition of $X_j^{\star\star}$.
\begin{align*}
    \|\nu - A_j^{\star\star}\nu\|_2^2 &= \|\nu - \sum_i \sqrt{1-\lambda_i} v_i (v_i\T \nu)\|^2_2 = \sum_i (1- \sqrt{1-\lambda_i})^2 (v_i\T \nu)^2 \\&= \sum_i \frac{\lambda_i^2}{(1+ \sqrt{1-\lambda_i})^2} (v_i\T \nu)^2 
    \leq \sum_i \lambda_i (v_i\T \nu)^2 =  \nu\T X_j^{\star \star} \nu \lesssim (\tilde d_j/c)^2,
\end{align*}
where we used $0 \leq \lambda_i \leq 1$. Next, $\|\tilde \mu_{j + 1} -( \mu - \hat \mu_j)/2\|_2 \leq \|A_j^{\star\star} Y/2 - A_{j}^{\star\star}\mu/2\|_2 + \|A_j^{\star\star}(\mu - \hat \mu_j)/2 -( \mu - \hat \mu_j)/2 \|_2 \lesssim \tilde d_j/(2T) + (\tilde d_j/c) \leq \tilde d_j/L$, with probability at least $1 - \exp(-\tilde d_j^2/(C^{'2} \sigma^2)) - (1 + \tilde d_j^3/\sigma^3)^{-1}$ by Lemma \ref{sub:gaussian:vector:estimator:1}. This completes the proof.
\end{proof}

Finally we can describe the iteration of our algorithm. 

\begin{breakablealgorithm}
    \renewcommand{\algorithmicrequire}{\textbf{Input:}}
    \renewcommand{\algorithmicensure}{\textbf{Output:}}
    \caption{Near-optimal Estimation Algorithm}
    \label{algorithm:main}
   \begin{algorithmic}[1]
        \REQUIRE Convex body $K$, $r, R \in \RR_+$ such that $rB_2 \subseteq K \subseteq R B_2$, observation $Y = \mu + \xi$, constant $\tilde L = L/(\sqrt{3} + 1)$ where $L$ is from Lemma \ref{sub:gaussian:vector:estimator}; We assume $\tilde L > 2(\sqrt{3} + 1)$\footnote{This assumption can be met upon selecting $C$ from the definition of $m$ and $c$ from the definition of local entropy sufficiently large.}
        \ENSURE Estimator $\hat \mu$ of $\mu$ 
        
        \STATE Set $\tilde d_1 = 2R$, $K_{(1)} := K$, number of iterations $M = \lceil\log_{\tilde L/(2(\sqrt{3} + 1))} R/r\rceil$
                \STATE Set $\hat \mu_1 := 0$
        \IF{$\sigma \leq r/\sqrt{n}$}
            \RETURN $Y$;
        \ENDIF

        \FOR{ j = 1, \ldots, M}
            \STATE Find $m = \lceil \tilde d_j^2/(C^2\sigma^2)\rceil \wedge n$, and calculate the matrix $A_j^{\star \star}$ for the set $K_{(j)}$
            \STATE Find $\tilde \mu_{j + 1} = (A_j^{\star \star}Y - A_j^{\star \star}\hat \mu_j)/2$
                        \STATE Set $\bar \mu_{j + 1} = \Pi_{K_{(j)}} \tilde \mu_{j + 1}$\footnote{Here we set $\epsilon = \tilde d_j/L$ in Lemma \ref{weak:projection:lemma} when evaluating the weak projection.}
            \STATE Lift back to $K$: $\hat \mu_{j + 1} = \Pi_K(2 \bar \mu_{j + 1} + \hat \mu_j)$\footnote{Here we set $\epsilon = 2\tilde d_j/\tilde L$ in Lemma \ref{weak:projection:lemma} when evaluating the weak projection.}
            \STATE Set $\tilde d_{j + 1} = 2(\sqrt{3} + 1)\tilde d_j/\tilde L$
            \STATE Set $K_{(j + 1)} := K\cap B_2(0, \tilde d_{j+1}/2)$
            \IF{$\tilde d_{j + 1} \leq (2r) \vee (C \sigma)$}
                \RETURN $\hat \mu := \hat \mu_{j + 1}$
            \ENDIF
        \ENDFOR
        \RETURN $\hat \mu := \hat \mu_{M + 1}$
    \end{algorithmic} 
\end{breakablealgorithm}

We make several comments about the algorithm. Our first observation is that Algorithm \ref{algorithm:main} does not require knowledge of $\eta^\star$, which may be difficult to evaluate precisely. Instead, the algorithm continues contracting for as long as further contraction is possible, thereby progressively reducing the distance between the current iterate and the signal $\mu$. The various reasons for which this contraction procedure may terminate are addressed by the preceding lemmas, which together guarantee that the resulting estimator is near-minimax optimal in expectation.

Note that, since we can evaluate the Minkowski gauge of $K$, we have a (weak) separation oracle for $K$ \cite{grotschel2012geometric}, which in turn enables (weak) projection in polynomial-time (see Theorem 2.5.9 in \cite{dadush2012integer}, see also \cite{lee2018efficient}). The properties of this weak projection are summarized in Lemma \ref{weak:projection:lemma}, where it is observed that it can output a vector that is in $K$ which is closer (up to an absolute constant factor) to any vector in the set $K$ than the original vector. When we write $\Pi_{K}$ we mean using this weak projection. In addition, for any $j$ the set $K_{(j)}$ also admits weak projection. This is because the Minkowski gauge of $K_{(j)}$ satisfies 
\begin{align}\label{intersection:minkowski:gauge}
    \rho_{K_{(j)}}(x) = \rho_K(x) \vee (2/\tilde d_j \cdot \|x\|_2),
\end{align}
and can therefore be explicitly evaluated. In addition $K_{(j)}$ is well balanced since $r B_2 \subseteq K_{(j)} \subseteq \tilde d_j/2 B_2$. Thus when we write $\Pi_{K_{(j)}}$ we mean the weak projection on $K_{(j)}$. Since the projections are ``weak'' they need not shrink distance to any point in $K$. Hence the need for scaling factors: $(\sqrt{3} + 1)$ and $2(\sqrt{3} + 1)$ in our algorithm.

We remark that our estimator is not proper, because when $r > \sigma\sqrt{n}$ the estimator $Y$ need not be in $K$. This can be remedied by outputting $\Pi_K Y$ in the case when $r > \sigma\sqrt{n}$, with $\epsilon = \sigma\sqrt{n}$ in Lemma \ref{weak:projection:lemma} or setting $r =  r \wedge \sqrt{n}\sigma$ which circumvents this issue altogether.

Next, to ensure that the matrix $A_j^{\star\star}$ is calculable, we need to show the existence of the oracle $\cO_{K_{(j)}}$ for any $j = 1,\ldots, M$. We will do so in the following section, Section \ref{QFM:section}, but for now we will show that if the original set $K$ is type-2, the sets $K_{(j)}$ remain type-2. This is so since its Minkowski gauge satisfies \eqref{intersection:minkowski:gauge}, and is thus type-2 with constant at most 
\begin{align*}
    T_2(K_{(j)}) \leq \sqrt{2}\max(T_2(K), T_2(B_2(0, \tilde d_j/2))) = \sqrt{2}\max(T_2(K), 1) = \sqrt{2} T_2(K).       
\end{align*}
Furthermore, as $K_{(j)} \subseteq K$, it is clear that the local entropy $\log M^{\operatorname{loc}}_{K_{(j)}} (\delta) \leq \log M^{\operatorname{loc}}_{K} (\delta)$ for any $\delta$. What is more, the latter implies that we have $\eta^\star(K_{(j)}) \leq \eta^{\star}(K)$. 

Before we go on to show that our algorithm is nearly minimax optimal we first show the following lemma, justifying why one can output $Y$ when $\sigma$ is small. 
\begin{lemma}\label{lemma:r:sigma}
    Let $\sigma \gtrsim r/\sqrt{n}$. Then the minimax rate $\eta^\star \gtrsim r$. On the other hand if $\sigma \lesssim r/\sqrt{n}$ we have $\eta^{\star} \asymp \sqrt{n}\sigma$.
\end{lemma}

\begin{proof}
    Recall that the minimax rate $\eta^\star$ satisfies $\eta^\star = \sup \{\eta \geq 0 : \eta^{ 2}/\sigma^2 \leq \log M^{\operatorname{loc}}(\eta)\}$. Observe that for $\eta \lesssim r$ we have $\log M^{\operatorname{loc}}(\eta) \asymp n$. On the other hand since $\sigma \gtrsim r/\sqrt{n}$, the ratio $(\eta/\sigma)^2 \lesssim n$. Thus $\eta^\star \geq \eta$ by definition. Since $\eta$ can be taken as $\eta \asymp r$ the proof follows.

    For the second part observe that since $\sigma \lesssim r/\sqrt{n}$ then $\sqrt{n}\sigma \lesssim r$ and therefore $\log M^{\operatorname{loc}}(\sqrt{n}\sigma) \asymp n$. Thus for $\eta \asymp \sqrt{n}\sigma$ it follows that $\eta^2/\sigma^2 \asymp \log M^{\operatorname{loc}}(\eta)$. On the other hand, $\eta^{\star}$ cannot be larger than $\sqrt{n} \sigma$ since otherwise $\eta^{\star2}/\sigma^2 \gtrsim n \gtrsim \log M^{\operatorname{loc}}(\eta^\star)$.
\end{proof}

We now state and prove the main result of the paper:

\begin{theorem}\label{main:result:GSM}
    Set $\overline \kappa := \max_{k \in [M]}\kappa(K_{(k)})$ and suppose $\log(R/r)$ is polynomial in $n$, while $\overline \kappa$ and $T_2(K)$ are polylogarithmic in $n$. Define $\Upsilon(K) := (1 \vee\left(2\sqrt{2}\,T_2(K)\overline\kappa^{1/2}\right)\vee\left(2C\sqrt{\Lambda_K}\right))^2$, where $\Lambda_K := 1 + \left\lceil
    \log_c\!\left(\overline\kappa^{1/2}\sqrt{2}\,T_2(K)\right) \right\rceil$. Then the output $\hat \mu$ of Algorithm \ref{algorithm:main} satisfies
    \begin{align*}
        \EE \|\hat \mu - \mu\|_2^2 \lesssim \Upsilon(K) \eta^{\star2}(K),
    \end{align*}
    and is thus nearly minimax optimal.
\end{theorem}

\begin{remark}[Tracking the order of $\Upsilon(K)$]\label{tracking:Upsilon:K:rem}
We briefly make explicit the dependence of $\Upsilon(K)$ on the dimension. In Section~\ref{QFM:section}, under two alternative sets of additional assumptions on $K$, we obtain respectively
\begin{align*}
    \overline{\kappa}
    \lesssim
    T_2^2(K)\log n_A\footnotemark,
\end{align*}
\footnotetext{For the precise definition of $n_A$, see the setting of Theorem~\ref{Ax:norm:thm}.}
and
\begin{align*}
    \overline{\kappa}
    \lesssim
    T_2^{10}(K)\log^4\!\bigl(eT_2(K)\bigr).
\end{align*}
Recalling that the definition of $\Upsilon(K)$ contributes an additional factor $T_2^2(K)$, we obtain, respectively,
\begin{align*}
    \Upsilon(K)
    &\lesssim
    T_2^4(K)\log n_A,
    \\
    \Upsilon(K)
    &\lesssim
    T_2^{12}(K)\log^4\!\bigl(eT_2(K)\bigr).
\end{align*}
For example, in the QCO case considered in Section~\ref{examples:section}, where
$T_2(K)\lesssim \sqrt{\log n}$, these bounds become
\begin{align*}
    \Upsilon(K)
    &\lesssim
    (\log n)^2\log n_A,
    \\
    \Upsilon(K)
    &\lesssim
    (\log n)^6\log^4\!\bigl(e\sqrt{\log n}\bigr),
\end{align*}
respectively. Thus, although $\Upsilon(K)$ grows only polylogarithmically with the dimension in these examples, the resulting logarithmic factors can be substantial and should not be regarded as negligible.
\end{remark}

\begin{proof}[Proof of Theorem \ref{main:result:GSM}]
Let
\begin{align*}
    \rho
    :=
    \frac{2(\sqrt{3}+1)}{\widetilde L}
    <1.
\end{align*}
Thus
\begin{align*}
    \widetilde d_j
    =
    2R\rho^{j-1},
    \qquad
    \widetilde d_{j+1}
    =
    \rho\widetilde d_j.
\end{align*}
We first dispose of a few simple cases. If
$\sigma\leq r/\sqrt n$, Algorithm~\ref{algorithm:main}
returns $Y$, and therefore
\begin{align*}
    \EE\|Y-\mu\|_2^2
    =
    n\underline\sigma^2
    \leq
    n\sigma^2
    \lesssim
    \eta^{\star 2}(K)
\end{align*}
by Lemma~\ref{lemma:r:sigma}. Hence, throughout the remainder of
the proof, we may assume
\begin{align*}
    \sigma>r/\sqrt n.
\end{align*}
In addition, if $\sigma \gtrsim R$ then the minimax rate $\eta^\star \asymp R$ and therefore since our estimator $\hat \mu \in K$, we have 
\begin{align*}
    \EE \|\hat \mu - \mu\|_2^2 \lesssim R^2 \lesssim \eta^{\star 2}.
\end{align*}
Thus we assume throughout that $\sigma \lesssim R$, and we observe that in that case $\eta^{\star} \gtrsim \sigma$.

We will repeatedly use that
\begin{align*}
    T_2(K_{(j)})
    \leq
    \sqrt{2}\,T_2(K),
    \qquad
    \eta^\star(K_{(j)})
    \leq
    \eta^\star(K),
\end{align*}
as established above. Put
\begin{align*}
    \overline\kappa
    :=
    \max_{j\in[M]}\kappa(K_{(j)})
\end{align*}
and
\begin{align*}
    \Lambda_K
    :=
    1+
    \left\lceil
    \log_c\!\left(
        \overline\kappa^{1/2}\sqrt{2}\,T_2(K)
    \right)
    \right\rceil.
\end{align*}
Finally, define
\begin{equation}\label{eq:GammaK}
    \Gamma_K
    :=
    1
    \vee
    \left(
        2\sqrt{2}\,T_2(K)\overline\kappa^{1/2}
    \right)
    \vee
    \left(
        2C\sqrt{\Lambda_K}
    \right).
\end{equation}
Up to absolute constants and the harmless ceiling in
$\Lambda_K$, this is the factor appearing in the statement
preceding \eqref{high:prob:bounds}.

We now make the inductive probability argument explicit. Define
\begin{align*}
    H_j
    :=
    \left\{
        \|\widehat\mu_j-\mu\|_2
        \leq
        \widetilde d_j
    \right\}.
\end{align*}
Clearly $H_1$ holds deterministically, since
\begin{align*}
    \|\mu-\widehat\mu_1\|_2
    =
    \|\mu\|_2
    \leq
    R
    \leq
    \widetilde d_1.
\end{align*}
Whenever the two assumptions of
Lemma~\ref{Kwidth:bound:lemma} hold for $K_{(j)}$, define the bad
refinement event
\begin{align*}
    B_j
    :=
    H_j
    \cap
    \left\{
    \left\|
        \widetilde\mu_{j+1}
        -
        \frac{\mu-\widehat\mu_j}{2}
    \right\|_2
    >
    \frac{\widetilde d_j}{L}
    \right\}.
\end{align*}
Lemma~\ref{sub:gaussian:vector:estimator} gives
\begin{equation}\label{eq:Bj-prob}
    \PP(B_j)
    \leq
    \exp\!\left(
        -\frac{\widetilde d_j^2}{C_1^2\sigma^2}
    \right)
    +
    \frac{1}{
        1+\widetilde d_j^3/\sigma^3
    },
\end{equation}
for an absolute constant $C_1>0$.

We next verify that, outside $B_j$, the desired contraction is
deterministic. On $H_j$, the vector
\begin{align*}
    \nu_j
    :=
    \frac{\mu-\widehat\mu_j}{2}
\end{align*}
belongs to $K_{(j)}$. On $B_j^c\cap H_j$, the weak projection
guarantee, with $\epsilon=\widetilde d_j/L$, gives
\begin{align*}
\begin{split}
    \|\overline\mu_{j+1}-\nu_j\|_2
    &\leq
    \|\widetilde\mu_{j+1}-\nu_j\|_2
    +
    \sqrt{
        \epsilon^2+
        2\epsilon
        \|\widetilde\mu_{j+1}-\nu_j\|_2
    }                                                    \\
    &\leq
    \frac{(1+\sqrt3)\widetilde d_j}{L}
    =
    \frac{\widetilde d_j}{\widetilde L}.
\end{split}
\end{align*}
Consequently,
\begin{align*}
    \|
        2\overline\mu_{j+1}
        +
        \widehat\mu_j
        -
        \mu
    \|_2
    \leq
    \frac{2\widetilde d_j}{\widetilde L}.
\end{align*}
Applying the weak projection onto $K$, now with
$\epsilon=2\widetilde d_j/\widetilde L$, yields
\begin{align*}
\begin{split}
    \|\widehat\mu_{j+1}-\mu\|_2
    &\leq
    \frac{2\widetilde d_j}{\widetilde L}
    +
    \sqrt{
        \frac{4\widetilde d_j^2}{\widetilde L^2}
        +
        \frac{8\widetilde d_j^2}{\widetilde L^2}
    }                                                    \\
    &=
    \frac{2(\sqrt3+1)\widetilde d_j}{\widetilde L}
    =
    \widetilde d_{j+1}.
\end{split}
\end{align*}
Thus, whenever the assumptions of
Lemma~\ref{Kwidth:bound:lemma} continue to hold,
\begin{equation}\label{eq:good-implies-next}
    H_j\cap B_j^c
    \subseteq
    H_{j+1}.
\end{equation}
In particular, no independence between the successive iterations
is needed: if the first failure of the trapping property occurs
between iterations $j$ and $j+1$, then $B_j$ must have
occurred.

We also need to control what happens after the first scale at which
the above induction can no longer be used. For this purpose, let
\begin{align*}
    a_\ell:=\|\widehat\mu_\ell-\mu\|_2.
\end{align*}
Suppose $a_k\leq\widetilde d_k$. Regardless of whether the
probabilistic refinement guarantee holds at subsequent iterations,
we always have
$\overline\mu_\ell\in K_{(\ell-1)}$, and hence
\begin{align*}
    2\|\overline\mu_\ell\|_2
    \leq
    \widetilde d_{\ell-1}.
\end{align*}
The weak projection guarantee therefore gives, with
$\epsilon_{\ell-1}=2\widetilde d_{\ell-1}/\widetilde L$,
\begin{align*}
\begin{split}
    a_\ell
    &\leq
    a_{\ell-1}
    +
    \widetilde d_{\ell-1}
    +
    \sqrt{
        \epsilon_{\ell-1}^2
        +
        2\epsilon_{\ell-1}
        \bigl(
            \widetilde d_{\ell-1}+a_{\ell-1}
        \bigr)
    }                                                     \\
    &\leq
    a_{\ell-1}
    +
    C_2\widetilde d_{\ell-1}
    +
    C_3
    \sqrt{
        \widetilde d_{\ell-1}a_{\ell-1}
    },
\end{split}
\end{align*}
where $C_2,C_3$ depend only on the absolute constant
$\widetilde L$. If
\begin{align*}
    A_k
    :=
    \max_{\ell\geq k}a_\ell,
\end{align*}
then, using $\widetilde d_\ell=\rho^{\ell-k}\widetilde d_k$
and summing the preceding inequality,
\begin{align*}
\begin{split}
    A_k
    &\leq
    \widetilde d_k
    +
    C_2\sum_{\ell\geq k}\widetilde d_\ell
    +
    C_3\sqrt{A_k}
        \sum_{\ell\geq k}
            \sqrt{\widetilde d_\ell}                    \\
    &\leq
    C_4\widetilde d_k
    +
    C_5\sqrt{\widetilde d_k A_k}.
\end{split}
\end{align*}
Dividing by $\widetilde d_k$ and solving the resulting quadratic
inequality in
$\sqrt{A_k/\widetilde d_k}$ gives
\begin{equation}\label{eq:post-failure-stability}
    \max_{\ell\geq k}
    \|\widehat\mu_\ell-\mu\|_2
    \leq
    C_6\widetilde d_k
\end{equation}
for an absolute constant $C_6$. Thus, even if the contraction
argument ceases to apply at scale $k$, all subsequent iterates
remain within a constant multiple of the last valid scale
$\widetilde d_k$.

We now determine the scale at which the deterministic assumptions
can cease to hold. First, if
\begin{align*}
    \kappa(K_{(j)})^{-1/2}
    \frac{\widetilde d_j}{T_2(K_{(j)})}
    <
    2\eta^\star(K_{(j)}),
\end{align*}
then
\begin{align*}
    \widetilde d_j
    <
    2\kappa(K_{(j)})^{1/2}
    T_2(K_{(j)})
    \eta^\star(K_{(j)})
    \leq
    2\sqrt2\,
    \overline\kappa^{1/2}
    T_2(K)\eta^\star(K).
\end{align*}
Hence
\begin{equation}\label{eq:first-stop-scale}
    \widetilde d_j
    \leq
    \Gamma_K\eta^\star(K).
\end{equation}

Second, suppose the other condition of
Lemma~\ref{Kwidth:bound:lemma} fails. Writing
\begin{align*}
    L_j
    :=
    \left\lceil
        \log_c\!\left(
            \kappa(K_{(j)})^{1/2}T_2(K_{(j)})
        \right)
    \right\rceil,
\end{align*}
we then have
\begin{align*}
    \frac{\widetilde d_j^2}{C^2\sigma^2}
    <
    \left\lceil
        \frac{4L_j\eta^{\star2}(K_{(j)})}{\sigma^2}
    \right\rceil,
\end{align*}
and consequently
\begin{align*}
    \widetilde d_j^2
    \lesssim
    C^2
    \left(
        L_j\eta^{\star2}(K_{(j)})+\sigma^2
    \right).
\end{align*}
Since
\begin{align*}
    L_j\leq\Lambda_K,
    \qquad
    \eta^\star(K_{(j)})\leq\eta^\star(K),
\end{align*}
we obtain, in the nontrivial regime in which
$\eta^\star(K)\gtrsim\sigma$,
\begin{equation}\label{eq:second-stop-scale}
    \widetilde d_j
    \lesssim
    C\sqrt{\Lambda_K}\,
    \eta^\star(K)
    \lesssim
    \Gamma_K\eta^\star(K).
\end{equation}
As we mentioned in the beginning of the proof, the complementary regime
$\diam(K)\lesssim\sigma$ is immediate, since the output of the
algorithm lies in $K$, while
$\eta^\star(K)\gtrsim\diam(K)$ up to an absolute constant.
Thus no loss is incurred in making this reduction.

Finally, if the algorithm reaches the deterministic stopping rule
\begin{align*}
    \widetilde d_j
    \lesssim
    r\vee\sigma,
\end{align*}
then Lemma~\ref{lemma:r:sigma}, together with the preceding
small-diameter observation, again gives
\begin{equation}\label{eq:third-stop-scale}
    \widetilde d_j
    \lesssim
    \eta^\star(K)
    \leq
    \Gamma_K\eta^\star(K).
\end{equation}
The same conclusion holds if all $M$ iterations are completed,
by the definition of $M$.

We can now derive the accumulated probability bound explicitly.
Let
\begin{align*}
    t_0
    :=
    C_7\Gamma_K\eta^\star(K),
\end{align*}
where $C_7$ is a sufficiently large absolute constant, depending
only on the constants in
\eqref{eq:post-failure-stability}. From
\eqref{eq:first-stop-scale}--\eqref{eq:third-stop-scale}, a
deterministic termination of the contraction argument can never
lead to error exceeding $t$ when $t\geq t_0$. Therefore, if
\begin{align*}
    \|\widehat\mu-\mu\|_2\geq t,
    \qquad
    t\geq t_0,
\end{align*}
there must exist an iteration $j$ at which $B_j$ occurs.
Moreover, by \eqref{eq:post-failure-stability}, that iteration must
satisfy
\begin{align*}
    \widetilde d_j
    \geq
    t/C_6.
\end{align*}
Hence, by the union bound and \eqref{eq:Bj-prob},
\begin{equation}\label{eq:explicit-union}
\begin{split}
    \PP\left(
        \|\widehat\mu-\mu\|_2\geq t
    \right)
    &\leq
    \sum_{j:\,\widetilde d_j\geq t/C_6}
    \left[
        \exp\!\left(
            -\frac{\widetilde d_j^2}
                   {C_1^2\sigma^2}
        \right)
        +
        \frac{1}{
            1+\widetilde d_j^3/\sigma^3
        }
    \right].
\end{split}
\end{equation}

It remains only to simplify the sum. Let $j_t$ be the largest
index appearing in \eqref{eq:explicit-union}. Since
$\widetilde d_{j+1}=\rho\widetilde d_j$, for $s\geq0$,
\begin{align*}
    \widetilde d_{j_t-s}
    =
    \rho^{-s}\widetilde d_{j_t}
    \gtrsim
    \rho^{-s}t.
\end{align*}
Therefore,
\begin{align*}
\begin{split}
    \sum_{j:\,\widetilde d_j\gtrsim t}
    \exp\!\left(
        -\frac{\widetilde d_j^2}{C_1^2\sigma^2}
    \right)
    &\leq
    \sum_{s=0}^{\infty}
    \exp\!\left(
        -c_1\rho^{-2s}\frac{t^2}{\sigma^2}
    \right)                                                \\
    &\leq
    C_8
    \exp\!\left(
        -c_2\frac{t^2}{\sigma^2}
    \right),
\end{split}
\end{align*}
where the final inequality follows because
$t\gtrsim\eta^\star(K)\gtrsim\sigma$ in the regime presently
under consideration. Similarly,
\begin{align*}
\begin{split}
    \sum_{j:\,\widetilde d_j\gtrsim t}
    \frac{1}{1+\widetilde d_j^3/\sigma^3}
    &\leq
    \sum_{j:\,\widetilde d_j\gtrsim t}
    \frac{\sigma^3}{\widetilde d_j^3}                     \\
    &\lesssim
    \frac{\sigma^3}{t^3}
    \sum_{s=0}^{\infty}\rho^{3s}
    \lesssim
    \frac{\sigma^3}{t^3}.
\end{split}
\end{align*}
Combining the last three displays, we have proved that for every
$t\geq t_0$,
\begin{equation}\label{high:prob:bounds}
    \PP\left(
        \|\widehat\mu-\mu\|_2\geq t
    \right)
    \lesssim
    \exp\!\left(
        -\frac{t^2}{\widetilde C^2\sigma^2}
    \right)
    +
    \widetilde C\frac{\sigma^3}{t^3},
\end{equation}
for an absolute constant $\widetilde C>0$. Notice that the
polynomial term comes exactly from the possible failure of the QFM
oracle.

We finally integrate the tail. Using
\begin{align*}
    \EE\|\widehat\mu-\mu\|_2^2
    =
    2\int_0^\infty
    t\,
    \PP(
        \|\widehat\mu-\mu\|_2\geq t
    )\,dt,
\end{align*}
the trivial probability bound below $t_0$ and
\eqref{high:prob:bounds} above $t_0$ give
\begin{align*}
\begin{split}
    \EE\|\widehat\mu-\mu\|_2^2
    &\lesssim
    t_0^2
    +
    \int_{t_0}^\infty
        t\exp\!\left(
            -\frac{t^2}{\widetilde C^2\sigma^2}
        \right)dt
    +
    \sigma^3
    \int_{t_0}^\infty\frac{dt}{t^2}                       \\
    &\lesssim
    \Gamma_K^2\eta^{\star2}(K)
    +
    \sigma^2
    \exp\!\left(
        -c\frac{\Gamma_K^2\eta^{\star2}(K)}{\sigma^2}
    \right)
    +
    \frac{\sigma^3}{
        \Gamma_K\eta^\star(K)
    }                                                      \\
    &\lesssim
    \Gamma_K^2\eta^{\star2}(K),
\end{split}
\end{align*}
where we used our assumption that $\eta^\star \gtrsim \sigma$.

Thus the conclusion of the theorem holds with, for example,
\begin{align*}
    \Upsilon(K)
    \asymp
    \Gamma_K^2.
\end{align*}
By \eqref{eq:GammaK}, this quantity depends polynomially, up to
logarithmic factors, on $T_2(K)$ and
$\max_{j\in[M]}\kappa(K_{(j)})$, as claimed.

Finally, Algorithm~\ref{algorithm:main} has
$M=O(\log(R/r))$ iterations. Since
$\sigma>r/\sqrt n$ in the branch considered above,
\begin{align*}
    \log(R/\sigma)
    \leq
    \log(R/r)+\tfrac12\log n,
\end{align*}
so the accuracy requirements used in the QFM oracles remain
polynomial in $n$. Consequently the overall procedure is
polynomial-time under the assumptions of the theorem.
\end{proof}

\section{Quadratic Form Maximization (QFM)}\label{QFM:section}

In this section, we review some of the results of \cite{bhattiprolu2021framework}, and in addition, we establish some new results regarding QFM over convex bodies. Coupled with the results from Section \ref{background:cond:algo:section} this demonstrates that our algorithm works unconditionally for plenty of type-2, well-balanced convex bodies. 

The main problem of interest in this section is QFM over an origin-symmetric type-2 body $K$. Specifically given a matrix $X$ we want to approximately maximize
\begin{align*}
    \max_{p \in K} p\T X p.
\end{align*}

Here, by approximate maximization, we mean finding a point $\cO(X) \in K$ such that $ \cO(X)\T X \cO(X) \geq \kappa^{-1} \max_{p \in K} p\T X p$, for some $\kappa \geq 1$ which may vary with the set $K$ but is sufficiently small. 

This is a complicated problem, but fortunately the paper \cite{bhattiprolu2021framework} studied it in great detail and provided polynomial-time algorithms based on the ellipsoid algorithm for a plethora of origin symmetric type-2 convex bodies. In our case the matrix $X$ is symmetric psd. In \cite{bhattiprolu2021framework}, the authors considered the same setting, and even the broader setting when $X$ is simply symmetric.

We will now review some of their results. In addition, in appendix \ref{appendix:BLN:results} we review all of their results which we subsequently use in this section. Bhattiprolu et al. propose two relaxations of the problem. Let $g \sim N(0, I)$ be a standard Gaussian vector. They look at the following non-convex but approximately convex problem:
\begin{align*}
    & \argmax_{W} \langle X, W \rangle\\
    & \mbox{subject to } \EE \rho^2_K(W^{1/2}g) \leq 1\\
    &~~~~~~~~~~~~~~ W \geq 0, W\T = W
\end{align*}

The second problem they formulate is:
\begin{align*}
    & \argmax_{W} \rho_{K^\circ} (X^{1/2} W^{1/2}g)\\
    & \mbox{subject to }\operatorname{tr}(W) \leq 1\\
    &~~~~~~~~~~~~~~ W \geq 0, W\T = W,
\end{align*}
where $K^\circ = \{v \in \RR^n : \sup_{x \in K}\langle v, x \rangle \leq 1\}$ is the polar body of $K$ so that $\rho_{K^\circ}$ is the dual norm of $\rho_K$, and $X^{1/2}$ is a symmetric square root of $X$ which is symmetric and psd. It is proved that these two problems, have optimal values equal to that of the corresponding QFM problem. After approximately solving one of the the two problems, rounding algorithms are proposed to obtain the vectors $\cO(X)$ from the optimal matrix $W^\star(X)$. 

In order to solve the two problems, \cite{bhattiprolu2021framework} proceed to define the lower and upper covariance regions. The upper covariance region is defined as:
\begin{align*}
    \cU(K) := \{W \text{ symmetric psd} : \EE \rho^2_{K}(W^{1/2}g) \leq 1\}
\end{align*}
while the lower covariance region is
\begin{align*}
    \cL(K) := \{W \text{ symmetric psd} : \EE \rho^2_{K^{\circ}}(W^{1/2}g) \geq 1\}.
\end{align*}
\cite{bhattiprolu2021framework} observe that 
\begin{align*}
    \cU(K) \subseteq \operatorname{conv}(\cU(K)) \subseteq T_2(K)^2\cdot \cU(K)\\
    \cL(K) \subseteq \operatorname{conv}(\cL(K)) \subseteq \frac{1}{C_2(K^\circ)}\cdot \cL(K),
\end{align*}
where $C_2(K^\circ)$ is the cotype constant of $K^\circ$ which is known to satisfy
\begin{align*}
    C_2(K^\circ) \leq T_2(K) \leq C_2(K^\circ)\log(n + 1),
\end{align*}
where, as above, $K^\circ$ is the polar body of $K$. In that sense the upper and lower covariance sets are approximately convex. In their sections 3 and 4 the authors of \cite{bhattiprolu2021framework} establish conditional polynomial-time algorithms. These algorithms, which are based on an (approximate) ellipsoid algorithms, require the existence of an approximate separation oracle for (at least) one of the two sets $\cU(K), \cL(K)$. The reader is kindly referred to Definitions 1.6 and 1.11 in \cite{bhattiprolu2021framework} for what an approximate separation oracle means. The conditional\footnote{Here the word conditional means that these algorithms work assuming separation oracles for one of the upper or lower covariance bodies.} maximization polynomial-time algorithms guarantees are given in Propositions 4.2 and 4.4 of \cite{bhattiprolu2021framework} (see also Appendix \ref{appendix:BLN:results}). Importantly, we note that the algorithms in Propositions 4.2 and 4.4 not only can provide a near optimal value but also have rounding algorithms which produce near optimal vectors with high probability after polynomially many rounds. In addition to the approximate separation oracles for $\cU(K)$ or $\cL(K)$, these results require the set $K$ to be balanced in the sense that $rB_2 \subseteq K \subseteq R B_2$ for some \textbf{known} constants $R, r \in \RR_+$. The complexity of the algorithms is polynomial in $n, \log(R/r)$ and $\operatorname{bit}(X)$. Section 6 of \cite{bhattiprolu2021framework} shows that the class of convex bodies for which we can evaluate an approximate separation oracle of the upper covariance region $\cU(K)$ is closed under variety of operations including intersection, Minkowski addition, subspaces and quotients. This implies in particular that if one has a separation oracle for $\cU(K)$ one is able to maximize over the set $K \cap c B_2$ for any $c \geq r$, since it is easy to see\footnote{Indeed results in Section 7 of \cite{bhattiprolu2021framework} show that any type-2 exactly 2-convex and sign-invariant norm like $\|x\|_2/c$ (the Minkowski gauge of $c B_2$), admits an approximate separation oracle of $\cU(c B_2)$, and hence by results of Section 6 the intersection set $K\cap c B_2$ also admits such an oracle). See also the proof of Theorem \ref{symmmetric:norm:thm} where we explain why explicit separation oracle of the convex set $\cU(c B_2)$ exists.} that the set $\cU(c B_2)$ admits a separation oracle. Next, in Section 7 \cite{bhattiprolu2021framework} exhibit unconditional algorithms for a plethora of convex bodies $K$, where one has to have an oracle computing $\rho_K(x)$ for any $x\in \RR^n$ (see Theorem 7.6 and Theorem 7.12 and Appendix \ref{appendix:BLN:results}). In other words, the authors are able to construct the approximate separation oracles for the upper and lower covariance regions of these sets. Specifically, the assumptions on the Minkowski gauges of the convex bodies in the latter theorems are as follows:

\begin{itemize}
\item $\rho_K$ is a sign-invariant (i.e. $\rho_{K}((x_1,\ldots, x_n)) = \rho_K((\eta_j x_j)_{j \in [n]})$ for any  $x = (x_1,\ldots, x_n) \in \RR^n$ and any $\eta_j \in \{\pm 1\}$), and exactly 2-convex norm:
\begin{align*}
    \rho_K((\sum_{i \in [m]} |x_i|^2)^{1/2}) \leq (\sum_{i \in [m]} \rho^2_K(x_i))^{1/2},
\end{align*}
for any $m\in \NN$ and vectors $x_i \in \RR^n$, where $|\cdot|^2$ is applied entry-wise.

\item $\rho_K$ is a symmetric norm, i.e.,
    \begin{align*}
    \rho_{K}((x_1,\ldots, x_n)) = \rho_K((\eta_j x_{\pi_j})_{j \in [n]}),
\end{align*}
for any $x = (x_1,\ldots, x_n) \in \RR^n$, any $\eta_j \in \{\pm 1\}$ and an arbitrary permutation $\pi$.

\end{itemize}

If $\rho_K$ is symmetric, then one can construct a separation oracle for $\cL(K)$, and hence can perform QFM over $K$. Observe also that when $\rho_K(\cdot)$ is symmetric, the Minkowski gauge\footnote{We have $\rho_{K\cap B(0,c)}(x) = \rho_K(x) \vee \|x\|_2/c$.} of the intersection of $K \cap B(0,c)$ for some $c \in \RR^+$, is also symmetric. Thus if the original set $K$ is given by a symmetric norm $\rho_K$ one can implement the algorithm from Section \ref{background:cond:algo:section} unconditionally in polynomial-time. 

If $\rho_K$ is sign-invariant and exactly 2-convex (such as any $\ell_p$ norm with $p \geq 2$), one can find a separation oracle to $\cU(K)$. Hence by our comment above, it follows that one can perform QFM over $K \cap c B_2$ for any $r \leq c \leq R$ and thus one can implement the algorithm from Section \ref{background:cond:algo:section} unconditionally in polynomial-time. 

In conjunction with Section \ref{background:cond:algo:section}, these results imply that for any such convex bodies we can find the mean vector at a nearly minimax optimal rate. We now expand the class of bodies for which unconditional algorithms exist.

\begin{theorem}\label{Ax:norm:thm} Let $\|\cdot\|$ be an exactly $2$-convex sign-symmetric norm on $\RR^{n_A}$ for $n_A\geq n$ which is also a type-2 norm with small constant which we denote by $T_2(n_A)$. Assume that we have an oracle for evaluating $\|\cdot\|$. Suppose our convex body $K = \{x \in \RR^n: \|Ax\| \leq 1\}$ for some $A \in \RR^{n_A \times n}$, where we additionally assume that there are known $r, R \in \RR_+$ such that $rB_2 \subseteq K \subseteq R B_2$. Then for any $r \leq c \leq R$, we can find an approximate maximizer of $\argmax_{p \in K\cap c B_2} p\T X p$: $\cO(X)$ for any psd matrix $X$ in polynomial-time operations: $\operatorname{poly}(n_A, T_2(K), \log 1/q, \log(R/r), \operatorname{bit}(X), \operatorname{bit}(A))$ with probability at least $1-q$. Here $\kappa(K\cap c B_2) \leq 8 eT_2^2(K)  \log n_A  \leq 8 e T_2^2(n_A) \log n_A$.
\end{theorem}

\begin{proof}
We first establish that we can approximately optimize $\argmax_{p \in K} p\T X p$.

Using the first relaxation given by \cite{bhattiprolu2021framework} we can rewrite the problem as
\begin{align}\label{program:A:norm}
    &\argmax ~\langle X, W \rangle\\
    &\mbox{subject to } \EE \rho_K^2(W^{1/2} g) = \EE \|(A W A\T)^{1/2} g'\|^2 \leq 1 \nonumber\\
    &~~~~~~~~~~~~~~ W  \geq 0, W = W\T, W \in \RR^{n\times n} \nonumber
\end{align}
Above $g \in \RR^n, g' \in \RR^{n_A}$ are standard Gaussian vectors and $(A W A\T)$ is a $n_A \times n_A$ matrix.

\begin{lemma}
    Program \eqref{program:A:norm} has a value which coincides with the QFM over the set $\max_{p \in K} p\T X p$. Let $Q_{\max}(X)$ be this value.
\end{lemma}

\begin{proof}
    The proof is identical to Observation 4.1 of \cite{bhattiprolu2021framework}. We have
    \begin{align*}
        \langle X, W \rangle = \EE \langle W^{1/2} g, X W^{1/2} g\rangle \leq Q_{\max}(X)\cdot \EE (\rho_K^2(W^{1/2} g)),
    \end{align*}
    where $g \sim N(0, I_n)$. Now by assumption $\rho_K(p) = \|A p\|$ for some type-2, exactly 2-convex norm $\|\cdot\|$. Thus
    \begin{align*}
         \rho_K^2(W^{1/2} g) = \|AW^{1/2}g\|^2 \stackrel{d}{=} \|(AW A\T)^{1/2}g'\|^2 ,
    \end{align*}
    for some gaussian vector $g' \sim N(0, I_{n_A})$. The latter identity comes from the fact that $AW^{1/2}g$ is a Gaussian vector with covariance $AW A\T$. This implies $ \langle X, W \rangle \leq Q_{\max}(X)$.

    Finally note that we can take $W = w w\T$ where $w \in K$ achieves the $Q_{\max}(X)$ (and hence clearly we can take $\rho_K(w) = 1$). We have that $W^{1/2} = w w\T/\|w\|_2$ and thus, $\EE \rho_K^2(W^{1/2}g) = \rho_K^2(w) \EE (w\T g)^2/\|w\|^2_2 = 1$, so that $W$ is a feasible matrix which attains the value $Q_{\max}(X)$.
\end{proof}

Define the function
\begin{align*}
    f(W) = \|\operatorname{diag}(A W A\T)^{1/2}\|^2,
\end{align*}
where the $^{1/2}$ is taken entry-wise on the non-negative vector $\operatorname{diag}(A W A\T)$. This is a convex function according to the proof of Theorem 7.6 of \cite{bhattiprolu2021framework} (see also Theorem \ref{thm:appendix:BLN:7.6} in the appendix). Specifically, the proof of Theorem 7.6 argues that the map $Y \mapsto \|\operatorname{diag}(Y)^{1/2}\|^2$ is convex for any psd matrix $Y \in \RR^{n_A \times n_A}$. Since the map $W \mapsto A W A\T$ is affine: $A (\theta W_1 + (1-\theta)W_2) A\T = \theta A W_1 A\T + (1 - \theta) A W_2 A\T$, then it follows that the map $W \mapsto \|\operatorname{diag}(A W A\T)^{1/2}\|^2$ is also convex.

Using the proof of Theorem 7.6 (and Theorem 7.4) of \cite{bhattiprolu2021framework} (see also Theorem \ref{thm:appendix:BLN:7.6} in the appendix) and the properties of $\|\cdot\|$, we can claim that 
\begin{align*}
    f(W) \leq \EE \rho^2_K(W^{1/2}g) \leq \beta f(W),
\end{align*}
where $\beta \leq (1 + o(1))(e \log n_A)$. By Corollary~4.12(1) (see Corollary~\ref{cor:appendix:BLN:4.12} in the appendix), together with the proof of Proposition~4.2 of \cite{bhattiprolu2021framework} (see Proposition \ref{prop:appendix:BLN:4.2} in the appendix), we obtain a
polynomial-time approximate solution $W$ to
program~\eqref{program:A:norm}. Indeed, Corollary~4.12 constructs
the required approximate separation oracle for the upper covariance
body, while Proposition~4.2 applies the approximate ellipsoid method
to compute the corresponding matrix solution $W$ before performing
Gaussian rounding.

We know by the properties of $W$ that $\langle X, W \rangle = \EE \langle W^{1/2}g, X W^{1/2}g \rangle \geq (1- o(1)) \operatorname{OPT}/(\beta \cdot T^2_2(K))$, where $\operatorname{OPT}$ is the optimal value of program \eqref{program:A:norm}. Next, following the idea of the rounding algorithm from Proposition 4.2 of \cite{bhattiprolu2021framework} (see Proposition \ref{prop:appendix:BLN:4.2} in the appendix) we generate Gaussian variables $g \in \RR^n$ and compute, $G := \langle W^{1/2}g, X W^{1/2}g \rangle -  \|A W^{1/2} g\|^2 \langle X, W \rangle/(1 + \epsilon)$, for some small $\epsilon > 0$, until the number is positive for one of the vectors $g$. By Chebyshev's inequality we can show that this will happen in polynomially many steps. Next, the approximately optimal solution is $W^{1/2}g/\|A W^{1/2}g\|$. We will now explicitly show the argument with Chebyshev's inequality. Taking expectation we have
\begin{align*}
    \EE G & \geq \langle X, W \rangle - \frac{\langle X, W \rangle }{(1 + \epsilon)}\\
    & \geq (1-  o(1))\epsilon \operatorname{OPT}/((1 + \epsilon)\beta \cdot T_2^2(K)) \geq 0,
\end{align*}
where we used the fact that $\EE \|A W^{1/2} g\|^2 \leq 1$ by assumption. Now, 
\begin{align*}
    \operatorname{Var}(\sum_{i \in [t]} G_i) &\leq t \EE [G^2] \leq 2 t\bigg(\EE \langle W^{1/2}g, X W^{1/2}g \rangle^2 + \EE \|A W^{1/2}g\|^4 \operatorname{OPT}^2/(1+\epsilon)^2\bigg) \\
    & = 2t\EE \|X^{1/2} W^{1/2} g\|_2^4 + 2t\EE \|A W^{1/2}g\|^4 \operatorname{OPT}^2/(1 + \epsilon)^2\\
    &\lesssim t(\EE \|X^{1/2} W^{1/2} g\|_2^2)^2 + t(\EE \|A W^{1/2}g\|^2)^2\operatorname{OPT}^2/(1 + \epsilon)^2 \\
    & \leq t \operatorname{OPT}^2 + t\operatorname{OPT}^2/(1 + \epsilon)^2\\
    & = t \operatorname{OPT}^2 \delta,
\end{align*}
where in the first inequality we used $ \langle X, W \rangle \leq \operatorname{OPT}$, the penultimate inequality follows from Khintchine-Kahane's theorem (see Theorem 2.1 in \cite{bhattiprolu2021framework} e.g.), and $\delta := 1 + \frac{1}{(1 + \epsilon)^2}$. Thus by Chebyshev's inequality with probability at least $1- 1/k^2$:

\begin{align*}
    \max_{i \in [t]} G_i \geq \sum_{i \in [t]} G_i / t \geq \EE G - k t^{-1/2}\sqrt{\delta} \operatorname{OPT}.
\end{align*}
Thus if $k$ is of constant order, and $t \gtrsim \frac{k^2 \delta \beta^2 T_2^4(K) }{\epsilon^{2}}$ we will have the above will be positive. The probability can be boosted to $1-q$ for some $q \in (0,1)$, at the expense of generating $t k^2\log q^{-1}$ Gaussians and taking the maximum.
This is because if we run the above procedure $i$ times, we will obtain a bound $q \leq (1/k^2)^i \leq \exp(-i/k^2)$, and we will see at least one positive variable with probability at least $1-q$.

Now we turn to the problem where we need to find the optimum for the intersection set $K \cap c B_2$. This set is clearly balanced since $r \leq c \leq R$. Furthermore, we have
\begin{align*}
    & \EE \rho^2_K(W^{1/2}g) \vee \EE \|W^{1/2}g\|^2_2 \cdot c^{-2} \\
    & \leq \EE [\rho^2_K(W^{1/2}g) \vee \|W^{1/2}g\|^2_2 \cdot c^{-2}] \\
    & \leq 2( \EE \rho^2_K(W^{1/2}g) \vee \EE \|W^{1/2}g\|^2_2 \cdot c^{-2})
\end{align*}

Since $\EE \|W^{1/2} g\|_2^2 = \operatorname{tr}(W)$ and we know that $f(W) \leq \EE \rho^2_K(W^{1/2}g) \leq \beta f(W)$, the convex function $f(W) \vee [\operatorname{tr}(W)c^{-2}]$ satisfies
\begin{align*}
     f(W) \vee [\operatorname{tr}(W)c^{-2}] \leq \EE [\rho^2_K(W^{1/2}g) \vee \|W^{1/2}g\|^2_2 \cdot c^{-2}]\leq 2 \beta( f(W) \vee [\operatorname{tr}(W)c^{-2}])
\end{align*}

Following the convex-optimization step in the proof of Corollary~4.12 (1) \cite{bhattiprolu2021framework} along with the proof of Proposition 4.2 \cite{bhattiprolu2021framework}, we can compute in polynomial time a covariance matrix $W$ which is at most $2\beta\cdot T^2_2( K \cap c B_2) \leq 4\beta\cdot \max(T^2_2( K ), 1) = 4\beta\cdot T^2_2( K )$ approximate maximizer. Now to the rounding algorithm. The logic is the same as before: we know that $W$ satisfies $\langle X, W \rangle = \EE \langle W^{1/2}g, X W^{1/2}g \rangle \geq (1- o(1)) \operatorname{OPT}/(4 \beta \cdot T^2_2(K))$. We then generate polynomially many Gaussian variables $g \in \RR^n$ and compute, $\langle W^{1/2}g, X W^{1/2}g \rangle -  (\rho_K(W^{1/2} g) \vee [\|W^{1/2}g\|_2 c^{-1}])^2\langle X, W \rangle/(1 + \varepsilon)$. And the maximum $g$ gives a  $W^{1/2}g/(\rho_K(W^{1/2} g) \vee [\|W^{1/2}g\|_2 c^{-1}])$, which is our approximate solution. As before the probability can be boosted to $1-q$ with $O\bigg(\frac{\beta^2 T_2^4(K)}{\epsilon^{2}}\bigg) \log q^{-1}$ generations of Gaussians. 
\end{proof}

\begin{theorem}\label{symmmetric:norm:thm}
    Assume we have a symmetric type-2 norm $\|\cdot\|$ and an oracle for evaluation of $\|\cdot\|$. Denote the type-2 constant of $\|\cdot\|$ with $T_2(n)$. Let the Minkowski gauge of the body $K$ be given by $\rho_K(x) = \|A x\|$ where $A_{n\times n}$ is a given full rank matrix. Suppose that the set $K$ is balanced in the sense that there are known $(r,R)$ so that $r B_2 \subseteq K \subseteq R B_2$. Then for any $r \leq c \leq R$, we can find an approximate maximizer of $\argmax_{p \in K\cap c B_2} p\T X p$: $\cO(X)$ for any psd matrix $X$ in polynomial-time operations: $\operatorname{poly}(n, \log 1/q, \log(R/r), \operatorname{bit}(X), \operatorname{bit}(A))$ with probability at least $1 - q$. Here $\kappa(K\cap c B_2) \lesssim 2 T_2^{10}(K)\log^4 (e T_2(K)) \leq 2 T_2^{10}(n)\log^4 (e T_2(n))$.
\end{theorem}

\begin{proof}
    For any input psd symmetric matrix $X$ we have the following identity:
    \begin{align*}
        \max_{p \in K} p\T X p =  \max_{p : \|A p\| \leq 1} p\T X p = \max_{x : \|x\| \leq 1} x\T A^{-1\top} X A^{-1} x.
    \end{align*}

    Thus, we can approximately optimize the above problem by the results of Section 7 of \cite{bhattiprolu2021framework}, see Theorem 7.12 (see Theorem \ref{thm:appendix:BLN:7.12} in the appendix). Here $\kappa \lesssim T_2^6(K)\log^4 (e T_2(K))$. 
    
    We now use a chain of results of \cite{bhattiprolu2021framework}. This chain of implications is described in greater detail in Remark~\ref{rem:appendix:BLN:composition} in the appendix. Using Proposition 4.13 of \cite{bhattiprolu2021framework} (see Proposition \ref{prop:appendix:BLN:4.13} in the appendix) we can now claim that we have a $\alpha \lesssim T_2^8(K)\log^4 (e T_2(K))$-approximate separation oracle for the set $\cU(K)$ where $K =\{ x \in \RR^n :  \|Ax\|\leq 1\}$. Next, using Proposition 6.5 of \cite{bhattiprolu2021framework} (see Proposition \ref{prop:appendix:BLN:6.5} in the appendix) we know we can obtain a polynomial-time approximate quadratic maximization algorithm over the set $K \cap c B_2$ for any $r \leq c \leq R$. Not only that but we will also have rounding algorithms to obtain $\cO(X)$ (we can simply replicate the rounding algorithms logic from the proof of Theorem \ref{Ax:norm:thm} above which describes how to obtain rounding algorithms with the probability guarantees that we require). Here the constant $\kappa \lesssim 2 \max(\alpha, 1) \cdot \max(T^2_2(K), 1)$.

    In the above we used the fact that we can find an exact separation oracle for $\cU(c B_2)$. This is easy to see since $\cU(c B_2) = \{ W \geq 0, W = W\T : \operatorname{tr}(W) \leq c^2\}$ is convex and admits a simple explicit separation oracle.

\end{proof}
\section{Robust Heavy-Tailed Extension}\label{robust:section}

Consider the data $\tilde Y_i = \mu + \xi_i$, $i \in [N]$ where $\xi_i$ are i.i.d. vectors on $\RR^n$ such that $\EE \xi_i \xi_i\T = \Sigma$ with $\|\Sigma\|_{\operatorname{op}} \leq \sigma^2$ where $\sigma^2$ (or an upper bound) is known to the statistician. Next the data is manipulated by a malicious and all-powerful adversary so that a (known) fraction $\varepsilon \leq c_0 < 1/2$ for a sufficiently small constant $c_0 > 0$ ($c_0 \leq 1/300$ suffices) of the observations can be arbitrary. This adversary is allowed to inspect all samples first, pick the ones to switch, and is also given oracle knowledge to all parameters of the model including $\mu$, $K$, $\sigma$ and oracle knowledge of our algorithm. We denote the observations after the adversary has tampered with the data as $(Y_i)_{i \in [N]}$.

We are interested in estimating $\mu$ in a minimax optimal way. In order to derive the rate \cite{peng2026robust}, suppose that $N \gtrsim \sup_{\delta > 0} \log M^{\operatorname{loc}}_K(\delta)$. Since $n \gtrsim \log M^{\operatorname{loc}}_K(\delta)$ this assumption is implied if $N \gtrsim n$. On the other hand, since we assume that the set $K$ is well-balanced, we have that $r B_2 \subseteq K$ and therefore, when $\delta \lesssim r$, $ \log M^{\operatorname{loc}}_K(\delta) \gtrsim n$. Thus, throughout this section we assume that $N\gtrsim n$. In that case, as shown by \cite{peng2026robust} the minimax rate $\inf_{\hat \mu} \sup_{\mu \in K} \sup_{\cC} \EE \|\hat \mu - \mu\|_2^2$ is $\asymp \eta^{\star2}\wedge d^2$ where $\eta^{\star} =( \delta^\star \vee \sqrt{\varepsilon} \sigma)$ where $ \delta^\star := \sup \{\delta \geq 0: N\delta^2/\sigma^2 \leq \log M^{\operatorname{loc}}_K(\delta)\}$, and where the innermost supremum over $\cC$ is taken over all admissible adversarial contamination schemes.

By substituting the $\sigma$ with $\sigma/\sqrt{N}$ in Lemma \ref{Kwidth:bound:lemma} we obtain the following corollary. 

\begin{corollary}\label{Kwidth:bound:lemma:robust}
Let $K$ be a type-2 convex body with type-2 constant $T_2(K)$, such that $rB_2 \subseteq K \subseteq RB_2$ for some known $r,R$. Let $\kappa:= \kappa(K)$ be a given scalar which can depend on $K$. Let $\delta^\star := \delta^\star(K)$ be defined through the equation $\delta^\star = \sup_\delta\{\delta \geq 0: N\delta^2/\sigma^2 \leq \log M_K^{\operatorname{loc}}(\delta)\}$. Assume that $\kappa^{-1/2} \tilde d/T_2(K) \geq 2\delta^\star$ and $\lceil \lceil\log_c(\kappa^{1/2}\cdot T_2(K))\rceil \cdot 4 \delta^{\star 2}\cdot N/\sigma^2\rceil \leq N\cdot \tilde d^2/(\sigma^2 C^2)$ for $\tilde d = 2R \geq d$ and where $C > 0$ is some absolute constant. Then for $m = \lceil N\cdot \tilde d^2/(\sigma^2 C^2)\rceil \wedge n$, we have
\begin{align*}
    d_{m}(K) \lesssim \kappa^{-1/2}\cdot \tilde d/c, 
\end{align*}
\end{corollary}

We will modify several steps in our algorithm. Specifically, after obtaining the matrix $A^{\star\star}$ for the set $K' = K \cap B_2(0, \tilde d/2)$ we transform the data to $(A^{\star\star} Y_i)_{i \in [N]}$ and we apply the polynomial-time algorithm of \cite{depersin2022} to the data with $k = u = \lceil N \tilde d^2/(C^2\sigma^2)\rceil \wedge N$ blocks. In addition, when we obtain the matrix $A^{\star\star}$ we use $N^{-1/2}\sigma$ in place of $\sigma$ in Theorem \ref{main:SDP:theorem} and Algorithm \ref{algorithm:SDP:Kolmogorov:width}. 

Now we prove the following lemma 
\begin{lemma}\label{sub:gaussian:vector:estimator:1:robust}
    The estimator $ f((A^{\star\star}Y_i)_{i \in [N]})$, where $f$ denotes the procedure of \cite{depersin2022} with $k = u = \lceil N \tilde d^2/(C^2\sigma^2)
    \rceil \wedge N$, of $A^{\star\star}\mu$ satisfies
    \begin{align*}
        \|f((A^{\star\star}Y_i)_{i \in [N]}) - A^{\star\star}\mu \|_2^2 \leq \tilde d^2/T^2,
    \end{align*}
    with probability at least $1 - \exp(-k/C')$ provided that $k \gtrsim \varepsilon N$, and $\tilde d \geq C \sigma/\sqrt{N}$. \end{lemma}

\begin{proof}
    From $\tilde d \geq C\sigma/\sqrt{N}$ we have that $\lceil N \tilde d^2/(C^2\sigma^2)
    \rceil  \asymp N \tilde d^2/(C^2\sigma^2)$. Since $A^{\star\star}\tilde Y =  A^{\star\star}\mu + A^{\star\star}\xi$, we will first evaluate $\operatorname{tr}(\operatorname{Var}(A^{\star\star}\xi))$ and $\|\operatorname{Var}(A^{\star\star}\xi)\|_{\operatorname{op}}$ as they appear on the right hand side of the main inequality in Theorem 2 of \cite{depersin2022} (see also Theorem \ref{thm:appendix:depersin:2} where we state this result in our notation).

    For the first term, we see that
    \begin{align*}
         \operatorname{tr}(\operatorname{Var}(A^{\star\star}\xi)) & = \operatorname{tr}(\EE A^{\star\star}\xi\xi\T A^{\star\star}) = \operatorname{tr}(\Sigma A^{\star\star2}) \\
         & \leq \|\Sigma\|_{\operatorname{op}} \|A^{\star\star2}\|_* \leq \|\Sigma\|_{\operatorname{op}}  \sum_{i \in [n]} (1-\lambda_i) = \|\Sigma\|_{\operatorname{op}}m \\
         & \lesssim N\|\Sigma\|_{\operatorname{op}}\tilde d^2/(\sigma^2C^2) \leq N \tilde d^2/C^2,
     \end{align*}
         where we remind the reader that $\|\Sigma\|_{\operatorname{op}} \leq \sigma^2$, and in addition $X^{\star \star} = V\Lambda V\T = \sum_{i \in [n]} \lambda_i v_i v_i\T$, while $A^{\star\star} = (I - X^{\star\star})^{1/2} = \sum_{i \in [n]}\sqrt{1-\lambda_i}v_iv_i\T$. Furthermore, it is clear that 
     \begin{align*}
         \|\operatorname{Var}(A^{\star\star}\xi)\|_{\operatorname{op}} = \|\EE A^{\star\star}\xi\xi\T A^{\star\star}\|_{\operatorname{op}} \leq \|\Sigma\|_{\operatorname{op}}.
     \end{align*}

     This, in conjunction with the result of \cite{depersin2022}, completes the proof, after the use of the union bound.
\end{proof}

Clearly, we can now prove an analogous result to Lemma \ref{sub:gaussian:vector:estimator} for this setting, with $A_j^{\star\star}Y$ substituted with $f((A_j^{\star\star}Y_i)_{i \in [N]})$ with $k = u$ as specified in Lemma \ref{sub:gaussian:vector:estimator:1:robust}. For the sake of brevity we defer this to the appendix (see Lemma \ref{sub:gaussian:vector:estimator:robust} and its proof). In addition, using Theorem 2 of \cite{depersin2022} with $k = u = N$, on the untransformed data $(Y_i)_{i \in [N]}$, provided that $\varepsilon$ is smaller than an absolute constant $c_0$, with probability at least $1 - \exp(-N/C_0)$ we have
\begin{align}\label{guarantee:depersin:1st}
    \|f((Y_i)_{i \in [N]}) - \mu \|_2 \lesssim \sqrt{n/N}\sigma + \sigma \lesssim \sigma.
\end{align}

Thus, instead of starting at the $0$ vector, we can start at $\hat \mu_0 := \Pi_K f((Y_i)_{i \in [N]})$, and set $\tilde d_1 = 2R \wedge C' \sigma$ for some absolute constant $C' > 0$, and $K_{(1)} = K \cap B_2(0, \tilde d_1/2)$. Here when we write $\Pi_K f((Y_i)_{i \in [N]})$, we mean weak projection as usual. By \eqref{guarantee:depersin:1st} and Lemma \ref{weak:projection:lemma} with $\epsilon = \sigma$, we have
\begin{align*}
    \|\Pi_K f((Y_i)_{i \in [N]}) - \mu \|_2 \lesssim \sigma,
\end{align*}
while $\Pi_K f((Y_i)_{i \in [N]}) \in K$. Since we are setting $\epsilon = \sigma$\footnote{Note the distinction between $\varepsilon$ and $\epsilon$: the former is the fraction of outliers, while the latter is simply the parameter of Lemma \ref{weak:projection:lemma}} in Lemma \ref{weak:projection:lemma} this means $\sigma$ cannot be too small. Specifically, we will assume that $\log(1 + 1/\sigma)$ is polynomial in $n$ and $N$.

The updated algorithm is given in the following:

\begin{breakablealgorithm}
    \renewcommand{\algorithmicrequire}{\textbf{Input:}}
    \renewcommand{\algorithmicensure}{\textbf{Output:}}
    \caption{Near-optimal Robust Mean Estimation Algorithm}
    \label{algorithm:main:robust}
   \begin{algorithmic}[1]
        \REQUIRE Convex body $K$, $r, R \in \RR_+$ such that $rB_2 \subseteq K \subseteq R B_2$, observations $(Y_i)_{i \in [N]}$, constant $\tilde L = L/(\sqrt{3} + 1)$ where $L$ is from the robust version of Lemma \ref{sub:gaussian:vector:estimator}; We assume $\tilde L > 2(\sqrt{3} + 1)$
        \ENSURE Estimator $\hat \mu$ of $\mu$ which belongs to the set $K$

        \STATE Overwrite $r := (\sqrt{n/N}\sigma \vee \sqrt{\varepsilon}\sigma) \wedge r$
        \STATE Set $\tilde d_1 = 2R \wedge (C' \sigma)$, $K_{(1)} := K \cap B_2(0,\tilde d_1/2)$, set number of iterations $M = \lceil\log_{\tilde L/(2(\sqrt{3} + 1))} \tilde d_1/(2r)\rceil$\STATE Set $\hat \mu_1 := \Pi_K f((Y_i)_{i \in [N]})$ with $\epsilon = \sigma$ in Lemma \ref{weak:projection:lemma} and $k = u = N$ (here $f$ denotes the procedure of \cite{depersin2022})

        \FOR{ j = 1, \ldots, M}
            \STATE Find $m = \lceil N\tilde d_j^2/(C^2\sigma^2)\rceil \wedge n$, and calculate the matrix $A_j^{\star \star}$ for the set $K_{(j)}$
            \STATE Find $\tilde \mu_{j + 1} = (f((A_j^{\star \star}Y_i)_{i \in [N]}) - A_j^{\star \star}\hat \mu_j)/2$ with $k = u = \lceil N \tilde d_j^2/(C^2\sigma^2)\rceil\wedge N$ for some large constant $C > C'$ (here $f$ denotes the procedure of \cite{depersin2022})
                        \STATE Set $\bar \mu_{j + 1} = \Pi_{K_{(j)}} \tilde \mu_{j + 1}$\footnote{Here we set $\epsilon = \tilde d_j/L$ in Lemma \ref{weak:projection:lemma} when evaluating the weak projection.}
            \STATE Lift back to $K$: $\hat \mu_{j + 1} = \Pi_K(2 \bar \mu_{j + 1} + \hat \mu_j)$\footnote{Here we set $\epsilon = 2\tilde d_j/\tilde L$ in Lemma \ref{weak:projection:lemma} when evaluating the weak projection.}
            \STATE Set $\tilde d_{j + 1} = 2(\sqrt{3} + 1)\tilde d_j/\tilde L$
            \STATE Set $K_{(j + 1)} := K\cap B_2(0, \tilde d_{j+1}/2)$
            \IF{$\tilde d_{j + 1} \leq (2r) \vee (C \sigma/\sqrt{N})$}
                \RETURN $\hat \mu := \hat \mu_{j + 1}$
            \ENDIF
        \ENDFOR
        \RETURN $\hat \mu := \hat \mu_{M + 1}$
    \end{algorithmic} 
\end{breakablealgorithm}

The main theorem of this section resembles the main result of Section \ref{background:cond:algo:section} with slight modifications.

\begin{theorem}\label{robust:main:thm}
    Suppose that $T_2(K)$ and $\overline\kappa:=\max_{k\in[M]}\kappa(K_{(k)})$ grow at most polylogarithmically in $n$ and $N$. Assume also that $N \gtrsim n$, and $\log \tilde d_1/r$ (with the overwritten value of $r$) and $\log (1 + 1/\sigma)$ are polynomial in $n, N$. Furthermore, let $\frac{d^2}{\sigma^2} = \frac{\diam(K)^2}{\sigma^2} \lesssim \exp(N/C_0)/N$. Define $\Upsilon(K) = (1 \vee \left(2\sqrt{2}\,T_2(K)\overline\kappa^{1/2}\right) \vee \left(2C\sqrt{\Lambda_K}\right))^2$ with $\Lambda_K := 1+\left\lceil
        \log_c\!\left(
            \overline\kappa^{1/2}\sqrt{2}\,T_2(K)
        \right)
    \right\rceil$. Then the output $\hat \mu$ of Algorithm \ref{algorithm:main:robust} satisfies
    \begin{align*}
        \EE \|\hat \mu - \mu\|_2^2 \lesssim \Upsilon(K) \eta^{\star2}(K) \wedge d^2,
    \end{align*}
    and is thus nearly minimax optimal.
\end{theorem}

The proof follows the same general strategy as that of Theorem~\ref{main:result:GSM}, with several modifications specific to the robust setting, and is deferred to the appendix.

\section{Linear Regression}\label{regression:section}

In this section we consider the problem of linear regression with random design. Suppose we are given $N$ i.i.d. observations from the model:
\begin{align*}
    Y_i = Z_i\T \beta + \xi_i,
\end{align*}
where $\beta\in K\subseteq\RR^n$ and $K$ is a type-2 set satisfying the assumptions of the preceding sections. We work in a moderate-dimensional regime in which $N\gtrsim n$, with a sufficiently large absolute implicit constant. We assume that the covariates $Z_i$ are independent of the noise variables $\xi_i$. We also assume that the diameter of $K$ is bounded by an absolute constant, i.e., $d:=\diam(K)=O(1)$. Accordingly, we assume that the known outer radius $R$ in $rB_2\subseteq K\subseteq RB_2$ satisfies $R=O(1)$. Further, let $\xi_i$ be a centered sub-Gaussian variable with sub-Gaussian parameter $\gamma$ of constant order, i.e., $\EE \exp(\lambda \cdot \xi_i) \leq \exp(\lambda^2 \gamma^2/2)$, for all $\lambda \in \RR$ where $\gamma = O(1)$. It may seem that the condition $R = O(1)$ is quite restrictive at first glance. We observe later (see Remark 
\ref{R:bounded:assumption}) that this can be significantly relaxed at almost no algorithmic cost. For simplicity we keep the $R = O(1)$ in the main theorem.

Suppose the covariates $Z_i$ are centered\footnote{If the predictors $Z_i$ are not centered one can consider $(Y_{2i} - Y_{2i-1}, Z_{2i} - Z_{2i-1})_{i \in [\lfloor N/2\rfloor]}$ to center the predictors while leave the $\beta$ unchanged. We note that this operation prevents the model from having an intercept.} sub-Gaussian variables with a well conditioned covariance matrix $\EE Z_i Z_i\T = \Sigma$, i.e., $c' \leq \lambda_{\min}(\Sigma) \leq \lambda_{\max}(\Sigma) \leq C'$ for some $c',C' > 0$, and sub-Gaussian parameter of constant order, i.e. $\sup_{v \in S^{n-1}} \EE \exp(\lambda v\T Z_i) \leq \exp(\lambda^2 \zeta^2/2)$ for $\zeta = O(1)$. Then the minimax rate $\eta^{\star2}$, in the case when $\xi_i \sim N(0,1)$, is given by the entropy equation $\eta^\star = \sup_\eta\{\eta \geq 0: N\eta^2 \leq \log M_K^{\operatorname{loc}}(\eta)\}$, as shown in \cite{prasadan2025characterizingminimaxratenonparametric}. 

Returning to Lemma \ref{Kwidth:bound:lemma}, and plugging in $N^{-1/2}$ for $\sigma$ we obtain the following corollary. 

\begin{corollary}\label{Kwidth:bound:lemma:regression}
Let $K$ be a type-2 convex body with type-2 constant $T_2(K)$, such that $rB_2 \subseteq K \subseteq RB_2$ for some known $r,R$ with $R = O(1)$. Let $\kappa:= \kappa(K)$ be a given scalar which can depend on $K$. Let $\eta^\star := \eta^\star(K)$ be the minimax rate defined by equation $\eta^\star = \sup_\eta\{\eta \geq 0: N\eta^2 \leq \log M_K^{\operatorname{loc}}(\eta)\}$. Assume that $\kappa^{-1/2} \tilde d/T_2(K) \geq 2\eta^\star$ and $\lceil \lceil\log_c(\kappa^{1/2}\cdot T_2(K))\rceil \cdot 4 \eta^{\star 2}\cdot N\rceil \leq N\cdot \tilde d^2/C^2$ for $\tilde d = 2R \geq d$ and where $C > 0$ is some absolute constant. Then for $m = \lceil N\cdot \tilde d^2/C^2\rceil \wedge n$, we have
\begin{align*}
    d_{m}(K) \lesssim \kappa^{-1/2}\cdot \tilde d/c, 
    \end{align*}
\end{corollary}

We now state a lemma analogous to Lemma \ref{sub:gaussian:vector:estimator}. However, the proof is more involved and uses an estimator adapted to the regression setting. We also mention that we are using Theorem \ref{main:SDP:theorem} and Algorithm \ref{algorithm:SDP:Kolmogorov:width} with $N^{-1/2}$ for $\sigma$.

\begin{lemma}\label{sub:gaussian:vector:estimator:regression} Suppose that we have an estimator $\hat \beta_j$ of $\beta$ such that $\|\hat \beta_j - \beta\|_2 \leq \tilde d_j$ happens with high probability, for some $\tilde d_j \in \RR_+$ such that $\tilde d_j \geq (2r) \vee (C/\sqrt{N})$. Assume also $N \gtrsim n$ for some absolute constant. Let $K_{(j)} = K \cap B(0, \tilde d_j/2)$. Obtain $m$ and the matrix $A_j^{\star\star}$ from Corollary \ref{Kwidth:bound:lemma:regression} and Algorithm \ref{algorithm:SDP:Kolmogorov:width} for the set $K_{(j)}$. Let 
\begin{align}\label{strange:optimization:program}
    \beta'_{j + 1} := \argmin_{\nu \in K_{(j)}} \sum_{i \in [N]} ((Y_i - Z_i\T\hat \beta_{j})/2 - Z_i\T A_j^{\star\star}\nu)^2
\end{align}
and define the estimator $\tilde \beta_{j + 1} := A^{\star\star}_j \beta'_{j + 1}$. Then solving program \eqref{strange:optimization:program} can be approximated in polynomial-time and the estimator $\tilde \beta_{j + 1}$ satisfies:
\begin{align*}
    \PP(\|\tilde \beta_{j + 1} -( \beta - \hat \beta_j)/2\|_2 \geq \tilde d_j / L, \|\hat \beta_j - \beta\|_2 \leq \tilde d_j) \leq \exp(-Nc'\tilde d_j^2) + (1 + N^{3/2} \tilde d_j^{3})^{-1}.
\end{align*}
for some absolute constant $L$ which can be made large provided that the constant $C$ is selected to be sufficiently large. 
\end{lemma}

Due to space constraints the proof is relegated to the appendix. Now that we have established Lemma \ref{sub:gaussian:vector:estimator:regression}, all that remains is to state the algorithm for estimating $\beta$. The algorithm is very similar to the one in the mean estimation setting:

\begin{breakablealgorithm}
    \renewcommand{\algorithmicrequire}{\textbf{Input:}}
    \renewcommand{\algorithmicensure}{\textbf{Output:}}
    \caption{Near-optimal Regression Estimation Algorithm}
    \label{algorithm:main:regression}
   \begin{algorithmic}[1]
        \REQUIRE Convex body $K$, $r, R \in \RR_+$ such that $rB_2 \subseteq K \subseteq R B_2$, with $R = O(1)$; observations $(Y_i, Z_i)_{i \in [N]}$, constant $\tilde L = L/(\sqrt{3} + 1)$ where $L$ is from Lemma \ref{sub:gaussian:vector:estimator:regression}; We assume $\tilde L > 2(\sqrt{3} + 1)$
        \ENSURE Estimator $\hat \beta$ of $\beta$ 

        \STATE Set $r = r \wedge (\frac{1}{2} \sqrt{n/N})$
        \STATE Set $\tilde d_1 = 2R$, $K_{(1)} := K$,
        number of iterations $M = \lceil\log_{\tilde L/(2(\sqrt{3} + 1))} R/r\rceil$
        \STATE Set $\hat \beta_1 := 0$
        \FOR{ j = 1, \ldots, M}
            \STATE Find $m = \lceil N \tilde d_j^2/C^2\rceil \wedge n$, and calculate the matrix $A_j^{\star \star}$ for the set $K_{(j)}$
            \STATE Find $\tilde \beta_{j + 1}$, the estimator from Lemma \ref{sub:gaussian:vector:estimator:regression}
            \STATE Set $\bar \beta_{j + 1} = \Pi_{K_{(j)}} \tilde \beta_{j + 1}$\footnote{Here we set $\epsilon = \tilde d_j/L$ in Lemma \ref{weak:projection:lemma} when evaluating the weak projection.}
            \STATE Lift back to $K$: $\hat \beta_{j + 1} = \Pi_K(2 \bar \beta_{j + 1} + \hat \beta_j)$\footnote{Here we set $\epsilon = 2\tilde d_j/\tilde L$ in Lemma \ref{weak:projection:lemma} when evaluating the weak projection.}
            \STATE Set $\tilde d_{j + 1} = 2(\sqrt{3} + 1)\tilde d_j/\tilde L$
            \STATE Set $K_{(j + 1)} := K\cap B_2(0, \tilde d_{j+1}/2)$
            \IF{$\tilde d_{j + 1} \leq (2r) \vee (C/\sqrt{N})$}
                \RETURN $\hat \beta := \hat \beta_{j + 1}$
            \ENDIF
        \ENDFOR
        \RETURN $\hat \beta := \hat \beta_{M + 1}$
    \end{algorithmic} 
\end{breakablealgorithm}

We can now state the main theorem of this section.

\begin{theorem}\label{regression:theorem:main}
Suppose that the assumptions of Section~\ref{regression:section} hold. Set
\begin{align*}
\overline\kappa
:=
\max_{j\in[M]}\kappa(K_{(j)}),
\end{align*}
and suppose that $T_2(K)$ and $\overline\kappa$ are polylogarithmic in
$n,N$, while $\log(R/r)$ is polynomial in $n,N$. Assume also that
$N\gtrsim n$ for a sufficiently large absolute constant. Then there exists
a factor
\begin{align*}
\Upsilon(K)
=
\left(1
    \vee
    \left(
        2\sqrt2\,T_2(K)\overline\kappa^{1/2}
    \right)
    \vee
    \left(
        2C\sqrt{\Lambda_K}
    \right)\right)^2,
\end{align*}
with $\Lambda_K = 1+\left\lceil\log_c\!\left(\overline\kappa^{1/2}\sqrt2\,T_2(K)\right)\right\rceil$, such that the output
$\widehat\beta$ of Algorithm~\ref{algorithm:main:regression} satisfies
\begin{align*}
\EE\|\widehat\beta-\beta\|_2^2
\lesssim
\Upsilon(K)\eta^{\star2}(K).
\end{align*}
In particular, Algorithm~\ref{algorithm:main:regression} is nearly minimax
rate optimal and runs in polynomial time.
\end{theorem}

The proof of Theorem \ref{regression:theorem:main} is deferred to the appendix since it follows a similar strategy to the proof of Theorem \ref{main:result:GSM}. Since ellipsoids are type-2 sets with small type-2 constant, we note that our method provides efficient implementation for linear regression problem under ellipsoidal constraints considered in the recent paper \cite{pathak2024noisy}, provided that our further assumptions on the constraint set, covariates and noise vector are satisfied. Finally we remark how the assumption $R = O(1)$ can be significantly relaxed at almost no algorithmic cost.

\begin{remark}[On the bounded-diameter assumption]
\label{rem:regression:bounded:diameter}\label{R:bounded:assumption}
The assumption $\diam(K)=O(1)$ in this section is mainly a convenient
normalization and can be relaxed by a preliminary localization step.  We
describe one way to do this.  The sample splitting below is useful because it
ensures that, conditional on the pilot estimator, the observations used in the
main procedure still form an independent linear-regression sample.

Split the observations into two subsamples $I_0$ and $I_1$ of comparable
cardinality, say $|I_0|=\lfloor N/2\rfloor$ and
$|I_1|=N-|I_0|$.  Using the first subsample, form the ordinary least-squares
estimator
\begin{align*}
    \widehat\beta_{\rm OLS}
    =
    \left(\sum_{i\in I_0} Z_iZ_i^\top\right)^{-1}
    \sum_{i\in I_0}Z_iY_i.
\end{align*}
Under the assumptions of this section, standard concentration for
sub-Gaussian random matrices implies that, provided $N\gtrsim n$,
\begin{align*}
    \lambda_{\min}\!\left(
        \frac1{|I_0|}\sum_{i\in I_0}Z_iZ_i^\top
    \right)
    \gtrsim 1
\end{align*}
with probability at least $1-\exp(-cN)$.  Conditional on the design,
\begin{align*}
    \widehat\beta_{\rm OLS}-\beta
    =
    \left(\sum_{i\in I_0}Z_iZ_i^\top\right)^{-1}
    \sum_{i\in I_0}Z_i\xi_i.
\end{align*}
Since the $\xi_i$ are independent centered sub-Gaussian random variables with
sub-Gaussian parameter $\gamma=O(1)$, an $\varepsilon$-net argument, together
with the preceding lower bound on the empirical covariance, gives
\begin{align*}
    \PP\!\left(
        \|\widehat\beta_{\rm OLS}-\beta\|_2>C_0\gamma
    \right)
    \leq
    \exp(-c_0N)
\end{align*}
for sufficiently large absolute constants $C_0,c_0>0$.  More precisely, one
has the usual bound
\begin{align*}
    \|\widehat\beta_{\rm OLS}-\beta\|_2
    \lesssim
    \gamma\sqrt{\frac nN}+t
\end{align*}
with probability at least
$1-\exp(-cNt^2/\gamma^2)-\exp(-cN)$, and the displayed constant-radius
bound follows from $N\gtrsim n$ and $\gamma=O(1)$. For more details on the last two statements please refer to Lemma \ref{lem:ols:pilot:bound} in the appendix. Let
\begin{align*}
    \overline\beta
    :=
    \Pi_K\widehat\beta_{\rm OLS},
\end{align*}
where, as elsewhere in the paper, $\Pi_K$ may be taken to be a weak Euclidean
projection with accuracy of order $\gamma$.  Lemma
\ref{weak:projection:lemma} then implies that, on the preceding event,
\begin{align*}
    \|\overline\beta-\beta\|_2
    \leq
    C_1\gamma
\end{align*}
for an absolute constant $C_1$, while
$\overline\beta\in K$.  Define
\begin{align*}
    \theta
    :=
    \frac{\beta-\overline\beta}{2}.
\end{align*}
Since $K$ is origin-symmetric and convex,
$\beta,-\overline\beta\in K$, and therefore
\begin{align*}
    \theta
    =
    \frac{\beta+(-\overline\beta)}2
    \in K.
\end{align*}
Moreover,
\begin{align*}
    \|\theta\|_2
    \leq
    C_1\gamma/2.
\end{align*}
Thus, on the pilot-good event, the new coefficient belongs to the bounded
convex body
\begin{align*}
    K_0
    :=
    K\cap B_2(0,C_2\gamma),
\end{align*}
where $C_2$ is a sufficiently large absolute constant.  Notice that $K_0$ is
again origin-symmetric and convex, and
\begin{align*}
    (r\wedge C_2\gamma)B_2
    \subseteq
    K_0
    \subseteq
    C_2\gamma B_2.
\end{align*}
Furthermore, by the same intersection argument used earlier in the paper,
\begin{align*}
    T_2(K_0)
    \leq
    \sqrt{2}\max\{T_2(K),1\},
\end{align*}
and
\begin{align*}
    \log M_{K_0}^{\operatorname{loc}}(\delta)
    \leq
    \log M_K^{\operatorname{loc}}(\delta).
\end{align*}
Hence the corresponding entropy rate satisfies
$\eta^\star(K_0)\lesssim \eta^\star(K)$, where replacing $N$ by
$|I_1|\asymp N$ changes the entropy fixed point by at most an absolute
constant factor.  Indeed, if
\begin{align*}
    \eta_N^\star
    :=
    \sup\{\eta:N\eta^2\leq
    \log M_K^{\operatorname{loc}}(\eta)\},
\end{align*}
then monotonicity of the local entropy immediately gives
$\eta_{N/2}^\star\leq\sqrt{2}\,\eta_N^\star$.
For the classes considered in Section~\ref{QFM:section}, the QFM assumptions
are also preserved under intersection with a Euclidean ball.

Now use only the second subsample $I_1$ and define
\begin{align*}
    Y_i'
    :=
    \frac{Y_i-Z_i^\top\overline\beta}{2},
    \qquad i\in I_1.
\end{align*}
Conditional on the first subsample,
\begin{align*}
    Y_i'
    =
    Z_i^\top\theta+\frac{\xi_i}{2},
    \qquad i\in I_1,
\end{align*}
and the observations on the right-hand side are independent of the pilot
estimator.  Thus they satisfy exactly the regression assumptions of this
section, but with parameter $\theta\in K_0$, a constraint set of bounded
diameter.  Applying the bounded-diameter procedure to this transformed
problem produces an estimator $\widehat\theta$ satisfying, on the pilot-good
event,
\begin{align*}
    \EE\!\left[
        \|\widehat\theta-\theta\|_2^2
        \,\middle|\,
        (Y_i,Z_i)_{i\in I_0}
    \right]
    \lesssim
    \Upsilon(K_0)\eta^{\star2}(K_0)
    \lesssim
    \Upsilon(K)\eta^{\star2}(K),
\end{align*}
up to the same logarithmic factors as in the main theorem.

Finally, lift the estimator back to the original parameterization by setting
\begin{align*}
    \widetilde\beta
    :=
    2\widehat\theta+\overline\beta
\end{align*}
and weakly projecting $\widetilde\beta$ onto $K$.  Taking the accuracy of this
last weak projection to be, for example,
\begin{align*}
    \epsilon_{\rm proj}
    \lesssim
    r\wedge N^{-1/2},
\end{align*}
does not affect the statistical rate.  Indeed, since $rB_2\subseteq K$,
\begin{align*}
    \eta^\star(K)
    \gtrsim
    r\wedge\sqrt{\frac nN}
    \gtrsim
    r\wedge N^{-1/2},
\end{align*}
while Lemma~\ref{weak:projection:lemma} shows that this final projection changes
the estimation error by at most an absolute constant factor at the relevant
scale.  The resulting estimator is proper, i.e., it belongs to $K$.

The only additional contribution to the expected risk comes from failure of
the OLS localization event.  Since the final estimator belongs to $K$, its
squared error is always at most $d^2$, where $d=\diam(K)$.  Consequently the
preceding reduction yields
\begin{align*}
    \EE\|\widehat\beta-\beta\|_2^2
    \lesssim
    \Upsilon(K)\eta^{\star2}(K)
    +
    d^2\exp(-c_0N).
\end{align*}
Therefore the bounded-diameter assumption may be removed whenever
\begin{equation}
    d^2\exp(-c_0N)
    \lesssim
    \eta^{\star2}(K).
    \label{eq:regression:pilot:diameter:condition}
\end{equation}
A completely explicit sufficient condition follows from
$rB_2\subseteq K$, which implies
\begin{align*}
    \eta^{\star2}(K)
    \gtrsim \eta^{\star2}(r B_2) \asymp
    r^2\wedge\frac nN.
\end{align*}
Thus it suffices to assume
\begin{align*}
    d^2\exp(-c_0N)
    \lesssim
    r^2\wedge\frac nN.
\end{align*}
\end{remark}

\section{Some Remarks}\label{examples:section}

In this section, we give an example of a type of convex bodies that obey all conditions previously required. It has long been known \cite{donoho1990minimax} that for convex Quadratically Convex Orthosymmetric (QCO) sets, the minimax rate is given by an equation involving the Kolmogorov widths. In fact, the truncated series estimator is rate optimal over any QCO set. However, it was not previously known how one can compute a minimax optimal estimator (including the truncated series one) for QCO sets in generality. Here we argue that any orthosymmetric QCO set induces a norm which is of type-2 with small type-2 constant and is exactly 2-convex. 

\begin{lemma}\label{lemma:QCO}
Let $ K \subset \mathbb{R}^n $ be a \emph{quadratically convex orthosymmetric (QCO)} set, i.e.
\begin{enumerate}
    \item $K$ is convex, origin-symmetric, and sign-invariant;
    \item the set $ K^{(2)} := \{x^2 : x \in K\}$ is convex, where squaring is done coordinatewise.
\end{enumerate}
Then the norm $ \rho_K(\cdot) $ induced by $ K $ satisfies a \emph{type-2 inequality} with constant $\sqrt{ C \log n } $, i.e., $T_2(K)\lesssim \sqrt{\log n}$:
\begin{align*}
    \mathbb{E}\rho_K(\sum_{i=1}^m \varepsilon_i x_i)^2
    \;\le\;
    C (\log n) \sum_{i=1}^m \rho_K(x_i)^2,
\end{align*}
for all $x_1,\dots,x_m \in \mathbb{R}^n$, where $\varepsilon_i$ are i.i.d.\ Rademacher signs.
\end{lemma}

We note furthermore that the Minkowski gauges of QCO sets are necessarily 2-convex. To show this define
\begin{align*}
    K^{2} := \{\operatorname{sign}(x)\circ x^2 : x \in K\}.
\end{align*}
Since $K$ is convex, sign-invariant and quadratically convex, $K^2$ is convex and sign-invariant (see Lemma \ref{K2:is:convex} for a short proof).  
By direct inspection,
\begin{align*}
    |x| \in tK \iff x^2 \in t^2 K^2,
\end{align*}
and therefore
\begin{equation*}
    \rho_K(x) = \sqrt{\rho_{K^2}(x^2)}.
\end{equation*}
Thus, by the triangle inequality:
\begin{align*}
    \rho_K^2((\sum_{i \in [m]} |x_i|^2)^{1/2}) = \rho_{K^2}(\sum_{i \in [m]} |x_i|^2) \leq \sum_{i \in [m]} \rho_{K^2}(|x_i|^2) = \sum_{i \in [m]}\rho_K^2(x_i)
\end{align*}

It follows immediately that if $K \subset \RR^m$ is a QCO set (with $m$ being ``small''), we can implement our algorithm over bounded, well-balanced sets $K' = \{x \in \RR^n : \rho_K(Ax) \leq 1\}$ for some full rank matrix $A \in \RR^{m \times n}$ (provided that we can evaluate $\rho_K$). Clearly in the special case when $A = I_{n}$, this implies that we can implement our algorithm for any well-balanced QCO set with oracle evaluating its Minkowski gauge. Examples of QCO include ellipsoids, hyperrectangles and sets of the type $\{ x \in \RR^n : \sum_{i \in [n]} \theta_i\psi(x_i^2) \leq 1\}$ for some convex function $\psi$ and constants $\theta_i \in \RR_+$.

It is easy to see that not all bodies with small type-2 constant are QCO sets. We attach the following proposition which verifies this fact. 

\begin{proposition}\label{not:all:type2:sets:are:QCO}
Let $n \ge 2$ and define, for $x \in \mathbb{R}^n$,
\begin{align*}
\|x\| = \Big( \|x\|_2^2 + \frac{1}{n}\|x\|_1^2 \Big)^{1/2}.
\end{align*}
Then $\|\cdot\|$ is a norm which is convex, sign--symmetric, permutation--invariant,
and has type--2 constant bounded by an absolute constant.
However, its unit ball
\begin{align*}
K = \bigl\{x \in \mathbb{R}^n : \|x\| \le 1\bigr\}
\end{align*}
is \emph{not quadratically convex}.
\end{proposition}

Hence our algorithm can be implemented for well-balanced sets $K' = \{x \in \RR^n : \|A x\| \leq 1\}$ where $A \in \RR^{n \times n}$ is a full rank matrix. We conclude that our algorithm applies not only to transformed QCO sets, but also to a broader class of constraint sets.

\section{Discussion}\label{discussion:section}

To the best of our knowledge, this paper proposed the first computationally efficient algorithm for constrained mean estimation in a Gaussian sequence model setting, heavy-tailed robust mean estimation and linear regression under type-2 convex constraints, where the constraint set satisfies further regularity assumptions as outlined in Sections \ref{background:cond:algo:section} and \ref{QFM:section}. There is one central open question that remains: can one develop near optimal polynomial-time algorithms for other (convex) constraints. Notably our algorithm does not work for $\ell_p$ balls for $1 \leq p < 2$, for which efficient near-minimax optimal soft and hard thresholding algorithms do exist \cite{johnstone2011gaussian}. We further had to assume that the set $K$ is origin-symmetric. The hypothetical efficient algorithm will likely have to use a fundamentally different approach as we are no longer guaranteed that the Kolmogorov widths decay as in the origin-symmetric type-2 sets considered in this paper.

Another interesting open question is whether our algorithm is meaningful for misspecified scenarios, i.e., when $\mu$ or $\beta$ do not belong to the set $K$.

Yet another question is as follows. In the work of \cite{neykov2022minimax} when the constraint set is convex, one can show that the estimator is adaptive to the location of the true mean $\mu \in K$. However, we suspect that the estimator provided in this paper is not adaptive. It will be interesting to look into this rigorously: either establish that the estimator proposed here is adaptive, or refute it and provide an adaptive estimator.

Another avenue for future research is to try to remove the condition $d^2/\sigma^2 \lesssim \exp(N/C_0)/N$ in the robust algorithm. We believe this assumption is not necessary, but we will likely need to develop a new algorithm which converges without the need of the pilot estimator based on the work of \cite{depersin2022}.

In the regression setting we had to assume the dimension is smaller than the sample size up to an absolute constant. It will be interesting if one can remove this condition. We crucially used in the proof of Lemma \ref{sub:gaussian:vector:estimator:regression}, to show that the eigenvalues of the sample covariance matrix are upper and lower bounded, and we do not know presently how to remove it.

Finally, our mean estimation procedures require knowledge of $\sigma$, or an accurate upper bound. Whether one can develop near optimal adaptive computationally tractable algorithms is an open problem.

\section{Acknowledgments}

The author is grateful to Sourav Chatterjee for encouraging him to think deeply about this problem. He also benefited from valuable discussions with Andrea Montanari, Siva Balakrishnan, and Ilias Diakonikolas, as well as from exploratory conversations with ChatGPT. Finally, the author would like to thank Vijay Bhattiprolu for a helpful email exchange, where Vijay clarified several points of the central paper to our work \cite{bhattiprolu2021framework}. Furthermore, the author thanks Daniel Dadush for clarifying a point in Theorem 2.5.9 in his PhD thesis.

\bibliographystyle{abbrv}
\bibliography{polytime}

\newpage 

\appendix

\section{External results used in the proofs}\label{appendix:external:results}

For the reader's convenience, in this appendix we collect the external results that are
used as black boxes in the paper.  We state them in the notation used throughout the
manuscript and, when the statement in the cited source uses a slightly different
normalization, we explicitly record the translation.  We only include the portions of
the cited results that are used in our arguments.

\subsection{Quadratic-form maximization results of Bhattiprolu--Lee--Naor}
\label{appendix:BLN:results}

We first collect the results from \cite{bhattiprolu2021framework}.  For an
origin-symmetric convex body $K\subset\RR^n$ and a symmetric matrix
$X\in\Smat^n$, write
\begin{align}
    Q_K(X):=\max_{p\in K}p\T Xp.
    \label{eq:appendix:QK:def}
\end{align}
Thus, a $\kappa$-approximate \emph{search} oracle for quadratic-form maximization
(QFM) over $K$ is an oracle that, on input $X$, returns a point
$\cO_K(X)\in K$ satisfying
\begin{align}
    \cO_K(X)\T X\cO_K(X)
    \geq \kappa^{-1}Q_K(X).
    \label{eq:appendix:QFM:oracle}
\end{align}
This is the convention used in \eqref{oracle:equation} in the main text.

Let $g\sim N(0,I_n)$.  Recall the upper and lower covariance regions
\begin{align}
    \cU(K)
    &:=%
    \left\{W\in\Smat^n:\ W\geq0,\ 
    \EE\rho_K^2(W^{1/2}g)\leq1\right\},
    \label{eq:appendix:upper:covariance}\\
    \cL(K)
    &:=%
    \left\{W\in\Smat^n:\ W\geq0,\ 
    \EE\rho_{K^\circ}^2(W^{1/2}g)\geq1\right\}.
    \label{eq:appendix:lower:covariance}
\end{align}
We use the notion of an approximate separation oracle from
\cite{bhattiprolu2021framework}.  Concretely, if a star-shaped set $S$ satisfies
$\operatorname{conv}(S)\subseteq \alpha S$, then an $\alpha$-approximate separation
oracle, on input $x$, either returns ``Inside,'' in which case $x\in\alpha S$, or
returns an affine hyperplane separating $x$ from $S$.  For an outward-closed set
$T$ satisfying $\operatorname{conv}(T)\subseteq \alpha^{-1}T$, the analogous
convention is that ``Inside'' certifies $x\in\alpha^{-1}T$, while a returned
hyperplane separates $x$ from $T$.  These are Definitions~1.6 and~1.11 of
\cite{bhattiprolu2021framework}, respectively.

The source \cite{bhattiprolu2021framework} formulates several of the results below
using the Gaussian type-$2$ constant $\widetilde T_2$.  Our paper uses the
Rademacher type-$2$ constant $T_2(K)$ from \eqref{type2:constant}.  By the standard
Gaussian--Rademacher comparison (equivalently, Kahane--Khintchine), these two
constants are equivalent up to an absolute multiplicative constant.  Accordingly,
when translating a bound involving $\widetilde T_2$ into our convention, we write
$\lesssim T_2(K)$ and absorb only a universal constant.

\begin{remark}\label{complexity:of:r:and:R:rem}
Although the results of \cite{bhattiprolu2021framework} state the running time in terms of $\log R$ and $\log(1/r)$ for a body satisfying
\[
rB_2\subseteq K\subseteq RB_2,
\]
the geometric dependence may equivalently be expressed in terms of $\log(R/r)$. Indeed, consider the rescaled body
\[
\widetilde K:=r^{-1}K,
\]
for which
\[
B_2\subseteq \widetilde K\subseteq (R/r)B_2.
\]
Moreover,
\[
\rho_{\widetilde K}(x)=r\rho_K(x),
\]
so an evaluation oracle for the Minkowski gauge of $K$ immediately yields one for $\widetilde K$. Quadratic-form maximization is also invariant under this normalization up to the deterministic scaling
\[
\max_{z\in\widetilde K} z^\top Xz
=
r^{-2}\max_{p\in K}p^\top Xp,
\]
and hence approximation factors are unchanged after mapping an output $z\in\widetilde K$ back to $rz\in K$. Finally, the type-$2$ constant is invariant under scalar dilations. Thus one may apply the corresponding results to $\widetilde K$ with inner radius $1$ and outer radius $R/r$, so that the dependence on the balancing radii is through $\log(R/r)$ alone, apart from the ordinary bit complexity of the numerical input parameters. This observation was also confirmed by Vijay Bhattiprolu in personal communication.
\end{remark}

\begin{proposition}[\cite{bhattiprolu2021framework}, Proposition~4.2]
\label{prop:appendix:BLN:4.2}
Suppose that $K$ is well-balanced, so that the radii $r,R>0$ are known and
\[
    rB_2\subseteq K\subseteq RB_2.
\]
Assume that an $\alpha$-approximate separation oracle for $\cU(K)$ is available.
Then, for every symmetric matrix $X\in\Smat^n$, there is a randomized algorithm
running in
\[
    \operatorname{poly}\!\left(n,\log R,\log(1/r),\operatorname{bit}(X)\right)
\]
time which, with probability at least $1-2^{-\Omega(n)}$, returns $p\in K$ such
that
\begin{align}
    p\T Xp
    \geq
    \frac{1}{(1+o(1))\alpha}\,Q_K(X).
    \label{eq:appendix:BLN:4.2}
\end{align}
In particular, the result is a search guarantee: the algorithm returns an
approximately optimal vector, not only an approximation to the optimal value.
\end{proposition}

\begin{proposition}[\cite{bhattiprolu2021framework}, Proposition~4.4]
\label{prop:appendix:BLN:4.4}
Suppose that $rB_2\subseteq K\subseteq RB_2$ with known $r,R>0$, and assume that
an $\alpha$-approximate separation oracle for $\cL(K)$ is available.  Then, for
every symmetric psd matrix $X\in\Smat^n$, there is a randomized algorithm running
in
\[
    \operatorname{poly}\!\left(n,\log R,\log(1/r),\operatorname{bit}(X)\right)
\]
time which, with probability at least $1-2^{-\Omega(n)}$, returns $p\in K$ such
that
\begin{align}
    p\T Xp
    \geq
    \frac{1}{(1+o(1))\alpha}\,Q_K(X).
    \label{eq:appendix:BLN:4.4}
\end{align}
\end{proposition}

\begin{corollary}[\cite{bhattiprolu2021framework}, Corollary~4.12(1)]
\label{cor:appendix:BLN:4.12}
Suppose that $rB_2\subseteq K\subseteq RB_2$ with known $r,R>0$.  Assume that
there is an exactly computable convex function
\[
    f:\{W\in\Smat^n:W\geq0\}\to\RR
\]
and constants $\alpha_1,\alpha_2\geq1$ such that, for every $W\geq0$,
\begin{align}
    \alpha_1^{-1}f(W)
    \leq
    \EE\rho_K^2(W^{1/2}g)
    \leq
    \alpha_2 f(W).
    \label{eq:appendix:BLN:4.12:assumption}
\end{align}
Then for every symmetric $X\in\Smat^n$ one can, in
\[
    \operatorname{poly}\!\left(n,\log R,\log(1/r),\operatorname{bit}(X)\right)
\]
time, return with probability at least $1-2^{-\Omega(n)}$ a vector $p\in K$ such
that
\begin{align}
    p\T Xp
    \gtrsim
    \frac{1}{\alpha_1\alpha_2T_2(K)^2}\,Q_K(X).
    \label{eq:appendix:BLN:4.12:conclusion}
\end{align}
Equivalently, under our Rademacher type-$2$ convention the approximation factor is
$\lesssim \alpha_1\alpha_2T_2(K)^2$.  If one uses the Gaussian type-$2$ constant
$\widetilde T_2$ of \cite{bhattiprolu2021framework}, the source states the factor
exactly as $\alpha_1\alpha_2\widetilde T_2(K)^2$.
\end{corollary}

\begin{proposition}[\cite{bhattiprolu2021framework}, Proposition~4.13]
\label{prop:appendix:BLN:4.13}
Suppose that $rB_2\subseteq K\subseteq RB_2$ with known $r,R>0$, and suppose
that a $\kappa$-approximate QFM search oracle over $K$, in the sense of
\eqref{eq:appendix:QFM:oracle}, is available.  Then one can construct an approximate
separation oracle for $\cU(K)$ with approximation factor
\begin{align}
    \alpha
    \lesssim
    \kappa T_2(K)^2,
    \label{eq:appendix:BLN:4.13:factor}
\end{align}
whose running time is polynomial in
$n,\log R,\log(1/r)$ and the bit complexity of the input matrix.  More precisely,
in the Gaussian type convention of \cite{bhattiprolu2021framework} the factor is
$(1+o(1))\kappa\widetilde T_2(K)^2$.
\end{proposition}

\begin{proposition}[\cite{bhattiprolu2021framework}, Proposition~6.5]
\label{prop:appendix:BLN:6.5}
Let $K_1,K_2\subset\RR^n$ be origin-symmetric, well-balanced convex bodies, and
assume that both QFM problems admit $\alpha$-approximate search oracles.  Then the
intersection $K_1\cap K_2$, whose Minkowski gauge is
\[
    \rho_{K_1\cap K_2}(x)=\rho_{K_1}(x)\vee\rho_{K_2}(x),
\]
admits a randomized polynomial-time QFM search algorithm with approximation factor
\begin{align}
    \kappa(K_1\cap K_2)
    \lesssim
    2\alpha\max\left\{T_2(K_1)^2,T_2(K_2)^2\right\},
    \label{eq:appendix:BLN:6.5}
\end{align}
and success probability at least $1-2^{-\Omega(n)}$.  If the two available QFM
oracles have factors $\kappa_1$ and $\kappa_2$, the same conclusion follows after
setting $\alpha=\kappa_1\vee\kappa_2$.

In particular, for $K_2=cB_2$ we have $T_2(cB_2)=1$ and exact QFM over $cB_2$.
Thus, if $K$ admits a $\kappa$-approximate QFM search oracle, then
\begin{align}
    \kappa(K\cap cB_2)
    \lesssim
    2\kappa T_2(K)^2.
    \label{eq:appendix:intersection:ball}
\end{align}
Under the definition \eqref{eq:appendix:upper:covariance}, one also has the explicit
identity
\begin{align}
    \cU(cB_2)
    =
    \left\{W\in\Smat^n:\ W\geq0,\ \tr(W)\leq c^2\right\}.
    \label{eq:appendix:U:euclidean:ball}
\end{align}
\end{proposition}

We next state the Banach-lattice results used to construct the separation oracles
appearing above.  Suppose that $\rho_K$ is sign-invariant.  For $1\leq p<\infty$,
define the $p$-convexity constant $M^{(p)}(K)$ as the smallest $C$ such that
\begin{align}
    \rho_K\!\left(\left(\sum_i |x_i|^p\right)^{1/p}\right)
    \leq
    C\left(\sum_i\rho_K(x_i)^p\right)^{1/p}
    \label{eq:appendix:p:convexity}
\end{align}
for every finite collection $(x_i)$.  Similarly, define the $q$-concavity constant
$M_{(q)}(K)$ as the smallest $C$ such that
\begin{align}
    \rho_K\!\left(\left(\sum_i |x_i|^q\right)^{1/q}\right)
    \geq
    C^{-1}\left(\sum_i\rho_K(x_i)^q\right)^{1/q}.
    \label{eq:appendix:q:concavity}
\end{align}
All powers and absolute values inside the vector in
\eqref{eq:appendix:p:convexity}--\eqref{eq:appendix:q:concavity} are applied
coordinatewise.  The norm is exactly $p$-convex when $M^{(p)}(K)=1$.

\begin{theorem}[\cite{bhattiprolu2021framework}, Theorem~7.4]
\label{thm:appendix:BLN:7.4}
Let $\rho_K$ be sign-invariant, let $1\leq p\leq2\leq q<\infty$, let
$x_1,\ldots,x_m\in\RR^n$, and let $g_1,\ldots,g_m$ be independent standard
Gaussian random variables.  Put
\[
    \gamma_s:=\left(\EE|g_1|^s\right)^{1/s}.
\]
Then
\begin{align}
    \frac{\gamma_p}{M^{(p)}(K)}
    \rho_K\!\left(\left(\sum_{i=1}^m|x_i|^2\right)^{1/2}\right)
    &\leq
    \left(\EE\rho_K^2\!\left(\sum_{i=1}^m g_i x_i\right)\right)^{1/2}
    \nonumber\\
    &\leq
    \gamma_q M_{(q)}(K)
    \rho_K\!\left(\left(\sum_{i=1}^m|x_i|^2\right)^{1/2}\right).
    \label{eq:appendix:BLN:7.4}
\end{align}
In particular, if $\rho_K$ is exactly $2$-convex, then the constant in the left
inequality can be taken to be $1$.
\end{theorem}

\begin{theorem}[\cite{bhattiprolu2021framework}, Theorem~7.6]
\label{thm:appendix:BLN:7.6}
Suppose that $\rho_K$ is an exactly computable, sign-invariant, exactly
$2$-convex norm and that $K$ is well-balanced.  Then, for every symmetric matrix
$X\in\Smat^n$, there is a randomized polynomial-time QFM search algorithm over $K$
which succeeds with probability at least $1-2^{-\Omega(n)}$ and has approximation
factor
\begin{align}
    \kappa(K)
    =
    \inf_{2\leq q<\infty}M_{(q)}(K)^2\gamma_q^2.
    \label{eq:appendix:BLN:7.6}
\end{align}
In particular, for a sign-invariant norm on $\RR^{n_A}$, one may take
$q\asymp\log n_A$; using $M_{(\log n_A)}(K)\leq e$ gives the logarithmic-factor
bound used in Section~\ref{QFM:section}.
\end{theorem}

\begin{theorem}[\cite{bhattiprolu2021framework}, Theorem~7.12]
\label{thm:appendix:BLN:7.12}
Suppose that $\rho_K$ is a symmetric norm, i.e., it is invariant under arbitrary
coordinate permutations and sign changes, and assume that $\rho_K(x)$ can be
evaluated in polynomial time.  Then, for every symmetric matrix $X\in\Smat^n$,
there is a randomized polynomial-time QFM search algorithm which succeeds with
probability at least $1-2^{-\Omega(n)}$ and whose approximation factor, in the
type-$2$ convention of this paper, satisfies
\begin{align}
    \kappa(K)
    \lesssim
    T_2(K)^6\log^4\!\bigl(eT_2(K)\bigr).
    \label{eq:appendix:BLN:7.12}
\end{align}
The source gives the more refined intermediate bound in terms of the dual
cotype-$2$ constant and then bounds it by the right-hand side of
\eqref{eq:appendix:BLN:7.12}.
\end{theorem}

\begin{corollary}[Linear images of symmetric norm balls]
\label{cor:appendix:BLN:symmetric:linear:image}
Let $\|\cdot\|$ be a symmetric norm on $\RR^n$, let $A\in\RR^{n\times n}$ be
full rank, and let
\[
    K=\{x\in\RR^n:\ \|Ax\|\leq1\}.
\]
Then Theorem~\ref{thm:appendix:BLN:7.12} applies after the change of variables
$z=Ax$.  Indeed, for every symmetric matrix $X$,
\begin{align}
    Q_K(X)
    &=
    \max_{\|z\|\leq1}
    z\T A^{-\top}XA^{-1}z.
    \label{eq:appendix:symmetric:linear:image}
\end{align}
Moreover, the normed spaces $(\RR^n,\rho_K)$ and
$(\RR^n,\|\cdot\|)$ are isometric via $A$, and hence their type-$2$ constants
agree.  Consequently one obtains a QFM search oracle over $K$ with
\begin{align}
    \kappa(K)
    \lesssim
    T_2(K)^6\log^4\!\bigl(eT_2(K)\bigr),
    \label{eq:appendix:symmetric:linear:image:factor}
\end{align}
with polynomial running time (including the bit complexity of $A$ and $X$).
\end{corollary}

\begin{remark}[The composition used for symmetric norms]
\label{rem:appendix:BLN:composition}
The chain used in the proof of Theorem~\ref{symmmetric:norm:thm} is now explicit.
Theorem~\ref{thm:appendix:BLN:7.12} first gives
\[
    \kappa(K)
    \lesssim
    T_2(K)^6\log^4\!\bigl(eT_2(K)\bigr).
\]
Proposition~\ref{prop:appendix:BLN:4.13} then gives an approximate separation
oracle for $\cU(K)$ with factor
\[
    \alpha
    \lesssim
    T_2(K)^8\log^4\!\bigl(eT_2(K)\bigr).
\]
Finally, Proposition~\ref{prop:appendix:BLN:6.5}, applied to $K$ and $cB_2$,
yields
\begin{align}
    \kappa(K\cap cB_2)
    \lesssim
    2T_2(K)^{10}\log^4\!\bigl(eT_2(K)\bigr),
    \label{eq:appendix:BLN:composition}
\end{align}
which is precisely the approximation factor used in
Theorem~\ref{symmmetric:norm:thm}.
\end{remark}

\paragraph{Remark on numbering.}
In the version of \cite{bhattiprolu2021framework} used above, item~4.12 is a
\emph{corollary}, not a proposition.  We therefore refer to it as
Corollary~4.12 throughout this appendix.

\subsection{Weak convex optimization and weak Euclidean projection}
\label{appendix:dadush:results}

We next record the convex-optimization result from \cite{dadush2012integer} that
we use to implement weak Euclidean projections and the convex minimizations in the
regression section.

\begin{theorem}[\cite{dadush2012integer}, Theorem~2.5.9]
\label{thm:appendix:dadush:2.5.9}
Let $K\subset\RR^n$ be a convex body for which there are known $a_0\in\RR^n$ and
$r,R>0$ such that
\[
    a_0+rB_2\subseteq K\subseteq a_0+RB_2,
\]
and suppose that a weak membership oracle for $K$ is available.  Let
$f:\RR^n\to\RR$ be a convex $L$-Lipschitz function.  Assume that $f$ is available
through an evaluation oracle which, given $x\in\QQ^n$ and $\delta>0$, returns a
rational number $t$ satisfying
\[
    |t-f(x)|\leq\delta.
\]
Then, for every $\epsilon>0$, one can compute in polynomial time a rational number
$w$ and a point $p\in K$ such that
\begin{align}
    w-\epsilon
    \leq
    \min_{x\in K}f(x)
    \leq
    f(p)
    \leq
    w.
    \label{eq:appendix:dadush:2.5.9}
\end{align}
The running time is polynomial in the encoding lengths of the input parameters
and the corresponding oracle accuracies.
\end{theorem}

Since our sets are origin-symmetric and well-balanced, we take $a_0=0$.  Moreover,
our oracle for the Minkowski gauge $\rho_K$ gives a weak membership oracle for $K$;
in particular, membership is determined by the condition $\rho_K(x)\leq1$ up to the
requested numerical tolerance.

The following is the specialization used repeatedly in the paper.

\begin{corollary}[Weak Euclidean projection]
\label{cor:appendix:weak:projection}
Let $K$ satisfy the assumptions above.  Given $z\in\RR^n$ and $\epsilon>0$, one
can compute in polynomial time a rational number $w$ and $p\in K$ such that
\begin{align}
    w-\epsilon
    \leq
    \|z-\overline\Pi_K(z)\|_2
    \leq
    \|z-p\|_2
    \leq
    w,
    \label{eq:appendix:weak:projection:first}
\end{align}
where $\overline\Pi_K(z)$ denotes the exact Euclidean projection of $z$ onto $K$.
Furthermore, for every $\nu\in K$,
\begin{align}
    \|p-\nu\|_2
    \leq
    \|z-\nu\|_2
    +
    \sqrt{\epsilon^2+2\epsilon\|z-\nu\|_2}.
    \label{eq:appendix:weak:projection:second}
\end{align}
\end{corollary}

\begin{proof}
Apply Theorem~\ref{thm:appendix:dadush:2.5.9} to
$f(x)=\|z-x\|_2$, which is convex and $1$-Lipschitz.  This gives
\eqref{eq:appendix:weak:projection:first}.  Put
$z_0=\overline\Pi_K(z)$.  The projection inequality for convex sets implies
\[
    \|p-z_0\|_2^2+\|z-z_0\|_2^2
    \leq
    \|z-p\|_2^2.
\]
Using $\|z-p\|_2\leq\|z-z_0\|_2+\epsilon$ gives
\[
    \|p-z_0\|_2
    \leq
    \sqrt{\epsilon^2+2\epsilon\|z-z_0\|_2}
    \leq
    \sqrt{\epsilon^2+2\epsilon\|z-\nu\|_2}.
\]
Also, by the projection inequality,
$\|z_0-\nu\|_2\leq\|z-\nu\|_2$.  The triangle inequality therefore yields
\eqref{eq:appendix:weak:projection:second}.
\end{proof}

\subsection{The robust mean estimator of Depersin--Lecue}
\label{appendix:depersin:results}

We state the result in the notation of Section~\ref{robust:section}.  The theorem
appears as Theorem~2.1 in the published version of \cite{depersin2022} (and is the
result referred to as Theorem~2 in the main text).

\begin{theorem}[\cite{depersin2022}, Theorem~2.1]
\label{thm:appendix:depersin:2}
Let $\widetilde Y_1,\ldots,\widetilde Y_N\in\RR^n$ be independent random vectors
with common mean $\mu$ and covariance matrix $\Sigma$.  An adversary is allowed to
replace arbitrarily the observations indexed by a set
$\mathcal I_{\mathrm{out}}\subseteq[N]$, producing the observed vectors
$Y_1,\ldots,Y_N$.  No moment assumption beyond the existence of $\Sigma$ is
required.

Let $k\in\{1,\ldots,N\}$ be the number of equal-sized blocks in the
Depersin--Lecue procedure and assume
\begin{align}
    k\geq300|\mathcal I_{\mathrm{out}}|.
    \label{eq:appendix:depersin:outlier:condition}
\end{align}
Let $u\in\NN$ be the parameter of the covering SDP used at each descent step.
Then their algorithm runs in time
\begin{align}
    \widetilde O\!\left(Nn+kun\right)
    \label{eq:appendix:depersin:runtime}
\end{align}
and, with probability at least
\begin{align}
    1-\exp(-k/180000)-10^{-u},
    \label{eq:appendix:depersin:probability}
\end{align}
returns an estimator, which we denote by $f((Y_i)_{i\in[N]})$, satisfying
\begin{align}
    \bigl\|f((Y_i)_{i\in[N]})-\mu\bigr\|_2
    \leq
    808\left(
        1200\sqrt{\frac{\tr(\Sigma)}{N}}
        +
        \sqrt{\frac{1200\|\Sigma\|_{\operatorname{op}}k}{N}}
    \right).
    \label{eq:appendix:depersin:exact:bound}
\end{align}
\end{theorem}

The form used in our proofs is the immediate universal-constant consequence
\begin{align}
    \bigl\|f((Y_i)_{i\in[N]})-\mu\bigr\|_2
    \lesssim
    \sqrt{\frac{\tr(\Sigma)}{N}}
    +
    \sqrt{\frac{\|\Sigma\|_{\operatorname{op}}k}{N}}.
    \label{eq:appendix:depersin:simplified}
\end{align}
For example, under the standing assumption
$\|\Sigma\|_{\operatorname{op}}\leq\sigma^2$, we also have
$\tr(\Sigma)\leq n\sigma^2$.  Taking $k=u=N$ and assuming the contamination
fraction is smaller than a sufficiently small absolute constant therefore gives
\begin{align}
    \bigl\|f((Y_i)_{i\in[N]})-\mu\bigr\|_2
    \lesssim
    \sqrt{n/N}\,\sigma+\sigma
    \lesssim\sigma
    \label{eq:appendix:depersin:pilot}
\end{align}
with probability at least $1-\exp(-N/C_0)$ for an absolute constant $C_0>0$.
This is the pilot-estimator guarantee used in
\eqref{guarantee:depersin:1st}.  Likewise, applying
\eqref{eq:appendix:depersin:simplified} after the linear transformation
$A^{\star\star}$ gives the bound used in
Lemma~\ref{sub:gaussian:vector:estimator:1:robust}.

\subsection{Carl's inequality}
\label{appendix:carl:result}

We next record Theorem~5 of \cite{carl1985inequalities} and its specialization to
the identity embedding associated with $K$.  For a bounded linear operator
$S:E\to F$, let $s_j(S)$ denote either its $j$-th Kolmogorov number or its $j$-th
Gelfand number, and let
\begin{align}
    g_m(S):=\inf_{k\geq1}k^{1/m}\varepsilon_k(S),
    \label{eq:appendix:carl:g}
\end{align}
where $\varepsilon_k(S)$ is the non-dyadic entropy number, i.e., the smallest
radius for which $S(B_E)$ can be covered by $k$ translates of that radius times
$B_F$.

\begin{theorem}[\cite{carl1985inequalities}, Theorem~5, real case]
\label{thm:appendix:carl:5}
Let $S:E\to F$ be a bounded linear operator between real Banach spaces.  Suppose
that $E$ has Rademacher type $p$ with type constant $\tau_p(E)$ and $F^*$ has
Rademacher type $q$ with type constant $\tau_q(F^*)$.  Then, for every
$m\geq1$,
\begin{align}
    \left(\prod_{j=1}^m s_j(S)\right)^{1/m}
    \leq
    \rho\,\tau_p(E)\tau_q(F^*)\,
    m^{-1+1/p+1/q}g_m(S),
    \label{eq:appendix:carl:operator}
\end{align}
where $\rho$ is an absolute numerical constant in the real case.
\end{theorem}

We now take
\[
    E=(\RR^n,\rho_K),
    \qquad
    F=(\RR^n,\|\cdot\|_2),
    \qquad
    S=I.
\]
Then $S(B_E)=K$, $F^*=F$ is Hilbert and has type-$2$ constant one, and
$\varepsilon_k(S)=e_k(K)$ under Definition~\ref{def_entropy_numbers}.  Taking
$p=q=2$ in Theorem~\ref{thm:appendix:carl:5}, and using the equivalence between
Carl's first-moment type-$2$ convention and the second-moment convention
\eqref{type2:constant}, gives
\begin{align}
    \left(\prod_{j=1}^m d_j(K)\right)^{1/m}
    \lesssim
    T_2(K)\inf_{k\geq1}k^{1/m}e_k(K).
    \label{eq:appendix:carl:specialized}
\end{align}
Carl's operator-theoretic Kolmogorov numbers use subspaces of dimension strictly
less than the index, whereas Definition~\ref{def_Kolmogorov_width} uses a
$k$-dimensional subspace.  The resulting shift by one is absorbed by the
monotonicity of the widths.  Equation~\eqref{eq:appendix:carl:specialized} is
therefore exactly the form of Carl's inequality used in \eqref{carls:ineq}.

\subsection{The global-to-local packing inequality of Yang--Barron}
\label{appendix:yang:barron}

Finally, we record the elementary packing inequality used in the proof of
Lemma~\ref{Kwidth:bound:lemma}.  Recall that $M(\eta,K)$ denotes the maximal
cardinality of an $\eta$-packing of $K$, while for the fixed constant $c>1$ from
Definition~\ref{definition:local:metric:entropy},
\[
    M_K^{\operatorname{loc}}(\eta)
    =
    \sup_{\nu\in K}
    M\!\left(\eta/c,B_2(\nu,\eta)\cap K\right).
\]
The original statement in \cite{yang1999information} is written with the constant
$2$; the same proof gives the version with an arbitrary fixed $c>1$ used here.

\begin{lemma}[Yang--Barron global-to-local packing inequality]
\label{lem:appendix:yang:barron}
For every $\eta>0$,
\begin{align}
    \log M(\eta/c,K)-\log M(\eta,K)
    \leq
    \log M_K^{\operatorname{loc}}(\eta)
    \leq
    \log M(\eta/c,K).
    \label{eq:appendix:yang:barron}
\end{align}
In particular, for every integer $j\geq0$ and every $\Upsilon>0$,
\begin{align}
    \log M(\Upsilon c^{j-1},K)
    -
    \log M(\Upsilon c^j,K)
    \leq
    \log M_K^{\operatorname{loc}}(\Upsilon c^j),
    \label{eq:appendix:yang:barron:used}
\end{align}
which is the form used in the proof of Lemma~\ref{Kwidth:bound:lemma}.
\end{lemma}

\begin{proof}
Let $\{a_1,\ldots,a_M\}$ be a maximal $\eta$-packing of $K$, so that
$M=M(\eta,K)$.  Maximality implies that the balls
$B_2(a_i,\eta)$ cover $K$.  Next let $\{b_1,\ldots,b_L\}$ be an
$\eta/c$-packing of $K$ of maximal cardinality, so that
$L=M(\eta/c,K)$.  Assign every $b_\ell$ to an index $i$ for which
$b_\ell\in B_2(a_i,\eta)$.  By the pigeonhole principle, for at least one $i$,
$B_2(a_i,\eta)\cap K$ contains at least $L/M$ of the points $b_\ell$.  These
points remain $\eta/c$-separated, and therefore
\[
    M_K^{\operatorname{loc}}(\eta)
    \geq
    \frac{M(\eta/c,K)}{M(\eta,K)}.
\]
Taking logarithms gives the first inequality in
\eqref{eq:appendix:yang:barron}.  The second inequality is immediate because
$B_2(\nu,\eta)\cap K\subseteq K$ for every $\nu\in K$.
\end{proof}

\section{Proofs of Section \ref{background:cond:algo:section}}

\begin{proof}[Proof of Lemma \ref{thm:capped-simplex}]
We form the Lagrangian
\begin{align*}
\mathcal{L}(w,\theta,\alpha,\beta)
= \frac{1}{2}\|w - v\|_2^2
+ \theta\Big(\sum_{i=1}^n w_i - k\Big)
- \sum_{i=1}^n \alpha_i w_i
+ \sum_{i=1}^n \beta_i (w_i - 1),
\end{align*}
where $\alpha_i, \beta_i \ge 0$ are the multipliers for the box constraints
$w_i \ge 0$ and $w_i \le 1$, respectively.

Stationarity gives
\begin{align*}
0 = \nabla_w \mathcal{L} = w - v + \theta \mathbf{1} - \alpha + \beta,
\quad\text{i.e.}\quad
w_i = v_i - \theta + \alpha_i - \beta_i.
\end{align*}
The complementary slackness conditions are:
\begin{align*}
\alpha_i w_i = 0, \qquad \beta_i (w_i - 1) = 0.
\end{align*}

Hence, for each coordinate $i$:
\begin{itemize}
    \item If $0 < w_i < 1$, then $\alpha_i = \beta_i = 0$ and $w_i = v_i - \theta$.
    \item If $w_i = 0$, then $\alpha_i \ge 0$ and $v_i - \theta \le 0$.
    \item If $w_i = 1$, then $\beta_i \ge 0$ and $v_i - \theta \ge 1$.
\end{itemize}
Therefore,
\begin{align*}
w_i = \min\{1, \max\{ 0, v_i - \theta^\star \}\}.
\end{align*}

The equality constraint $\sum_i w_i = k$ then becomes
\begin{align*}
f(\theta) = \sum_{i=1}^n \min\{1, \max\{ 0, v_i - \theta \}\} = k.
\end{align*}
The function $f(\theta)$ is continuous, nonincreasing, and piecewise linear.
As $f(\min(v)-1) = n$ and $f(\max(v)) = 0$, by the intermediate value theorem
there exists a $\theta^\star$ such that $f(\theta^\star) = k$. Although $\theta^\star$ need not be unique, the resulting projected vector is unique because the Euclidean projection problem has a strictly convex objective over a convex feasible set. Substituting the value of $\theta^\star$ yields the claimed projection.
Convexity of both the objective and constraint set ensures uniqueness. The last implication is obvious due to the formula for $w_i$.
\end{proof}

\subsection{Removing the global diameter dependence by Gaussian sample splitting}
\label{appendix:linear:type2}

In Remark~\ref{remark:linear:estimators:over:type-2:sets}, the passage from
local metric entropy to global metric entropy introduces the factor
\begin{align*}
    1+\log\frac{\diam(K)}{\eta^\star}.
\end{align*}
For a general convex body this quantity need not be logarithmic in the
ambient dimension, since $\diam(K)$ may be arbitrarily large compared with
the noise level. In the Gaussian sequence model, however, one can remove
this dependence on the global diameter by a sample-splitting localization
argument. We give the details below.

Throughout this appendix, let
\begin{align*}
    Y=\mu+\sigma Z,\qquad Z\sim N(0,I_n),\qquad \mu\in K,
\end{align*}
where $K\subset\RR^n$ is a compact, origin-symmetric convex body. Recall that
\begin{align*}
    \eta^\star
    =
    \sup\left\{
        \eta\geq0:
        \frac{\eta^2}{\sigma^2}
        \leq
        \log M_K^{\operatorname{loc}}(\eta)
    \right\}.
\end{align*}
We assume $\sigma>0$ and that $K$ is not a singleton; the remaining cases
are immediate.

We first record two elementary observations.

\begin{lemma}[A crude lower bound on the critical radius]
\label{lemma:eta-star-diameter-lower-bound}
Let $d=\diam(K)$. Then
\begin{align*}
    \eta^\star\gtrsim d\wedge\sigma.
\end{align*}
\end{lemma}

\begin{proof}
To be self contained we spell the details of this lemma. The same result can be found in \cite{prasadan2024some}.
Since $K$ is compact and origin-symmetric, there exists $v\in K$ with
$\|v\|_2=d/2$. By convexity,
\begin{align*}
    [-v,v]\subseteq K.
\end{align*}
Consequently, for every
\begin{align*}
    0<t\leq d/2,
\end{align*}
the set $B_2(0,t)\cap K$ contains the two points
$\pm t\,v/\|v\|_2$. Their Euclidean distance is $2t$, so, for the
sufficiently large absolute constant $c$ appearing in the definition of
local metric entropy,
\begin{align*}
    M_K^{\operatorname{loc}}(t)\geq2.
\end{align*}
Choose
\begin{align*}
    t_0
    =
    \min\left\{
        \frac d4,\,
        \frac{\sqrt{\log2}}{2}\sigma
    \right\}.
\end{align*}
Then
\begin{align*}
    \frac{t_0^2}{\sigma^2}\leq\log2
    \leq
    \log M_K^{\operatorname{loc}}(t_0),
\end{align*}
and therefore $t_0\leq\eta^\star$. This proves
$\eta^\star\gtrsim d\wedge\sigma$.
\end{proof}

\begin{lemma}[Type-$2$ is preserved by Euclidean localization]
\label{lemma:type2-euclidean-localization}
For every $R>0$, define
\begin{align*}
    L_R:=2K\cap RB_2.
\end{align*}
Then $L_R$ is an origin-symmetric convex body and
\begin{align*}
    T_2(L_R)
    \leq
    \sqrt{T_2(K)^2+1}
    \lesssim T_2(K).
\end{align*}
\end{lemma}

\begin{proof}
The Minkowski gauge of $L_R$ is
\begin{align*}
    \rho_{L_R}(x)
    =
    \max\left\{
        \frac12\rho_K(x),\,
        \frac{\|x\|_2}{R}
    \right\}.
\end{align*}
Hence, for arbitrary $x_1,\ldots,x_N\in\RR^n$ and independent Rademacher
variables $\epsilon_1,\ldots,\epsilon_N$,
\begin{align*}
    \EE_\epsilon
    \rho_{L_R}^2\left(\sum_{i=1}^N\epsilon_i x_i\right)
    &\leq
    \frac14
    \EE_\epsilon
    \rho_K^2\left(\sum_{i=1}^N\epsilon_i x_i\right)
    +
    \frac1{R^2}
    \EE_\epsilon
    \left\|\sum_{i=1}^N\epsilon_i x_i\right\|_2^2\\
    &\leq
    \frac{T_2(K)^2}{4}
    \sum_{i=1}^N\rho_K^2(x_i)
    +
    \frac1{R^2}\sum_{i=1}^N\|x_i\|_2^2\\
    &\leq
    \bigl(T_2(K)^2+1\bigr)
    \sum_{i=1}^N\rho_{L_R}^2(x_i).
\end{align*}
The claim follows.
\end{proof}

We next record the deterministic consequence of the local-to-global entropy
argument that we will use. The point is that the critical radius in the
following lemma is the critical radius of the original body $K$, not that
of the localized body.

\begin{lemma}[A small width for the localized error body]
\label{lemma:localized-body-small-width}
Fix $R>0$ and let
\begin{align*}
    L_R=2K\cap RB_2,
    \qquad
    D_R=\diam(L_R).
\end{align*}
There exists an orthogonal projection $P_m$ of rank $m\in\{0,\ldots,n\}$
such that
\begin{align}
    m
    &\lesssim
    \frac{\eta^{\star2}}{\sigma^2}
    \left(
        1+\log_+\frac{D_R}{\eta^\star}
    \right),
    \label{eq:localized-rank-bound}\\
    \sup_{h\in L_R}\|(I-P_m)h\|_2
    &\lesssim
    T_2(L_R)\eta^\star.
    \label{eq:localized-width-bound}
\end{align}
Here $\log_+(x)=\max\{0,\log x\}$. If $D_R\lesssim\eta^\star$, the
conclusion is understood with $m=0$.
\end{lemma}

\begin{proof}
Since $L_R\subseteq2K$, for every $\delta>0$,
\begin{align*}
    M_{L_R}^{\operatorname{loc}}(\delta)
    \leq
    M_{2K}^{\operatorname{loc}}(\delta)
    =
    M_K^{\operatorname{loc}}(\delta/2).
\end{align*}
By the definition of $\eta^\star$ and the monotonicity of local entropy,
for every $\delta\geq4\eta^\star$,
\begin{align}
    \log M_{L_R}^{\operatorname{loc}}(\delta)
    &\leq
    \log M_K^{\operatorname{loc}}(\delta/2)\nonumber\\
    &\leq
    \log M_K^{\operatorname{loc}}(2\eta^\star)
    \leq
    \frac{4\eta^{\star2}}{\sigma^2}.
    \label{eq:localized-local-entropy}
\end{align}
The last inequality follows from the definition of $\eta^\star$ (with an
irrelevant adjustment of absolute constants if the supremum in the
definition is attained).

Apply the same telescoping argument used in the proof of
Lemma~\ref{Kwidth:bound:lemma}. More explicitly, Lemma~3 of
\cite{yang1999information} gives, with the same sufficiently large constant
$c$ as in the definition of local entropy,
\begin{align*}
    \log M_{L_R}(a c^{j-1})
    -
    \log M_{L_R}(a c^j)
    \leq
    \log M_{L_R}^{\operatorname{loc}}(a c^j).
\end{align*}
Choose $a=C_0\eta^\star$, where $C_0$ is a sufficiently large absolute
constant. Summing this inequality until $a c^j>D_R$ and using
\eqref{eq:localized-local-entropy} yields
\begin{equation}
    \log M_{L_R}(C_1\eta^\star)
    \lesssim
    \frac{\eta^{\star2}}{\sigma^2}
    \left(
        1+\log_+\frac{D_R}{\eta^\star}
    \right).
    \label{eq:localized-global-entropy}
\end{equation}
A maximal packing is also a covering, so the same estimate, up to an
absolute change in $C_1$, bounds the logarithm of the corresponding
covering number.

Set
\begin{align*}
    H
    :=
    C_2
    \frac{\eta^{\star2}}{\sigma^2}
    \left(
        1+\log_+\frac{D_R}{\eta^\star}
    \right).
\end{align*}
If $D_R\leq C_4\eta^\star$, for a sufficiently large absolute
constant $C_4$, then, since $L_R$ is origin-symmetric,
\begin{align*}
    d_0(L_R)=\sup_{h\in L_R}\|h\|_2=D_R/2\lesssim\eta^\star,
\end{align*}
and we may take $m=0$.

Suppose therefore that $D_R>C_4\eta^\star$. Since
$D_R\leq 2d$ and Lemma~\ref{lemma:eta-star-diameter-lower-bound} gives
$\eta^\star\gtrsim d\wedge\sigma$, choosing $C_4$ sufficiently large
forces $d\gtrsim\sigma$, and hence
\begin{align*}
    \eta^\star\gtrsim\sigma.
\end{align*}
In particular $H\gtrsim1$. By increasing the absolute constant in $H$ if
necessary, \eqref{eq:localized-global-entropy} implies that for some
$k\leq e^H$,
\begin{align*}
    e_k(L_R)\lesssim\eta^\star.
\end{align*}
Take $m=\lceil C_3H\rceil\wedge n$. If $m=n$, then $d_m(L_R)=0$ and there
is nothing to prove. Otherwise, Carl's inequality \eqref{carls:ineq}
gives
\begin{align*}
    d_m(L_R)
    &\leq
    \left(\prod_{j=1}^m d_j(L_R)\right)^{1/m}\\
    &\lesssim
    T_2(L_R)\,
    k^{1/m}e_k(L_R)
    \lesssim
    T_2(L_R)\eta^\star,
\end{align*}
where $C_3$ is chosen large enough that $k^{1/m}\lesssim1$.
This proves both claims.
\end{proof}

We can now give the sample-splitting construction.

\begin{proposition}[Gaussian sample-splitting localization]
\label{prop:gaussian-sample-splitting-localization}
There exists a randomized estimator $\widehat\mu$ such that
\begin{equation}
    \sup_{\mu\in K}
    \EE_\mu\|\widehat\mu-\mu\|_2^2
    \lesssim
    \left(
        T_2(K)^2+\log(en)
    \right)\eta^{\star2}.
    \label{eq:sample-split-risk}
\end{equation}
In particular, the potentially non-dimensional factor
$\log(\diam(K)/\eta^\star)$ from
Remark~\ref{remark:linear:estimators:over:type-2:sets} can be replaced, in
the Gaussian sequence model, by an $O(\log n)$ factor. The estimator
constructed below is not a single fixed linear estimator of $Y$; rather,
it is a sample-split estimator consisting of a pilot estimate followed by
a linear correction.
\end{proposition}

\begin{proof}
Let $G\sim N(0,I_n)$ be independent of $Y$, and define
\begin{equation}
    Y^{(1)}=Y+\sigma G,
    \qquad
    Y^{(2)}=Y-\sigma G.
    \label{eq:gaussian-random-split}
\end{equation}
Since
\begin{align*}
    Y=\mu+\sigma Z,
\end{align*}
we have
\begin{align*}
    Y^{(1)}-\mu=\sigma(Z+G),
    \qquad
    Y^{(2)}-\mu=\sigma(Z-G).
\end{align*}
The two vectors $Z+G$ and $Z-G$ are jointly Gaussian and have zero
cross-covariance. Therefore
\begin{equation}
    Y^{(1)},Y^{(2)}
    \stackrel{\operatorname{ind}}{\sim}
    N(\mu,2\sigma^2I_n).
    \label{eq:independent-gaussian-split}
\end{equation}

Use the first copy to form the Euclidean projection pilot
\begin{align*}
    \widetilde\mu:=\Pi_K(Y^{(1)}).
\end{align*}
Because Euclidean projection onto a closed convex set is non-expansive and
$\Pi_K(\mu)=\mu$,
\begin{equation}
    \|\widetilde\mu-\mu\|_2
    \leq
    \|Y^{(1)}-\mu\|_2.
    \label{eq:pilot-contraction}
\end{equation}
Moreover, since both $\widetilde\mu$ and $\mu$ belong to $K$,
\begin{equation}
    h:=\mu-\widetilde\mu\in K-K=2K.
    \label{eq:error-in-2K}
\end{equation}

Fix a sufficiently large absolute constant $A$ and set
\begin{align*}
    R=A\sqrt n\,\sigma,
    \qquad
    L_R=2K\cap RB_2.
\end{align*}
Let $P_m$ be the projection furnished by
Lemma~\ref{lemma:localized-body-small-width}, and define
\begin{equation}
    \widehat\mu
    :=
    \widetilde\mu
    +
    P_m\bigl(Y^{(2)}-\widetilde\mu\bigr).
    \label{eq:sample-split-estimator}
\end{equation}
Conditional on $Y^{(1)}$, the vector $h$ is fixed, while
$Y^{(2)}-\mu\sim N(0,2\sigma^2I_n)$ is independent of $Y^{(1)}$.
Therefore
\begin{align}
    \widehat\mu-\mu
    &=
    -(I-P_m)h
    +
    P_m(Y^{(2)}-\mu),
    \label{eq:sample-split-error}
\end{align}
and hence
\begin{equation}
    \EE\!\left[
        \|\widehat\mu-\mu\|_2^2
        \,\middle|\,
        Y^{(1)}
    \right]
    =
    \|(I-P_m)h\|_2^2
    +
    2\sigma^2m.
    \label{eq:conditional-sample-split-risk}
\end{equation}

Let
\begin{align*}
    \mathcal E:=\{\|h\|_2\leq R\}.
\end{align*}
On $\mathcal E$, \eqref{eq:error-in-2K} implies $h\in L_R$, and therefore
\begin{align*}
    \|(I-P_m)h\|_2^2
    \lesssim
    T_2(L_R)^2\eta^{\star2}.
\end{align*}
On $\mathcal E^c$, since $P_m$ is an orthogonal projection,
\begin{align*}
    \|(I-P_m)h\|_2^2\leq\|h\|_2^2.
\end{align*}
Taking expectations in \eqref{eq:conditional-sample-split-risk} gives
\begin{align}
    \EE_\mu\|\widehat\mu-\mu\|_2^2
    &\leq
    2\sigma^2m
    +
    C T_2(L_R)^2\eta^{\star2}
    +
    \EE_\mu\left[
        \|h\|_2^2\mathbf 1_{\{\|h\|_2>R\}}
    \right].
    \label{eq:risk-before-tail}
\end{align}

We first control the deterministic terms. Since
\begin{align*}
    D_R=\diam(L_R)
    \leq
    2\min\{d,R\},
\end{align*}
Lemma~\ref{lemma:eta-star-diameter-lower-bound} gives
\begin{align}
    \frac{D_R}{\eta^\star}
    &\lesssim
    \frac{\min\{d,\sqrt n\,\sigma\}}{d\wedge\sigma}
    \lesssim
    \sqrt n.
    \label{eq:localized-log-is-logn}
\end{align}
Consequently,
\begin{align*}
    1+\log_+\frac{D_R}{\eta^\star}
    \lesssim
    \log(en).
\end{align*}
Combining this with
Lemma~\ref{lemma:localized-body-small-width} and
Lemma~\ref{lemma:type2-euclidean-localization},
\begin{equation}
    2\sigma^2m
    +
    C T_2(L_R)^2\eta^{\star2}
    \lesssim
    \left(
        T_2(K)^2+\log(en)
    \right)\eta^{\star2}.
    \label{eq:good-part-risk}
\end{equation}

It remains only to control the tail term in \eqref{eq:risk-before-tail}.
By \eqref{eq:pilot-contraction},
\begin{align*}
    \|h\|_2
    \leq
    \|Y^{(1)}-\mu\|_2.
\end{align*}
If $d\leq R$, then $\|h\|_2\leq d\leq R$ deterministically, because
$h\in2K$ and
\begin{align*}
    \sup_{x\in2K}\|x\|_2=d.
\end{align*}
Thus the tail term vanishes.

Suppose now that $d>R$. By
Lemma~\ref{lemma:eta-star-diameter-lower-bound},
\begin{align}
    \eta^\star\gtrsim\sigma.
    \label{eq:eta-star-at-least-sigma-large-d}
\end{align}
Furthermore,
\begin{align*}
    \frac{Y^{(1)}-\mu}{\sqrt2\sigma}\sim N(0,I_n).
\end{align*}
A standard Gaussian norm tail bound implies, for $A$ sufficiently large,
\begin{equation}
    \EE_\mu\left[
        \|Y^{(1)}-\mu\|_2^2
        \mathbf 1_{\{\|Y^{(1)}-\mu\|_2>R\}}
    \right]
    \lesssim
    n\sigma^2e^{-c n}.
    \label{eq:gaussian-tail-second-moment}
\end{equation}
For completeness, \eqref{eq:gaussian-tail-second-moment} follows from
\begin{align*}
    \PP(\|Z\|_2\geq\sqrt n+t)\leq e^{-t^2/2},
\end{align*}
together with Cauchy--Schwarz and
$\EE\|Z\|_2^4=n(n+2)$.
Using \eqref{eq:pilot-contraction},
\eqref{eq:eta-star-at-least-sigma-large-d}, and the elementary bound
$n e^{-cn}\lesssim1$, we obtain
\begin{align}
    \EE_\mu\left[
        \|h\|_2^2\mathbf 1_{\{\|h\|_2>R\}}
    \right]
    &\leq
    \EE_\mu\left[
        \|Y^{(1)}-\mu\|_2^2
        \mathbf 1_{\{\|Y^{(1)}-\mu\|_2>R\}}
    \right]\nonumber\\
    &\lesssim
    n\sigma^2e^{-cn}
    \lesssim
    \sigma^2
    \lesssim
    \eta^{\star2}.
    \label{eq:tail-absorbed}
\end{align}
Combining \eqref{eq:risk-before-tail},
\eqref{eq:good-part-risk}, and \eqref{eq:tail-absorbed} proves
\eqref{eq:sample-split-risk}.
\end{proof}

\begin{corollary}[The auxiliary randomization can be removed]
\label{cor:sample-split-derandomization}
There exists a deterministic estimator $\overline\mu=\overline\mu(Y)$
satisfying the same risk bound \eqref{eq:sample-split-risk}.
\end{corollary}

\begin{proof}
Write the estimator in Proposition~
\ref{prop:gaussian-sample-splitting-localization} as
$\widehat\mu(Y,G)$, where $G\sim N(0,I_n)$ is the auxiliary Gaussian
randomization, independent of $Y$. Define
\begin{align*}
    \overline\mu(Y)
    :=
    \EE_G\!\left[\widehat\mu(Y,G)\mid Y\right].
\end{align*}
This is a genuine estimator depending only on the observed data $Y$; in
particular, the averaging distribution of $G$ does not involve the
unknown parameter $\mu$. Conditional Jensen's inequality gives, for every
$\mu\in K$,
\begin{align*}
    \|\overline\mu(Y)-\mu\|_2^2
    &=
    \left\|
        \EE_G\!\left[
            \widehat\mu(Y,G)-\mu
            \,\middle|\,
            Y
        \right]
    \right\|_2^2\\
    &\leq
    \EE_G\!\left[
        \|\widehat\mu(Y,G)-\mu\|_2^2
        \,\middle|\,
        Y
    \right].
\end{align*}
Taking expectation with respect to $Y$ proves the claim.

For the particular construction above one may write this deterministic
estimator more explicitly as
\begin{align*}
    \overline\mu(Y)
    =
    P_mY
    +
    (I-P_m)
    \EE_G\!\left[\Pi_K(Y+\sigma G)\right],
\end{align*}
because $\EE_G G=0$.
\end{proof}

\begin{remark}
The construction above should be viewed as a statistical localization
argument rather than as a proof that a single deterministic linear
estimator is near-minimax. Indeed, the estimator
\eqref{eq:sample-split-estimator} can be written as
\begin{align*}
    \widehat\mu
    =
    P_mY^{(2)}+(I-P_m)\widetilde\mu,
\end{align*}
and its offset depends on the independent pilot sample $Y^{(1)}$.
Thus the argument removes the dependence on the global diameter at the
price of leaving the class of fixed linear estimators. By
Corollary~\ref{cor:sample-split-derandomization}, the auxiliary
randomization itself is inessential, but the resulting deterministic
estimator remains nonlinear. The argument also does not, by itself,
address the computational problem of finding the appropriate
Kolmogorov-width projection $P_m$.
\end{remark}

\section{Proofs of Section \ref{robust:section}}

\begin{lemma}[Robust refinement lemma]
\label{sub:gaussian:vector:estimator:robust}
Suppose that $\hat\mu_j\in K$ is an estimator of $\mu$ and let
\begin{align*}
    H_j
    :=
    \left\{
        \|\hat\mu_j-\mu\|_2\leq \tilde d_j
    \right\}.
\end{align*}
Set
\begin{align*}
    K_{(j)}
    :=
    K\cap B_2(0,\tilde d_j/2),
    \qquad
    k_j
    :=
    \left\lceil
        \frac{N\tilde d_j^2}{C^2\sigma^2}
    \right\rceil\wedge N.
\end{align*}
Assume that the two conditions of
Corollary~\ref{Kwidth:bound:lemma:robust} hold for $K_{(j)}$, that
\begin{align*}
    \tilde d_j
    \geq
    C\frac{\sigma}{\sqrt N},
    \qquad
    k_j\gtrsim \varepsilon N,
\end{align*}
and obtain $X_j^{\star\star}$ from
Algorithm~\ref{algorithm:SDP:Kolmogorov:width}, with
Theorem~\ref{main:SDP:theorem} applied using $\sigma/\sqrt N$ in place of
$\sigma$.  Put
\begin{align*}
    A_j^{\star\star}
    :=
    (I-X_j^{\star\star})^{1/2},
\end{align*}
and define
\begin{equation}
    \tilde\mu_{j+1}
    :=
    \frac{
        f((A_j^{\star\star}Y_i)_{i\in[N]})
        -
        A_j^{\star\star}\hat\mu_j
    }{2},
    \label{eq:robust-refinement-estimator}
\end{equation}
where $f$ is the procedure of \cite{depersin2022} with
$k=u=k_j$.

Then there exist absolute constants $L,C_1>0$ such that
\begin{align}
    \PP\left(
        \left\|
            \tilde\mu_{j+1}
            -
            \frac{\mu-\hat\mu_j}{2}
        \right\|_2
        >
        \frac{\tilde d_j}{L},
        \ H_j
    \right)
    \leq
    \exp(-k_j/C_1)
    +
    \frac{1}{
        1+N^{3/2}\tilde d_j^3/\sigma^3
    }.
    \label{eq:robust-refinement-probability}
\end{align}
Moreover, by choosing the absolute constants $c$ in the definition of local
entropy and $C$ in the definition of $k_j$ sufficiently large, $L$ may be
taken to be any prescribed sufficiently large absolute constant.
\end{lemma}

\begin{proof}
Let
\begin{align*}
    \nu_j
    :=
    \frac{\mu-\hat\mu_j}{2}.
\end{align*}
On the event $H_j$, since $\mu,\hat\mu_j\in K$ and $K$ is convex and
origin-symmetric,
\begin{align*}
    \nu_j\in K,
    \qquad
    \|\nu_j\|_2\leq \tilde d_j/2,
\end{align*}
and hence
\begin{equation}
    \nu_j\in K_{(j)}.
    \label{eq:robust-nuj-in-Kj}
\end{equation}

We first control the deterministic approximation error produced by
$A_j^{\star\star}$.  Let
\begin{align*}
    X_j^{\star\star}
    =
    V\Lambda V^\top
    =
    \sum_{\ell=1}^n
        \lambda_\ell v_\ell v_\ell^\top,
\end{align*}
so that
\begin{align*}
    A_j^{\star\star}
    =
    (I-X_j^{\star\star})^{1/2}
    =
    \sum_{\ell=1}^n
        \sqrt{1-\lambda_\ell}\,
        v_\ell v_\ell^\top.
\end{align*}
On the event that the QFM calls in
Algorithm~\ref{algorithm:SDP:Kolmogorov:width} all succeed,
Theorem~\ref{main:SDP:theorem} and
Corollary~\ref{Kwidth:bound:lemma:robust} give
\begin{equation}
    \sup_{p\in K_{(j)}}
        p^\top X_j^{\star\star}p
    \lesssim
    \frac{\tilde d_j^2}{c^2}.
    \label{eq:robust-SDP-bias-bound}
\end{equation}
Since $0\leq\lambda_\ell\leq1$, we have
\begin{align}
    \|\nu_j-A_j^{\star\star}\nu_j\|_2^2
    &=
    \sum_{\ell=1}^n
        \left(
            1-\sqrt{1-\lambda_\ell}
        \right)^2
        (v_\ell^\top\nu_j)^2
        \nonumber\\
    &=
    \sum_{\ell=1}^n
        \frac{\lambda_\ell^2}{
            \left(
                1+\sqrt{1-\lambda_\ell}
            \right)^2
        }
        (v_\ell^\top\nu_j)^2
        \nonumber\\
    &\leq
    \sum_{\ell=1}^n
        \lambda_\ell
        (v_\ell^\top\nu_j)^2
    =
    \nu_j^\top X_j^{\star\star}\nu_j
    \lesssim
    \frac{\tilde d_j^2}{c^2},
    \label{eq:robust-linear-bias}
\end{align}
where in the last inequality we used
\eqref{eq:robust-nuj-in-Kj} and
\eqref{eq:robust-SDP-bias-bound}.  Therefore
\begin{equation}
    \|\nu_j-A_j^{\star\star}\nu_j\|_2
    \lesssim
    \frac{\tilde d_j}{c}.
    \label{eq:robust-linear-bias-root}
\end{equation}

We next control the stochastic error.  Conditional on
$A_j^{\star\star}$, the uncontaminated transformed observations have mean
$A_j^{\star\star}\mu$ and covariance
\begin{align*}
    A_j^{\star\star}\Sigma A_j^{\star\star}.
\end{align*}
Moreover, applying the deterministic linear map
$A_j^{\star\star}$ does not increase the fraction of contaminated
observations.  By
Lemma~\ref{sub:gaussian:vector:estimator:1:robust}, since
$k_j\gtrsim\varepsilon N$ and
$\tilde d_j\geq C\sigma/\sqrt N$, with probability at least
$1-\exp(-k_j/C_1)$,
\begin{equation}
    \left\|
        f((A_j^{\star\star}Y_i)_{i\in[N]})
        -
        A_j^{\star\star}\mu
    \right\|_2
    \leq
    \frac{\tilde d_j}{T},
    \label{eq:robust-projected-mean-error}
\end{equation}
after changing the absolute constant $C_1$ if necessary.  Here the additional
$10^{-u}$ failure probability in
Theorem~\ref{thm:appendix:depersin:2} is absorbed into
$\exp(-k_j/C_1)$ because $u=k_j$.

Using the definition \eqref{eq:robust-refinement-estimator}, we have
\begin{align}
    \tilde\mu_{j+1}-\nu_j
    &=
    \frac12
    \left[
        f((A_j^{\star\star}Y_i)_{i\in[N]})
        -
        A_j^{\star\star}\mu
    \right]
    +
    \left(
        A_j^{\star\star}\nu_j-\nu_j
    \right).
    \label{eq:robust-refinement-decomposition}
\end{align}
Consequently, on $H_j$, on the event that the SDP construction succeeds,
and on the event \eqref{eq:robust-projected-mean-error},
\begin{align}
    \left\|
        \tilde\mu_{j+1}-\nu_j
    \right\|_2
    &\leq
    \frac12
    \left\|
        f((A_j^{\star\star}Y_i)_{i\in[N]})
        -
        A_j^{\star\star}\mu
    \right\|_2
    +
    \|A_j^{\star\star}\nu_j-\nu_j\|_2
    \nonumber\\
    &\lesssim
    \frac{\tilde d_j}{2T}
    +
    \frac{\tilde d_j}{c}.
    \label{eq:robust-refinement-final-bound}
\end{align}
Choosing the absolute constants in the construction sufficiently large gives
\begin{align*}
    \left\|
        \tilde\mu_{j+1}
        -
        \frac{\mu-\hat\mu_j}{2}
    \right\|_2
    \leq
    \frac{\tilde d_j}{L}.
\end{align*}

It remains only to collect the failure probabilities.  Since
Theorem~\ref{main:SDP:theorem} is applied here with the effective noise scale
$\sigma/\sqrt N$, its failure probability is
\begin{align*}
    \frac{1}{
        1+\tilde d_j^3/(\sigma/\sqrt N)^3
    }
    &=
    \frac{1}{
        1+N^{3/2}\tilde d_j^3/\sigma^3
    }.
\end{align*}
The projected robust-mean step fails with probability at most
$\exp(-k_j/C_1)$.  A union bound therefore yields
\eqref{eq:robust-refinement-probability}.
\end{proof}

\begin{proof}[Proof of Theorem \ref{robust:main:thm}]
Let
\begin{align*}
    \rho:=\frac{2(\sqrt{3}+1)}{\widetilde L}<1,
\end{align*}
so that
\begin{align*}
    \widetilde d_j=\widetilde d_1\rho^{j-1},
    \qquad
    \widetilde d_{j+1}=\rho\widetilde d_j.
\end{align*}
Write
\begin{align*}
    \delta^\star:=\delta^\star(K),
    \qquad
    \eta^\star:=\eta^\star(K)
    =\delta^\star\vee\sqrt{\varepsilon}\,\sigma,
\end{align*}
and set
\begin{align*}
    \overline\kappa:=\max_{j\in[M]}\kappa(K_{(j)}).
\end{align*}
As in the proof of Theorem~\ref{main:result:GSM}, we repeatedly use
\begin{align*}
    T_2(K_{(j)})\leq \sqrt{2}\,T_2(K),
    \qquad
    \delta^\star(K_{(j)})\leq \delta^\star(K).
\end{align*}
The second inequality follows from
$K_{(j)}\subseteq K$ and monotonicity of the local entropy under inclusion.
Define
\begin{align*}
    \Lambda_K
    :=
    1+\left\lceil
        \log_c\!\left(
            \overline\kappa^{1/2}\sqrt{2}\,T_2(K)
        \right)
    \right\rceil
\end{align*}
and
\begin{equation}\label{eq:robust-GammaK}
    \Gamma_K
    :=
    1
    \vee
    \left(2\sqrt{2}\,T_2(K)\overline\kappa^{1/2}\right)
    \vee
    \left(2C\sqrt{\Lambda_K}\right).
\end{equation}
All constants below may be enlarged by absolute factors, so that the
additional constants appearing in the robust mean estimator and in the weak
projections are absorbed into $\Gamma_K$.

We will also use the elementary lower bound
\begin{equation}\label{eq:robust-parametric-floor}
    \eta^\star(K)
    \gtrsim
    \frac{\sigma}{\sqrt N}
    \wedge d.
\end{equation}
Indeed, since $K$ is convex and has diameter $d$, it contains a line
segment of length arbitrarily close to $d$. Hence, for every
$0<t\lesssim d$, a ball of radius comparable to $t$, centered at a
suitable point of this segment, contains two points of $K$ separated by a
constant multiple of $t$. Thus
$\log M_K^{\operatorname{loc}}(t)\gtrsim1$, and taking
$t\asymp (\sigma/\sqrt N)\wedge d$ in the defining entropy inequality
gives \eqref{eq:robust-parametric-floor}. In particular, if
$d\lesssim\sigma/\sqrt N$, then
$\eta^\star(K)\gtrsim d$, and the theorem is immediate from properness
of the estimator. We therefore assume from now on that
\begin{equation}\label{eq:robust-nontrivial-diameter}
    d\gtrsim\frac{\sigma}{\sqrt N},
    \qquad\text{so that}\qquad
    \eta^\star(K)\gtrsim\frac{\sigma}{\sqrt N}.
\end{equation}

We first control the initialization. Let
\begin{align*}
    G_0
    :=
    \left\{
        \|\widehat\mu_1-\mu\|_2\leq \widetilde d_1
    \right\}.
\end{align*}
By \eqref{guarantee:depersin:1st}, Lemma~\ref{weak:projection:lemma}, and the
choice of the absolute constant $C'$ in
$\widetilde d_1=2R\wedge C'\sigma$, we may choose $C'$ sufficiently large
so that
\begin{equation}\label{eq:robust-pilot}
    \PP(G_0^c)\leq \exp(-N/C_0).
\end{equation}
Indeed, if $\widetilde d_1=C'\sigma$, this follows directly from the pilot
bound after weak projection, while if $\widetilde d_1=2R$, then both
$\widehat\mu_1$ and $\mu$ belong to $K$, and hence
$\|\widehat\mu_1-\mu\|_2\leq \diam(K)\leq 2R=\widetilde d_1$
deterministically.

For $j\geq1$, define the trapping event
\begin{align*}
    H_j
    :=
    \left\{
        \|\widehat\mu_j-\mu\|_2\leq \widetilde d_j
    \right\}.
\end{align*}
Thus $G_0\subseteq H_1$. Suppose that $H_j$ holds and that the two
conditions of Corollary~\ref{Kwidth:bound:lemma:robust} hold for $K_{(j)}$.
Put
\begin{align*}
    k_j
    :=
    \left\lceil
        \frac{N\widetilde d_j^2}{C^2\sigma^2}
    \right\rceil\wedge N.
\end{align*}
The robust analogue of Lemma~\ref{sub:gaussian:vector:estimator}, obtained
from Lemma~\ref{sub:gaussian:vector:estimator:1:robust}, yields the good
refinement event
\begin{align*}
    E_j
    :=
    \left\{
        \left\|
            \widetilde\mu_{j+1}
            -\frac{\mu-\widehat\mu_j}{2}
        \right\|_2
        \leq
        \frac{\widetilde d_j}{L}
    \right\}.
\end{align*}
Whenever
\begin{equation}\label{eq:robust-operating-scale}
    \widetilde d_j
    \gtrsim
    \frac{\sigma}{\sqrt N}
    \vee
    \sqrt{\varepsilon}\,\sigma,
\end{equation}
we have $k_j\gtrsim \varepsilon N$, and therefore
\begin{equation}\label{eq:robust-Ej-prob}
\begin{split}
    \PP(E_j^c\cap H_j)
    &\leq
    \exp(-k_j/C_1)
    +
    \frac{1}{1+N^{3/2}\widetilde d_j^3/\sigma^3} \\
    &\lesssim
    \exp\!\left(
        -c_1\frac{N\widetilde d_j^2}{\sigma^2}
    \right)
    +
    \frac{1}{1+N^{3/2}\widetilde d_j^3/\sigma^3},
\end{split}
\end{equation}
for absolute constants $C_1,c_1>0$. In the second line we used
$\widetilde d_j\lesssim\sigma$, which follows from
$\widetilde d_j\leq\widetilde d_1\lesssim\sigma$; thus the truncation
$k_j\leq N$ changes only the absolute constant in the exponent.

We next check that the same deterministic contraction as in the Gaussian
sequence proof holds on $H_j\cap E_j$. On $H_j$,
\begin{align*}
    \nu_j:=\frac{\mu-\widehat\mu_j}{2}
    \in K\cap B_2(0,\widetilde d_j/2)=K_{(j)},
\end{align*}
because $K$ is convex and origin symmetric. Applying
Lemma~\ref{weak:projection:lemma} to the projection of
$\widetilde\mu_{j+1}$ onto $K_{(j)}$, with
$\epsilon=\widetilde d_j/L$, gives on $H_j\cap E_j$
\begin{align*}
\begin{split}
    \|\overline\mu_{j+1}-\nu_j\|_2
    &\leq
    \|\widetilde\mu_{j+1}-\nu_j\|_2
    +
    \sqrt{
        \epsilon^2
        +2\epsilon\|\widetilde\mu_{j+1}-\nu_j\|_2
    } \\
    &\leq
    \frac{(1+\sqrt3)\widetilde d_j}{L}
    =
    \frac{\widetilde d_j}{\widetilde L}.
\end{split}
\end{align*}
Consequently,
\begin{align*}
    \|2\overline\mu_{j+1}+\widehat\mu_j-\mu\|_2
    \leq
    \frac{2\widetilde d_j}{\widetilde L}.
\end{align*}
A second application of Lemma~\ref{weak:projection:lemma}, now to the
projection onto $K$ with
$\epsilon=2\widetilde d_j/\widetilde L$, yields
\begin{align*}
\begin{split}
    \|\widehat\mu_{j+1}-\mu\|_2
    &\leq
    \frac{2\widetilde d_j}{\widetilde L}
    +
    \sqrt{
        \frac{4\widetilde d_j^2}{\widetilde L^2}
        +
        \frac{8\widetilde d_j^2}{\widetilde L^2}
    } \\
    &=
    \frac{2(\sqrt3+1)\widetilde d_j}{\widetilde L}
    =
    \widetilde d_{j+1}.
\end{split}
\end{align*}
Hence, as long as the deterministic assumptions and
\eqref{eq:robust-operating-scale} continue to hold,
\begin{equation}\label{eq:robust-good-implies-next}
    H_j\cap E_j\subseteq H_{j+1}.
\end{equation}
In particular, no independence across iterations is required.

Exactly as in the proof of Theorem~\ref{main:result:GSM}, after the first
scale at which the contraction argument ceases to apply, subsequent iterates
cannot move farther than a constant multiple of that scale. For completeness,
let
\begin{align*}
    a_\ell:=\|\widehat\mu_\ell-\mu\|_2.
\end{align*}
If $a_k\leq\widetilde d_k$, then, regardless of what happens at later
iterations, $\overline\mu_\ell\in K_{(\ell-1)}$, and hence
$2\|\overline\mu_\ell\|_2\leq\widetilde d_{\ell-1}$. Therefore, with
$\epsilon_{\ell-1}=2\widetilde d_{\ell-1}/\widetilde L$,
\begin{align*}
    a_\ell
    \leq
    a_{\ell-1}
    +C_2\widetilde d_{\ell-1}
    +C_3\sqrt{\widetilde d_{\ell-1}a_{\ell-1}}.
\end{align*}
Since $\widetilde d_\ell$ decreases geometrically, the same geometric-series
argument as in Theorem~\ref{main:result:GSM} gives
\begin{equation}\label{eq:robust-post-failure-stability}
    \max_{\ell\geq k}
    \|\widehat\mu_\ell-\mu\|_2
    \leq
    C_4\widetilde d_k
\end{equation}
for an absolute constant $C_4$.

We now identify all scales at which the induction may stop for deterministic
reasons. First, if the first condition of
Corollary~\ref{Kwidth:bound:lemma:robust} fails, then
\begin{align*}
    \widetilde d_j
    <
    2\kappa(K_{(j)})^{1/2}T_2(K_{(j)})
    \delta^\star(K_{(j)})
    \leq
    2\sqrt2\,\overline\kappa^{1/2}T_2(K)\eta^\star(K)
    \leq
    \Gamma_K\eta^\star(K).
\end{align*}
Second, suppose the other condition of
Corollary~\ref{Kwidth:bound:lemma:robust} fails. If
\begin{align*}
    L_j
    :=
    \left\lceil
        \log_c\!\left(
            \kappa(K_{(j)})^{1/2}T_2(K_{(j)})
        \right)
    \right\rceil,
\end{align*}
then the failure implies
\begin{align*}
    \frac{N\widetilde d_j^2}{C^2\sigma^2}
    <
    \left\lceil
        \frac{4L_jN\delta^{\star2}(K_{(j)})}{\sigma^2}
    \right\rceil,
\end{align*}
and hence
\begin{align*}
    \widetilde d_j^2
    \lesssim
    C^2\left(
        L_j\delta^{\star2}(K_{(j)})
        +\frac{\sigma^2}{N}
    \right).
\end{align*}
Now $L_j\leq\Lambda_K$,
$\delta^\star(K_{(j)})\leq\eta^\star(K)$, and
\eqref{eq:robust-nontrivial-diameter} gives
$\sigma/\sqrt N\lesssim\eta^\star(K)$. Thus
\begin{align*}
    \widetilde d_j
    \lesssim
    C\sqrt{\Lambda_K}\,\eta^\star(K)
    \lesssim
    \Gamma_K\eta^\star(K).
\end{align*}
Third, if \eqref{eq:robust-operating-scale} fails, then directly
\begin{align*}
    \widetilde d_j
    \lesssim
    \frac{\sigma}{\sqrt N}
    \vee
    \sqrt{\varepsilon}\,\sigma
    \lesssim
    \eta^\star(K).
\end{align*}
In particular, the condition $k_j\gtrsim\varepsilon N$, needed for the
procedure of \cite{depersin2022}, can fail only after the algorithm has
already reached the minimax scale.

It remains to verify that the deterministic stopping rule in
Algorithm~\ref{algorithm:main:robust} also occurs at the minimax scale. Let
$r_0$ denote the original inner radius before Algorithm~\ref{algorithm:main:robust}
overwrites $r$, and put
\begin{align*}
    r_{\mathrm{alg}}
    :=
    r_0\wedge
    \left(
        \sqrt{n/N}\,\sigma
        \vee
        \sqrt{\varepsilon}\,\sigma
    \right).
\end{align*}
Since $r_0B_2\subseteq K$, the minimax risk over $K$ is at least the
minimax risk over $r_0B_2$. By the Euclidean-ball lower bound quoted above
from \cite{peng2026robust},
\begin{align*}
    \eta^{\star2}(K)\wedge d^2
    \gtrsim
    \left(
        \frac{n\sigma^2}{N}
        \vee
        \varepsilon\sigma^2
    \right)
    \wedge r_0^2.
\end{align*}
Consequently
\begin{equation}\label{eq:robust-r-minimax}
    r_{\mathrm{alg}}
    \lesssim
    \eta^\star(K)
\end{equation}
(up to the immaterial case in which the diameter itself is already the
smaller scale, in which case the theorem is immediate because the estimator
is proper). Therefore the rule
\begin{align*}
    \widetilde d_{j+1}
    \leq
    2r_{\mathrm{alg}}
    \vee
    \frac{C\sigma}{\sqrt N}
\end{align*}
can only stop the algorithm once
$\widetilde d_{j+1}\lesssim\Gamma_K\eta^\star(K)$. The same conclusion
holds if all $M$ iterations are completed, by the definition of $M$.

We now obtain the tail bound. Let
\begin{align*}
    t_0:=C_5\Gamma_K\eta^\star(K),
\end{align*}
where $C_5$ is a sufficiently large absolute constant, in particular large
enough relative to the constant in
\eqref{eq:robust-post-failure-stability}. On the pilot event $G_0$, the
preceding discussion shows that a deterministic termination of the induction
cannot produce final error at least $t$ for $t\geq t_0$. Hence, if
$t\geq t_0$, $G_0$ occurs, and
$\|\widehat\mu-\mu\|_2\geq t$, then for some iteration $j$ with
$\widetilde d_j\gtrsim t$ the event $E_j^c\cap H_j$ must occur. Therefore,
using \eqref{eq:robust-pilot}, \eqref{eq:robust-Ej-prob}, and a union bound,
for $t_0\leq t\leq d$,
\begin{equation}\label{eq:robust-tail-sum}
\begin{split}
    \PP(\|\widehat\mu-\mu\|_2\geq t)
    &\leq
    \exp(-N/C_0) \\
    &\quad+
    \sum_{j:\,\widetilde d_j\gtrsim t}
    \left[
        \exp\!\left(
            -c_1\frac{N\widetilde d_j^2}{\sigma^2}
        \right)
        +
        \frac{1}{1+N^{3/2}\widetilde d_j^3/\sigma^3}
    \right].
\end{split}
\end{equation}
Since the $\widetilde d_j$'s form a geometric sequence, exactly as in the
proof of Theorem~\ref{main:result:GSM}, the two sums are dominated by their
smallest-scale terms. Thus
\begin{equation}\label{eq:robust-tail}
    \PP(\|\widehat\mu-\mu\|_2\geq t)
    \lesssim
    \exp(-N/C_0)
    +
    \exp\!\left(-c_2\frac{Nt^2}{\sigma^2}\right)
    +
    \frac{\sigma^3}{N^{3/2}t^3},
    \qquad
    t_0\leq t\leq d.
\end{equation}
For $t>d$, the probability on the left is zero because
$\widehat\mu,\mu\in K$.

Finally, integrate the tail. Using the trivial probability bound below
$t_0$, \eqref{eq:robust-tail} on $[t_0,d]$, and properness of the
estimator, we obtain
\begin{align*}
\begin{split}
    \EE\|\widehat\mu-\mu\|_2^2
    &\lesssim
    t_0^2
    +
    d^2\exp(-N/C_0) \\
    &\quad+
    \int_{t_0}^{\infty}
        t\exp\!\left(-c_2\frac{Nt^2}{\sigma^2}\right)\,dt
    +
    \frac{\sigma^3}{N^{3/2}}
    \int_{t_0}^{\infty}\frac{dt}{t^2} \\
    &\lesssim
    \Gamma_K^2\eta^{\star2}(K)
    +
    d^2\exp(-N/C_0)
    +
    \frac{\sigma^2}{N}
        \exp\!\left(-c_3\frac{Nt_0^2}{\sigma^2}\right)
    +
    \frac{\sigma^3}{N^{3/2}t_0}.
\end{split}
\end{align*}
By the assumption
\begin{align*}
    \frac{d^2}{\sigma^2}
    \lesssim
    \frac{\exp(N/C_0)}{N},
\end{align*}
the pilot-failure contribution satisfies
\begin{align*}
    d^2\exp(-N/C_0)
    \lesssim
    \frac{\sigma^2}{N}
    \lesssim
    \eta^{\star2}(K).
\end{align*}
Moreover, \eqref{eq:robust-parametric-floor} and
$t_0\gtrsim\eta^\star(K)$ imply
\begin{align*}
    \frac{\sigma^3}{N^{3/2}t_0}
    \lesssim
    \eta^{\star2}(K),
\end{align*}
and the exponential integral is bounded by the same quantity. Consequently,
\begin{align*}
    \EE\|\widehat\mu-\mu\|_2^2
    \lesssim
    \Gamma_K^2\eta^{\star2}(K).
\end{align*}
Since the estimator is proper, we also have the deterministic bound
$\|\widehat\mu-\mu\|_2\leq d$, and hence
\begin{align*}
    \EE\|\widehat\mu-\mu\|_2^2
    \lesssim
    \Gamma_K^2\eta^{\star2}(K)\wedge d^2.
\end{align*}
Thus the theorem holds with
$\Upsilon(K)\asymp\Gamma_K^2$, which is polynomial up to logarithmic
factors in $T_2(K)$ and $\overline\kappa$.

Finally, Algorithm~\ref{algorithm:main:robust} uses
$M=O(\log(\widetilde d_1/r_{\mathrm{alg}}))$ iterations. The assumptions
that $\log(\widetilde d_1/r_{\mathrm{alg}})$ and
$\log(1+1/\sigma)$ are polynomial in $n,N$, together with the same QFM
complexity argument as in Theorem~\ref{main:result:GSM}, imply that the entire
procedure runs in polynomial time.
\end{proof}

\section{Proofs of Section \ref{regression:section}}

\begin{proof}[Proof of Lemma \ref{sub:gaussian:vector:estimator:regression}]
The QFM oracle will make a mistake in Algorithm \ref{algorithm:SDP:Kolmogorov:width} with probability at most $1/(1 + N^{3/2}\tilde d_j^3)$. We now place ourselves on the event that the matrix $A_j^{\star\star}$ satisfies the guarantee of Theorem \ref{main:SDP:theorem}. In the end we will use the union bound.

    Suppose we have an estimator such that $\|\hat \beta_{j} - \beta\|_2 \leq \tilde d_{j}$. Consider the set $K_{(j)} =  K \cap B_2(0, \tilde d_{j}/2)$ and compute $m$ as in Corollary \ref{Kwidth:bound:lemma:regression}. Denote the corresponding matrix $A_j^{\star\star} = (I - X_j^{\star\star})^{1/2}$. Let $X_j^{\star\star} = V\Lambda V\T$ be the eigendecomposition of $X_j^{\star\star}$.

We consider the estimate $\beta'_{j + 1}$ defined in \eqref{strange:optimization:program}.
We will argue that by Theorem 2.5.9 of \cite{dadush2012integer} (see also Theorem 1 in \cite{lee2018efficient}) that $\beta'_{j + 1}$ is computable. First, we note that the set $K_{(j)}$ admits a membership oracle since we can evaluate its Minkowski gauge. Furthermore, the convex map $\nu \mapsto \|((Y_i - Z_i\T \hat \beta_j) - Z_i\T A^{\star\star}\nu)_{i \in [N]}\|_2$ is Lipschitz. We have 
\begin{align*}
    \MoveEqLeft |\|((Y_i - Z_i\T \hat \beta_j) - Z_i\T A^{\star\star}\nu)_{i \in [N]}\|_2 - \|((Y_i - Z_i\T \hat \beta_j) - Z_i\T A^{\star\star}\nu')_{i \in [N]}\|_2|\\
    & \leq \|(Z_i\T A^{\star\star} (\nu - \nu'))_{i \in [N]}\|_2 \\
    & \leq \sqrt{ \lambda_{\max}\bigg(A^{\star\star} \sum_{i \in [N]} Z_i Z_i\T A^{\star\star} \bigg)} \|\nu - \nu'\|_2\\
    & \leq \sqrt{ \lambda_{\max}\bigg(\sum_{i \in [N]} Z_i Z_i\T }\bigg) \|\nu - \nu'\|_2
\end{align*}
We have by Theorem 5.39 of \cite{vershynin2010introduction}:
\begin{align*}
    \sqrt{ \lambda_{\max}\bigg(\sum_{i \in [N]}\Sigma^{-1/2} Z_i Z_i \T\Sigma^{-1/2}\bigg) } \leq \sqrt{N} + C\sqrt{n} + t,
\end{align*}
with probability at least $1 - 2 \exp(-c t^2)$. Recall that $\lambda_{\max}(\Sigma) = O(1)$. Setting $t = \sqrt{N}\tilde d_{j}$ we obtain that with probability at least $1 - 2 \exp(-c N \tilde d_{j}^2)$
\begin{align}\label{lambdamaxSigma}
    \sqrt{ \lambda_{\max}\bigg(\sum_{i \in [N]} Z_i Z_i\T \bigg)} \lesssim \sqrt{N} + C\sqrt{n} + \sqrt{N}\tilde d_{j},
\end{align}
and thus the map $\nu \mapsto \|((Y_i - Z_i\T \hat \beta_j) - Z_i\T A^{\star\star}\nu)_{i \in [N]}\|_2$ is Lipschitz with a constant $\lesssim  (\sqrt{N} + c \sqrt{n} + \sqrt{N}\tilde d_{j}) \asymp (\sqrt{N} + \sqrt{n})$, since $\tilde d_{j} = O(1)$ by assumption. Hence the above approximate minimization over $K_{(j)}$ is computable by Theorem 2.5.9 in \cite{dadush2012integer}.

Let $K'_{(j)} = \{ A_j^{\star\star} \nu : \nu \in K_{(j)}\}$. Denote with $\nu_{j} := (\beta - \hat \beta_{j} )/2 \in K_{(j)}$ by assumption. Let $\Pi_{K'_{(j)}}\nu_j = \argmin_{\nu \in K'_{(j)}} \|\nu - \nu_j\|_2$. We know that $\|\Pi_{K'_{(j)}}\nu_j - \nu_j\|_2 \leq \|A_j^{\star\star}\nu_j - \nu_j\|_2 \lesssim \tilde d_{j}/c$ where the last inequality follows from a similar argument as in the proof of Lemma \ref{sub:gaussian:vector:estimator}. Next using Theorem 2.5.9 \cite{dadush2012integer} with $\epsilon \asymp \sqrt{N} \tilde d_j^{2}$ we can obtain a near optimal solution $\beta'_{j+1}$. Then by the (near) optimality we have:
\begin{align}\label{least:squares:equations}
    \sum_{i \in [N]} ((Y_i - Z_i\T\hat \beta_{j})/2 - Z_i\T A_j^{\star\star} \beta'_{j+1})^2 & \leq \sum_{i \in [N]} ((Y_i - Z_i\T\hat \beta_{j})/2 - Z_i\T \Pi_{K'_{(j)}}\nu_j)^2 + \epsilon^2 \nonumber \\
    & + 2\|((Y_i - Z_i\T\hat \beta_{j})/2 - Z_i\T \Pi_{K'_{(j)}}\nu_j)_{i \in [N]}\|_2\epsilon.
\end{align}

Next using the shorthand $\hat \Sigma := N^{-1}\sum_{i \in [N]} Z_i Z_i\T$, we handle the last term of the right hand side:
\begin{align*}
    \|((Y_i - Z_i\T\hat \beta_{j})/2 - Z_i\T \Pi_{K'_{(j)}}\nu_j)_{i \in [N]}\|_2 \leq \sqrt{N(\Pi_{K'_{(j)}}\nu_j - \nu_j)\T \hat \Sigma (\Pi_{K'_{(j)}}\nu_j - \nu_j)} + \|\xi\|_2.
\end{align*}
Using \eqref{lambdamaxSigma}, the first term can be bounded as $\sqrt{N (\Pi_{K'_{(j)}}\nu_j - \nu_j)\T \hat \Sigma (\Pi_{K'_{(j)}}\nu_j - \nu_j)}\lesssim (\sqrt{N} + \sqrt{n}) \tilde d_j/c$. On the other hand, by Proposition 5.16 of \cite{vershynin2010introduction} we know that
\begin{align*}
    \PP(\|\xi\|^2_2 - \EE \|\xi\|^2_2 \geq t) \leq 2\exp(-C'' (t^2/N \wedge t)).
\end{align*}
Hence we can select $t = N\tilde d_j$ to obtain that with probability at least $1 - 2\exp(-C'' N \tilde d_j^2)$ we have $\|\xi\|^2_2 \leq \EE \|\xi\|_2^2 + N \tilde d_j \leq N \gamma^2 + N \tilde d_j \lesssim N$. Hence since $\epsilon \asymp \sqrt{N}\tilde d_j^{2}$, and $\tilde d_j = O(1)$ we have that:
\begin{align*}
    \epsilon^2 + 2\|((Y_i - Z_i\T\hat \beta_{j})/2 - Z_i\T \Pi_{K'_{(j)}}\nu_j)_{i \in [N]}\|_2\epsilon = O(N \tilde d_j^2).
\end{align*}
with high probability.

Then by rearranging \eqref{least:squares:equations} we have
\begin{align*}
    \MoveEqLeft (\tilde \beta_{j+1} - \nu_j)\T\hat \Sigma (\tilde \beta_{j+1} - \nu_j) \\
    & \leq  (\Pi_{K'_{(j)}}\nu_j - \nu_j)\T \hat \Sigma (\Pi_{K'_{(j)}}\nu_j - \nu_j)  + 2\bigg|N^{-1} \sum_{i \in [N]} \xi_i Z_i\T (\Pi_{K'_{(j)}}\nu_{j} - \tilde\beta_{j + 1})\bigg| + \tilde d_{j}^2/C^{'2},
\end{align*}
where $\tilde \beta_{j + 1} = A_j^{\star\star} \beta'_{j + 1}$. Since $\sqrt{a + b + c} \leq \sqrt{a} + \sqrt{b} + \sqrt{c}$ for any $a, b, c \geq 0$ we have
\begin{align}\label{Sigma:hat:bound}
    \|\hat \Sigma^{1/2}(\tilde \beta_{j + 1} - \nu_j)\|_2 \leq \|\hat \Sigma^{1/2}(\Pi_{K'_{(j)}}\nu_{j} - \nu_j)\|_2 + \sqrt{\bigg|2 N^{-1} \sum_{i \in [N]} \xi_i Z_i\T (\Pi_{K'_{(j)}}\nu_{j} - \tilde \beta_{j + 1})\bigg|} + \tilde d_j/C'.
\end{align}
We first control the second term on the right hand side. Let $\tilde \nu_j \in K_{(j)}$, so that $\Pi_{K'_{(j)}} \nu_j = A^{\star\star}_j\tilde \nu_j$. Also, let $A_j^{\star \star} = \sum_{k \in [n]} \sqrt{1-\lambda_k} v_k v_k\T $ be the eigendecomposition of $A_j^{\star\star}$.

We now control the second term on the right-hand side of \eqref{Sigma:hat:bound}. Let $\tilde \nu_j \in K_{(j)}$ be such that
\begin{align*}
\Pi_{K'_{(j)}}\nu_j=A_j^{\star\star}\tilde\nu_j,
\end{align*}
and recall that $\tilde \beta_{j+1}=A_j^{\star\star}\beta'_{j+1}$, where $\beta'_{j+1}\in K_{(j)}$. Hence
\begin{align*}
\|\tilde\nu_j-\beta'_{j+1}\|_2\leq \tilde d_j,
\end{align*}
since $K_{(j)}\subseteq B_2(0,\tilde d_j/2)$.

For notational simplicity, write $A=A_j^{\star\star}$ and define
\begin{align*}
S:=\sum_{i\in[N]}\xi_i Z_i.
\end{align*}
Then
\begin{align}
\left|
\sum_{i\in[N]}\xi_i Z_i\T
\left(
\Pi_{K'_{(j)}}\nu_j-\tilde\beta_{j+1}
\right)
\right|
&=
\left|
\left\langle
S,A(\tilde\nu_j-\beta'_{j+1})
\right\rangle
\right|
\nonumber\\
&\leq
\|AS\|_2\,
\|\tilde\nu_j-\beta'_{j+1}\|_2
\nonumber\\
&\leq
\tilde d_j\|AS\|_2.
\label{eq:regression-noise-cauchy}
\end{align}
We therefore bound $\|AS\|_2$ directly, rather than bounding its coordinates separately.

Let $\mathbf Z\in\RR^{N\times n}$ denote the design matrix whose $i$-th row is $Z_i\T$, and let
\begin{align*}
\xi=(\xi_1,\ldots,\xi_N)\T.
\end{align*}
Then
\begin{align*}
AS=A\mathbf Z\T\xi,
\end{align*}
so that
\begin{align}
\|AS\|_2^2
=
\xi\T B\xi,
\qquad
B:=\mathbf Z A^2\mathbf Z\T.
\label{eq:regression-noise-quadratic}
\end{align}
We condition on the design matrix $\mathbf Z$ and on the randomness used to construct $A$. Since the noise variables $\xi_i$ are independent of the design, conditional on $\mathbf Z$ the vector $\xi$ continues to have independent, centered, sub-Gaussian entries with sub-Gaussian parameter of constant order.

Recall that
\begin{align*}
A^2=I-X_j^{\star\star},
\qquad
0\leq A^2\leq I,
\end{align*}
and
\begin{align*}
\tr(A^2)=m,
\end{align*}
where
\begin{align*}
m=\left\lceil \frac{N\tilde d_j^2}{C^2}\right\rceil\wedge n.
\end{align*}
On the event
\begin{align*}
\|\hat\Sigma\|_{\operatorname{op}}\lesssim 1,
\qquad
\hat\Sigma=\frac1N\mathbf Z\T\mathbf Z,
\end{align*}
which, as shown above, occurs with probability at least
\begin{align*}
1-\exp(-cN\tilde d_j^2),
\end{align*}
we have
\begin{align}
\tr(B)
&=
\tr(A^2\mathbf Z\T\mathbf Z)
=
N\tr(A^2\hat\Sigma)
\nonumber\\
&\leq
N\|\hat\Sigma\|_{\operatorname{op}}\tr(A^2)
\lesssim Nm.
\label{eq:trace-B-regression}
\end{align}
Moreover,
\begin{align}
\|B\|_{\operatorname{op}}
&=
\|\mathbf Z A^2\mathbf Z\T\|_{\operatorname{op}}
\leq
\|\mathbf Z\|_{\operatorname{op}}^2
\|A^2\|_{\operatorname{op}}
\nonumber\\
&\leq
\|\mathbf Z\T\mathbf Z\|_{\operatorname{op}}
=
N\|\hat\Sigma\|_{\operatorname{op}}
\lesssim N.
\label{eq:op-B-regression}
\end{align}
Since $B$ is positive semidefinite,
\begin{align}
\|B\|_F^2
\leq
\|B\|_{\operatorname{op}}\tr(B)
\lesssim N^2m.
\label{eq:frob-B-regression}
\end{align}

We now apply the Hanson--Wright inequality \cite{rudelson2013hanson} conditionally on $\mathbf Z$ and $A$. Since
\begin{align*}
\EE_\xi[\xi\T B\xi\mid \mathbf Z,A]
=
\EE\xi_1^2\,\tr(B)
\lesssim Nm,
\end{align*}
we obtain
\begin{align}
\PP_\xi\left(
\xi\T B\xi
\geq
C_1Nm+t
\;\middle|\;
\mathbf Z,A
\right)
\leq
2\exp\left(
-c
\left[
\frac{t^2}{N^2m}
\wedge
\frac{t}{N}
\right]
\right)
\label{eq:HW-regression-projected-noise}
\end{align}
for suitable absolute constants $c,C_1>0$.

By assumption $\tilde d_j\geq C/\sqrt N$, and therefore
\begin{align*}
m
=
\left\lceil\frac{N\tilde d_j^2}{C^2}\right\rceil\wedge n
\lesssim
\frac{N\tilde d_j^2}{C^2}.
\end{align*}
Consequently,
\begin{align*}
Nm\lesssim\frac{N^2\tilde d_j^2}{C^2}.
\end{align*}
Taking
\begin{align*}
t=C_2\frac{N^2\tilde d_j^2}{C^2}
\end{align*}
for a sufficiently large absolute constant $C_2$, \eqref{eq:HW-regression-projected-noise} gives
\begin{align}
\PP_\xi\left(
\|AS\|_2^2
\gtrsim
\frac{N^2\tilde d_j^2}{C^2}
\;\middle|\;
\mathbf Z,A
\right)
\leq
2\exp(-c' N\tilde d_j^2),
\end{align}
where, as usual, constants depending on the sufficiently large tuning constant $C$ are absorbed into $c'>0$. Thus, with probability at least
\begin{align*}
1-2\exp(-c'N\tilde d_j^2),
\end{align*}
we have
\begin{align}
\left\|
A_j^{\star\star}
\sum_{i\in[N]}\xi_i Z_i
\right\|_2
\lesssim
\frac{N\tilde d_j}{C}.
\label{eq:projected-regression-noise}
\end{align}

Combining \eqref{eq:regression-noise-cauchy} and
\eqref{eq:projected-regression-noise}, we conclude that
\begin{align}
\left|
\sum_{i\in[N]}
\xi_i Z_i\T
\left(
\Pi_{K'_{(j)}}\nu_j-\tilde\beta_{j+1}
\right)
\right|
\lesssim
\frac{N\tilde d_j^2}{C}
\label{eq:final-regression-noise}
\end{align}
with probability at least $1-\exp(-c'N\tilde d_j^2)$, up to a change in the absolute constant $c'$.

Therefore,
\begin{align}
\sqrt{
\left|
2N^{-1}
\sum_{i\in[N]}
\xi_i Z_i\T
\left(
\Pi_{K'_{(j)}}\nu_j-\tilde\beta_{j+1}
\right)
\right|
}
\lesssim
\frac{\tilde d_j}{\sqrt C}.
\end{align}
Since $C$ can be chosen to be a sufficiently large absolute constant, this term can be made smaller than $\tilde d_j/T'$ for any prescribed absolute constant $T'>0$.

Now, by Theorem 5.39 of \cite{vershynin2010introduction}, applied to the isotropic sub-Gaussian vectors
$\Sigma^{-1/2}Z_i$, whose sub-Gaussian norms are bounded by an absolute constant since
$\lambda_{\min}(\Sigma)=\Omega(1)$, we have
\begin{align*}
1-C\sqrt{\frac nN}-\frac{t}{\sqrt N}
&\leq
\sqrt{\lambda_{\min}
\left(
\Sigma^{-1/2}\widehat\Sigma\Sigma^{-1/2}
\right)}
\\
&\leq
\sqrt{\lambda_{\max}
\left(
\Sigma^{-1/2}\widehat\Sigma\Sigma^{-1/2}
\right)}
\leq
1+C\sqrt{\frac nN}+\frac{t}{\sqrt N}
\end{align*}
with probability at least $1-2\exp(-ct^2)$.

Take
\begin{align*}
t=\frac{\sqrt N\,\widetilde d_j}{B},
\end{align*}
where $B>0$ is a sufficiently large absolute constant. Since
$\widetilde d_j\leq 2R=O(1)$, and since $N\gtrsim n$ with a sufficiently
large absolute implicit constant, we may choose the constants so that
\begin{align*}
C\sqrt{\frac nN}+\frac{\widetilde d_j}{B}\leq \frac12.
\end{align*}
Consequently, with probability at least
\begin{align*}
1-2\exp(-c'N\widetilde d_j^2),
\end{align*}
we have
\begin{align*}
1
\asymp
\lambda_{\min}
\left(
\Sigma^{-1/2}\widehat\Sigma\Sigma^{-1/2}
\right)
\leq
\lambda_{\max}
\left(
\Sigma^{-1/2}\widehat\Sigma\Sigma^{-1/2}
\right)
\asymp
1.
\end{align*}
Equivalently,
\begin{align*}
c_0\Sigma
\preceq
\widehat\Sigma
\preceq
C_0\Sigma
\end{align*}
for some absolute constants $c_0,C_0>0$. Therefore, for every
$v\in\mathbb R^n$,
\begin{align*}
c_0^{1/2}\|\Sigma^{1/2}v\|_2
\leq
\|\widehat\Sigma^{1/2}v\|_2
\leq
C_0^{1/2}\|\Sigma^{1/2}v\|_2.
\end{align*}
Combining this with \eqref{Sigma:hat:bound} and the preceding bounds gives
\begin{align*}
\lambda_{\min}(\Sigma)^{1/2}
\|\widetilde\beta_{j+1}-\nu_j\|_2
&\lesssim
\lambda_{\max}(\Sigma)^{1/2}
\|\Pi_{K'_{(j)}}\nu_j-\nu_j\|_2
+
\frac{\widetilde d_j}{T'}
+
\frac{\widetilde d_j}{C'}
\\
&\leq
\frac{\widetilde d_j}{L},
\end{align*}
after choosing the absolute constants sufficiently large.
Since $\lambda_{\min}(\Sigma)=\Omega(1)$, it follows that
\begin{align*}
\|\widetilde\beta_{j+1}-\nu_j\|_2
\leq
\frac{\widetilde d_j}{L}.
\end{align*}

Combining, our previous inequalities with \eqref{Sigma:hat:bound}, we can claim that
\begin{align*}
\lambda_{\min}(\Sigma)^{1/2}
\|\tilde\beta_{j+1}-\nu_j\|_2
&\lesssim
\lambda_{\max}(\Sigma)^{1/2}
\|\Pi_{K'_{(j)}}\nu_j-\nu_j\|_2
+
\frac{\tilde d_j}{T'}
+
\frac{\tilde d_j}{C'}\\
&\leq
\frac{\tilde d_j}{L},
\end{align*}
for some constants $T$ and $T'$ that depend on $C, \lambda_{\max}(\Sigma) = O(1), c, n/N$ and other constants defined in this proof. Thus since $\lambda_{\min}(\Sigma)^{1/2}  = \Omega(1)$ we have
\begin{align*}
    \|\tilde \beta_{j + 1} - \nu_j\|_2 \leq \tilde d_{j}/L,
\end{align*}
for some constant $L$ which can be made arbitrarily large by making the constant $C$ large. Upon grouping all of our probability bounds we obtain that this happens with probability at least $1 - \exp(- c' N \tilde d_{j}^2)$ for some constant $c'$. The union bound completes the proof. 
\end{proof}

\begin{proof}[Proof of Theorem \ref{regression:theorem:main}]
Let
\begin{align*}
    \rho
    :=
    \frac{2(\sqrt{3}+1)}{\widetilde L}
    <1.
\end{align*}
Thus, as long as the algorithm continues,
\begin{align*}
    \widetilde d_j
    =
    2R\rho^{j-1},
    \qquad
    \widetilde d_{j+1}
    =
    \rho\widetilde d_j.
\end{align*}
We write
\begin{align*}
    \eta^\star:=\eta^\star(K)
    =
    \sup\left\{
        \eta\geq0:
        N\eta^2\leq \log M_K^{\operatorname{loc}}(\eta)
    \right\},
\end{align*}
and put
\begin{align*}
    \overline\kappa
    :=
    \max_{j\in[M]}\kappa(K_{(j)}).
\end{align*}
Recall that
\begin{align*}
    K_{(j)}
    =
    K\cap B_2(0,\widetilde d_j/2).
\end{align*}
As in the proof of Theorem~\ref{main:result:GSM},
\begin{align*}
    T_2(K_{(j)})
    \leq
    \sqrt2\,T_2(K),
    \qquad
    \eta^\star(K_{(j)})
    \leq
    \eta^\star(K).
\end{align*}
The second inequality follows because $K_{(j)}\subseteq K$, and hence
\begin{align*}
    M_{K_{(j)}}^{\operatorname{loc}}(\delta)
    \leq
    M_K^{\operatorname{loc}}(\delta)
\end{align*}
for every $\delta>0$.

Define
\begin{align*}
    \Lambda_K
    :=
    1+
    \left\lceil
        \log_c\!\left(
            \overline\kappa^{1/2}\sqrt2\,T_2(K)
        \right)
    \right\rceil
\end{align*}
and
\begin{equation}\label{eq:regression-GammaK}
    \Gamma_K
    :=
    1
    \vee
    \left(
        2\sqrt2\,T_2(K)\overline\kappa^{1/2}
    \right)
    \vee
    \left(
        2C\sqrt{\Lambda_K}
    \right).
\end{equation}
All absolute constants below may be enlarged, so that harmless constants
coming from the weak projections and from
Lemma~\ref{sub:gaussian:vector:estimator:regression} are absorbed into
$\Gamma_K$.

We first record two elementary consequences of the entropy definition that
will be useful below. Let $d:=\diam(K)$. Since $K$ is convex and
origin-symmetric, it contains a line segment of length arbitrarily close to
$d$. Therefore, for a sufficiently small absolute constant $a>0$, setting
\begin{align*}
    \eta_0
    :=
    a\left(d\wedge N^{-1/2}\right)
\end{align*}
gives
\begin{align*}
    \log M_K^{\operatorname{loc}}(\eta_0)
    \geq
    \log 2,
    \qquad
    N\eta_0^2
    \leq
    \log 2.
\end{align*}
Consequently,
\begin{equation}\label{eq:regression-parametric-floor}
    \eta^\star(K)
    \gtrsim
    d\wedge N^{-1/2}.
\end{equation}
In particular, if $d\lesssim N^{-1/2}$, then the conclusion of the theorem is
immediate: every iterate returned by Algorithm~\ref{algorithm:main:regression}
belongs to $K$, and hence
\begin{align*}
    \|\widehat\beta-\beta\|_2
    \leq d
    \lesssim \eta^\star(K).
\end{align*}
We may therefore work throughout the remainder of the proof in the
nontrivial regime
\begin{equation}\label{eq:regression-nontrivial-diameter}
    d\gtrsim N^{-1/2},
    \qquad\text{and hence}\qquad
    \eta^\star(K)\gtrsim N^{-1/2}.
\end{equation}

We also need to justify the $r$-part of the stopping rule. Denote by $r_0$
the original inner radius and by
\begin{align*}
    r_{\rm alg}
    :=
    r_0\wedge \frac12\sqrt{\frac nN}
\end{align*}
the value after the overwrite in Algorithm~\ref{algorithm:main:regression}.
Then $r_{\rm alg}B_2\subseteq K$ and
\begin{align*}
    Nr_{\rm alg}^2\leq n/4.
\end{align*}
For a sufficiently small absolute constant $a_0>0$, the Euclidean ball
$B_2(0,a_0r_{\rm alg})\subseteq K$ admits an
$(a_0r_{\rm alg}/c)$-packing of cardinality at least $\exp(c_0n)$ for an
absolute constant $c_0>0$. Hence
\begin{align*}
    \log M_K^{\operatorname{loc}}(a_0r_{\rm alg})
    \gtrsim n
    \gtrsim
    N(a_0r_{\rm alg})^2,
\end{align*}
where $a_0$ is chosen sufficiently small. By the definition of $\eta^\star$,
\begin{equation}\label{eq:regression-r-below-eta}
    r_{\rm alg}
    \lesssim
    \eta^\star(K).
\end{equation}
Notice that this argument uses the overwritten radius. We do not claim that
the minimax risk is always bounded below by $r_0^2$; rather, the overwrite
ensures precisely that the Euclidean-ball local entropy at scale
$r_{\rm alg}$ is large enough for~\eqref{eq:regression-r-below-eta}.

We now make the inductive probability argument explicit. Define the trapping
events
\begin{align*}
    H_j
    :=
    \left\{
        \|\widehat\beta_j-\beta\|_2
        \leq
        \widetilde d_j
    \right\}.
\end{align*}
The initial event $H_1$ holds deterministically, since
$\widehat\beta_1=0$ and
\begin{align*}
    \|\beta\|_2
    \leq R
    \leq \widetilde d_1=2R.
\end{align*}
Whenever the two assumptions of
Corollary~\ref{Kwidth:bound:lemma:regression} hold for $K_{(j)}$ and
\begin{equation}\label{eq:regression-operating-scale}
    \widetilde d_j
    \geq
    (2r_{\rm alg})\vee(CN^{-1/2}),
\end{equation}
define the good refinement event
\begin{align*}
    E_j
    :=
    \left\{
        \left\|
            \widetilde\beta_{j+1}
            -
            \frac{\beta-\widehat\beta_j}{2}
        \right\|_2
        \leq
        \frac{\widetilde d_j}{L}
    \right\}.
\end{align*}
Lemma~\ref{sub:gaussian:vector:estimator:regression} gives
\begin{equation}\label{eq:regression-refinement-prob}
\begin{split}
    \PP(H_j\cap E_j^c)
    \leq
    \exp(-c_1N\widetilde d_j^2)
    +
    \frac{1}{1+N^{3/2}\widetilde d_j^3}
\end{split}
\end{equation}
for an absolute constant $c_1>0$. We stress that this is an unconditional
bound on the intersection $H_j\cap E_j^c$; no independence between
iterations is required.

We next verify that, on $H_j\cap E_j$, the same deterministic contraction as
in the proof of Theorem~\ref{main:result:GSM} holds. On $H_j$,
\begin{align*}
    \nu_j
    :=
    \frac{\beta-\widehat\beta_j}{2}
    \in
    K\cap B_2(0,\widetilde d_j/2)
    =
    K_{(j)}.
\end{align*}
Indeed, $\beta,-\widehat\beta_j\in K$ and $K$ is convex and
origin-symmetric. Applying Lemma~\ref{weak:projection:lemma} to the weak
projection of $\widetilde\beta_{j+1}$ onto $K_{(j)}$, with
$\epsilon=\widetilde d_j/L$, gives on $H_j\cap E_j$
\begin{align*}
\begin{split}
    \|\overline\beta_{j+1}-\nu_j\|_2
    &\leq
    \|\widetilde\beta_{j+1}-\nu_j\|_2
    +
    \sqrt{
        \epsilon^2
        +2\epsilon
        \|\widetilde\beta_{j+1}-\nu_j\|_2
    }\\
    &\leq
    \frac{(1+\sqrt3)\widetilde d_j}{L}
    =
    \frac{\widetilde d_j}{\widetilde L}.
\end{split}
\end{align*}
Consequently,
\begin{align*}
    \|2\overline\beta_{j+1}+\widehat\beta_j-\beta\|_2
    \leq
    \frac{2\widetilde d_j}{\widetilde L}.
\end{align*}
A second application of Lemma~\ref{weak:projection:lemma}, now to the weak
projection onto $K$ with
$\epsilon=2\widetilde d_j/\widetilde L$, yields
\begin{align*}
\begin{split}
    \|\widehat\beta_{j+1}-\beta\|_2
    &\leq
    \frac{2\widetilde d_j}{\widetilde L}
    +
    \sqrt{
        \frac{4\widetilde d_j^2}{\widetilde L^2}
        +
        \frac{8\widetilde d_j^2}{\widetilde L^2}
    }\\
    &=
    \frac{2(\sqrt3+1)\widetilde d_j}{\widetilde L}
    =
    \widetilde d_{j+1}.
\end{split}
\end{align*}
Thus, whenever the deterministic assumptions of the width corollary and
\eqref{eq:regression-operating-scale} continue to hold,
\begin{equation}\label{eq:regression-good-implies-next}
    H_j\cap E_j
    \subseteq
    H_{j+1}.
\end{equation}

As in the proof of Theorem~\ref{main:result:GSM}, we also need a deterministic
stability statement after the first scale at which the contraction argument
can no longer be invoked. Let
\begin{align*}
    a_\ell
    :=
    \|\widehat\beta_\ell-\beta\|_2.
\end{align*}
Suppose $a_k\leq\widetilde d_k$. Regardless of whether the probabilistic
refinement guarantee holds at later iterations,
$\overline\beta_\ell\in K_{(\ell-1)}$, so
\begin{align*}
    2\|\overline\beta_\ell\|_2
    \leq
    \widetilde d_{\ell-1}.
\end{align*}
Using Lemma~\ref{weak:projection:lemma} for the lift back to $K$, with
$\epsilon_{\ell-1}=2\widetilde d_{\ell-1}/\widetilde L$, gives
\begin{align*}
\begin{split}
    a_\ell
    &\leq
    a_{\ell-1}
    +
    \widetilde d_{\ell-1}
    +
    \sqrt{
        \epsilon_{\ell-1}^2
        +
        2\epsilon_{\ell-1}
        (\widetilde d_{\ell-1}+a_{\ell-1})
    }\\
    &\leq
    a_{\ell-1}
    +
    C_2\widetilde d_{\ell-1}
    +
    C_3\sqrt{\widetilde d_{\ell-1}a_{\ell-1}},
\end{split}
\end{align*}
where $C_2,C_3$ are absolute constants depending only on
$\widetilde L$. If
\begin{align*}
    A_k
    :=
    \max_{\ell\geq k}a_\ell,
\end{align*}
then, using
$\widetilde d_\ell=\rho^{\ell-k}\widetilde d_k$ and summing the preceding
inequality over the remaining iterations,
\begin{align*}
\begin{split}
    A_k
    &\leq
    \widetilde d_k
    +
    C_2\sum_{\ell\geq k}\widetilde d_\ell
    +
    C_3\sqrt{A_k}
        \sum_{\ell\geq k}\sqrt{\widetilde d_\ell}\\
    &\leq
    C_4\widetilde d_k
    +
    C_5\sqrt{\widetilde d_kA_k}.
\end{split}
\end{align*}
Solving the resulting quadratic inequality in
$\sqrt{A_k/\widetilde d_k}$ gives
\begin{equation}\label{eq:regression-post-failure-stability}
    \max_{\ell\geq k}
    \|\widehat\beta_\ell-\beta\|_2
    \leq
    C_6\widetilde d_k
\end{equation}
for an absolute constant $C_6$. Hence, after the first unsuccessful
refinement, or after the first deterministic scale at which the width
argument ceases to apply, all later iterates remain within a constant multiple
of the last valid scale.

We now identify the deterministic scales at which the induction can stop.
First, suppose the first condition of
Corollary~\ref{Kwidth:bound:lemma:regression} fails for $K_{(j)}$. Then
\begin{align*}
    \kappa(K_{(j)})^{-1/2}
    \frac{\widetilde d_j}{T_2(K_{(j)})}
    <
    2\eta^\star(K_{(j)}),
\end{align*}
and therefore
\begin{align*}
    \widetilde d_j
    <
    2\kappa(K_{(j)})^{1/2}
    T_2(K_{(j)})
    \eta^\star(K_{(j)})
    \leq
    2\sqrt2\,
    \overline\kappa^{1/2}T_2(K)\eta^\star(K).
\end{align*}
Thus
\begin{equation}\label{eq:regression-first-stop}
    \widetilde d_j
    \leq
    \Gamma_K\eta^\star(K).
\end{equation}

Second, suppose the second condition of
Corollary~\ref{Kwidth:bound:lemma:regression} fails. Set
\begin{align*}
    L_j
    :=
    \left\lceil
        \log_c\!\left(
            \kappa(K_{(j)})^{1/2}T_2(K_{(j)})
        \right)
    \right\rceil.
\end{align*}
Then
\begin{align*}
    \frac{N\widetilde d_j^2}{C^2}
    <
    \left\lceil
        4L_jN\eta^{*2}(K_{(j)})
    \right\rceil,
\end{align*}
so, using $\lceil x\rceil\leq x+1$,
\begin{align*}
    \widetilde d_j^2
    \lesssim
    C^2\left(
        L_j\eta^{*2}(K_{(j)})
        +
        \frac1N
    \right).
\end{align*}
Since $L_j\leq\Lambda_K$,
$\eta^\star(K_{(j)})\leq\eta^\star(K)$, and
$N^{-1/2}\lesssim\eta^\star(K)$ by
\eqref{eq:regression-nontrivial-diameter}, we conclude that
\begin{equation}\label{eq:regression-second-stop}
    \widetilde d_j
    \lesssim
    C\sqrt{\Lambda_K}\,\eta^\star(K)
    \lesssim
    \Gamma_K\eta^\star(K).
\end{equation}

Third, if the operating-scale condition
\eqref{eq:regression-operating-scale} ceases to hold, then by
\eqref{eq:regression-parametric-floor} and
\eqref{eq:regression-r-below-eta},
\begin{equation}\label{eq:regression-third-stop}
    \widetilde d_j
    \lesssim
    r_{\rm alg}\vee N^{-1/2}
    \lesssim
    \eta^\star(K)
    \leq
    \Gamma_K\eta^\star(K).
\end{equation}
The same conclusion holds if the algorithm terminates because
$\widetilde d_{j+1}\leq(2r_{\rm alg})\vee(CN^{-1/2})$; indeed
$\widetilde d_j=\rho^{-1}\widetilde d_{j+1}$ and $\rho$ is an absolute
constant. Finally, if all $M$ iterations are completed, then the definition
of $M$ gives
\begin{align*}
    \widetilde d_{M+1}
    \leq
    2r_{\rm alg},
\end{align*}
so the same conclusion follows from
\eqref{eq:regression-r-below-eta}.

We can now accumulate the failure probabilities. Let
\begin{align*}
    t_0
    :=
    C_7\Gamma_K\eta^\star(K),
\end{align*}
where $C_7$ is a sufficiently large absolute constant, depending only on the
constant in~\eqref{eq:regression-post-failure-stability}. By
\eqref{eq:regression-first-stop}--\eqref{eq:regression-third-stop}, a purely
deterministic termination of the contraction argument cannot lead to final
error exceeding $t$ for $t\geq t_0$. Hence, if
\begin{align*}
    \|\widehat\beta-\beta\|_2\geq t,
    \qquad t\geq t_0,
\end{align*}
then there must be a first iteration $j$ for which $H_j$ holds but $E_j$
fails. By~\eqref{eq:regression-post-failure-stability}, this iteration must
also satisfy
\begin{align*}
    \widetilde d_j
    \geq
    t/C_6.
\end{align*}
Consequently, the union bound and
\eqref{eq:regression-refinement-prob} imply
\begin{equation}\label{eq:regression-explicit-union}
\begin{split}
    \PP\left(
        \|\widehat\beta-\beta\|_2\geq t
    \right)
    \leq
    \sum_{j:\,\widetilde d_j\geq t/C_6}
    \left[
        \exp(-c_1N\widetilde d_j^2)
        +
        \frac{1}{1+N^{3/2}\widetilde d_j^3}
    \right].
\end{split}
\end{equation}
This is the regression analogue of the accumulated bound in the proof of
Theorem~\ref{main:result:GSM}.

It remains to simplify the geometric sum. Let $j_t$ be the largest index in
\eqref{eq:regression-explicit-union}. Since
$\widetilde d_{j+1}=\rho\widetilde d_j$, for $s\geq0$,
\begin{align*}
    \widetilde d_{j_t-s}
    =
    \rho^{-s}\widetilde d_{j_t}
    \gtrsim
    \rho^{-s}t.
\end{align*}
Because $t\geq t_0\gtrsim\eta^\star(K)\gtrsim N^{-1/2}$ in the regime under
consideration,
\begin{align*}
\begin{split}
    \sum_{j:\,\widetilde d_j\gtrsim t}
    \exp(-c_1N\widetilde d_j^2)
    &\leq
    \sum_{s=0}^{\infty}
    \exp(-c_2\rho^{-2s}Nt^2)\\
    &\lesssim
    \exp(-c_3Nt^2).
\end{split}
\end{align*}
Similarly,
\begin{align*}
\begin{split}
    \sum_{j:\,\widetilde d_j\gtrsim t}
    \frac{1}{1+N^{3/2}\widetilde d_j^3}
    &\leq
    \frac{1}{N^{3/2}}
    \sum_{j:\,\widetilde d_j\gtrsim t}
    \frac{1}{\widetilde d_j^3}\\
    &\lesssim
    \frac{1}{N^{3/2}t^3}.
\end{split}
\end{align*}
Thus, for every $t\geq t_0$,
\begin{equation}\label{eq:regression-high-probability-bound}
    \PP\left(
        \|\widehat\beta-\beta\|_2\geq t
    \right)
    \lesssim
    \exp(-c_3Nt^2)
    +
    \frac{C_8}{N^{3/2}t^3}.
\end{equation}
The second term is exactly the accumulated failure probability coming from
the randomized QFM oracle.

We finally integrate the tail. Since the estimator is proper,
$\|\widehat\beta-\beta\|_2\leq d$, so
\eqref{eq:regression-high-probability-bound} may harmlessly be integrated to
infinity. Using
\begin{align*}
    \EE\|\widehat\beta-\beta\|_2^2
    =
    2\int_0^\infty
    t\,
    \PP(\|\widehat\beta-\beta\|_2\geq t)\,dt,
\end{align*}
the trivial bound below $t_0$ and
\eqref{eq:regression-high-probability-bound} above $t_0$ give
\begin{align*}
\begin{split}
    \EE\|\widehat\beta-\beta\|_2^2
    &\lesssim
    t_0^2
    +
    \int_{t_0}^{\infty}t\exp(-c_3Nt^2)\,dt
    +
    \frac{1}{N^{3/2}}
    \int_{t_0}^{\infty}\frac{dt}{t^2}\\
    &\lesssim
    \Gamma_K^2\eta^{*2}(K)
    +
    \frac1N\exp(-c_4Nt_0^2)
    +
    \frac{1}{N^{3/2}t_0}.
\end{split}
\end{align*}
By~\eqref{eq:regression-nontrivial-diameter},
$N^{-1/2}\lesssim\eta^\star(K)$. Therefore
\begin{align*}
    \frac1N
    \lesssim
    \eta^{*2}(K),
    \qquad
    \frac{1}{N^{3/2}t_0}
    \lesssim
    \frac{1}{N^{3/2}\eta^\star(K)}
    \lesssim
    \eta^{*2}(K).
\end{align*}
It follows that
\begin{align*}
    \EE\|\widehat\beta-\beta\|_2^2
    \lesssim
    \Gamma_K^2\eta^{*2}(K).
\end{align*}
Thus the theorem holds with, for example,
\begin{align*}
    \Upsilon(K)
    \asymp
    \Gamma_K^2,
\end{align*}
which depends polynomially, up to logarithmic factors, on $T_2(K)$ and
$\max_{j\in[M]}\kappa(K_{(j)})$.

It remains only to note the computational complexity. Since $R=O(1)$ and
\begin{align*}
    r_{\rm alg}
    =
    r_0\wedge \frac12\sqrt{\frac nN},
\end{align*}
the number of iterations is
\begin{align*}
    M
    =
    O\!\left(\log\frac{R}{r_{\rm alg}}\right)
    =
    O\!\left(\log\frac{R}{r_0}+\log N\right),
\end{align*}
up to absolute constants. The weak projections, the SDP/QFM computations,
and the approximate convex optimization in
Lemma~\ref{sub:gaussian:vector:estimator:regression} are polynomial-time under
the assumptions of Section~\ref{regression:section}. The required oracle
failure probability at scale $\widetilde d_j$ has logarithm
$O(\!\log(1+N^{3/2}\widetilde d_j^3))$, which is polynomially bounded over
the active scales. Hence the entire procedure runs in polynomial time. This
completes the proof.
\end{proof}

\begin{lemma}[High-probability bound for the OLS pilot estimator]
\label{lem:ols:pilot:bound}
Consider the linear model
\begin{align*}
    Y_i = Z_i^\top \beta + \xi_i,
    \qquad i=1,\ldots,m,
\end{align*}
where the noise variables $\xi_i$ are independent of the design, centered, and
sub-Gaussian with parameter $\gamma^2$. Let
\begin{align*}
    G := \sum_{i=1}^m Z_i Z_i^\top,
    \qquad
    \widehat\beta_{\rm OLS}
    :=
    G^{-1}\sum_{i=1}^m Z_iY_i.
\end{align*}
Suppose that the event
\begin{align*}
    \mathcal E_Z
    :=
    \left\{
        \lambda_{\min}\!\left(
            \frac1m G
        \right)
        \geq c_0
    \right\}
\end{align*}
holds for some absolute constant $c_0>0$. Then, conditional on the design,
for every $t\geq 0$,
\begin{align*}
    \PP\!\left(
        \left.
        \|\widehat\beta_{\rm OLS}-\beta\|_2
        >
        C\left(
            \gamma\sqrt{\frac{n}{m}}+t
        \right)
        \,\right|\, Z_1,\ldots,Z_m
    \right)
    \leq
    2\exp\left(
        -c\frac{m t^2}{\gamma^2}
    \right)
\end{align*}
on $\mathcal E_Z$, where $c,C>0$ are absolute constants depending only on
$c_0$.

In particular, if
\begin{align*}
    \PP(\mathcal E_Z^c)\leq \delta,
\end{align*}
then
\begin{align*}
    \PP\!\left(
        \|\widehat\beta_{\rm OLS}-\beta\|_2
        >
        C\left(
            \gamma\sqrt{\frac{n}{m}}+t
        \right)
    \right)
    \leq
    2\exp\left(
        -c\frac{m t^2}{\gamma^2}
    \right)
    +\delta.
\end{align*}
\end{lemma}

\begin{proof}
Conditional on the design,
\begin{align*}
    \widehat\beta_{\rm OLS}-\beta
    =
    G^{-1}\sum_{i=1}^m Z_i\xi_i.
\end{align*}
Fix $u\in S^{n-1}$. Then
\begin{align*}
    \left\langle
        u,\widehat\beta_{\rm OLS}-\beta
    \right\rangle
    =
    \sum_{i=1}^m
    \left\langle G^{-1}u,Z_i\right\rangle \xi_i.
\end{align*}
Since the $\xi_i$ are independent centered sub-Gaussian random variables, the
right-hand side is conditionally sub-Gaussian with variance proxy
\begin{align*}
\begin{aligned}
    \gamma^2
    \sum_{i=1}^m
    \left\langle G^{-1}u,Z_i\right\rangle^2
    &=
    \gamma^2
    u^\top G^{-1}
    \left(\sum_{i=1}^m Z_iZ_i^\top\right)
    G^{-1}u \\
    &=
    \gamma^2 u^\top G^{-1}u.
\end{aligned}
\end{align*}
On $\mathcal E_Z$, we have $G\succeq c_0mI_n$, and hence
\begin{align*}
    u^\top G^{-1}u
    \leq
    \frac{1}{c_0m}.
\end{align*}
Therefore, for every $s\geq0$,
\begin{align*}
    \PP\!\left(
        \left.
        \left|
        \left\langle
            u,\widehat\beta_{\rm OLS}-\beta
        \right\rangle
        \right|
        >s
        \,\right|\,Z_1,\ldots,Z_m
    \right)
    \leq
    2\exp\left(
        -c\frac{ms^2}{\gamma^2}
    \right).
\end{align*}

Let $\mathcal N$ be a $1/2$-net of $S^{n-1}$ with
$|\mathcal N|\leq 5^n$. By the standard net inequality,
\begin{align*}
    \|x\|_2
    \leq
    2\max_{u\in\mathcal N}|\langle u,x\rangle|.
\end{align*}
Applying the preceding tail bound to every $u\in\mathcal N$ and taking a
union bound yields
\begin{align*}
    \PP\!\left(
        \left.
        \|\widehat\beta_{\rm OLS}-\beta\|_2
        >
        C\left(
            \gamma\sqrt{\frac{n}{m}}+t
        \right)
        \,\right|\,Z_1,\ldots,Z_m
    \right)
    \leq
    2\exp\left(
        -c\frac{mt^2}{\gamma^2}
    \right)
\end{align*}
on $\mathcal E_Z$. The unconditional statement follows by adding
$\PP(\mathcal E_Z^c)$.
\end{proof}

\noindent
Under the sub-Gaussian random-design assumptions of Section~\ref{regression:section},
standard concentration for sub-Gaussian sample covariance matrices gives
$\PP(\mathcal E_Z^c)\leq 2\exp(-c m)$ provided $m\gtrsim n$. Thus, for
$m\asymp N$,
\begin{align*}
    \|\widehat\beta_{\rm OLS}-\beta\|_2
    \lesssim
    \gamma\sqrt{\frac nN}+t
\end{align*}
with probability at least
\begin{align*}
    1
    -
    2\exp\left(-c\frac{Nt^2}{\gamma^2}\right)
    -
    2\exp(-cN).
\end{align*}

\section{Supplemental Results}

\begin{lemma}\label{weak:projection:lemma}

Given a convex body with a weak separation oracle, we can weakly project on $K$, in the sense that for any $z \in \RR^n$, and $\epsilon > 0$ we can find a rational number $w$ and a point $p \in K$ in polynomial-time, such that 
\begin{align*}
    w - \epsilon \leq \|z - \overline{\Pi}_K(z)\|_2 \leq \|z - p\|_2 \leq w,
\end{align*}
where $\overline{\Pi}_K(z)$ denotes the true Euclidean projection of $z$ on $K$. Moreover, for any point $\nu \in K$ we have 
\begin{align*}
    \|p - \nu\|_2 \leq \|z- \nu\|_2 + \sqrt{\epsilon^2 + 2\epsilon\|z-\nu\|}.
\end{align*}
\end{lemma}
\begin{corollary}
    Hence if $\|z - \nu\| \leq L $ and $\epsilon \leq L$ we have
    \begin{align*}
    \|p - \nu\|_2 \leq (\sqrt{3} + 1)L
\end{align*}
\end{corollary}

\begin{proof}
    This is a direct corollary of Theorem 2.5.9 of \cite{dadush2012integer}, see also \cite{lee2018efficient} using the fact that the $\ell_2$ norm is $1$-Lipschitz. 

    For the second part we first observe that
    \begin{align*}
        \|\nu - \Pi_K(z)\|_2 \leq \|z - \nu\|_2
    \end{align*}
    Next by the properties of the Euclidean projection on convex bodies we have:
    \begin{align*}
        \|p - \Pi_K(z)\|_2^2 + \|z - \Pi_K(z)\|_2^2 \leq  \|z - p\|_2^2 \leq  \|z - \Pi_K(z)\|_2^2 + \epsilon^2 + 2 \epsilon  \|z - \Pi_K(z)\|_2
    \end{align*}
    Thus 
    \begin{align*}
        \|p - \Pi_K(z)\|_2 \leq \sqrt{\epsilon^2 + 2 \epsilon  \|z - \Pi_K(z)\|_2} \leq \sqrt{\epsilon^2 + 2 \epsilon  \|z - \nu\|_2} \leq \sqrt{3} L,
    \end{align*}
    where we used the inequlity $ \|z -\nu\|_2 \vee \epsilon \leq L$.
    By the triangle inequality 
    \begin{align*}
        \|\nu - p\| \leq  \|\nu - \Pi_K(z)\|_2  + \|p - \Pi_K(z)\|_2  \leq \|z - \nu\|_2 + \sqrt{3} L \leq  (\sqrt{3}+1)L.
    \end{align*}

\end{proof}

\begin{proof}[Proof of Lemma \ref{lemma:QCO}]
\textbf{Step 1: Quadratic representation of the norm.} Let $\|\cdot\|_K = \rho_K(\cdot)$.
Define
\begin{align*}
    K^{2} := \{\operatorname{sign}(x)\circ x^2 : x \in K\}.
\end{align*}
Since $K$ is sign-invariant and quadratically convex, $K^2$ is convex, sign-invariant, and symmetric.  
By direct inspection,
\begin{align*}
    |x| \in tK \iff x^2 \in t^2 K^2,
\end{align*}
and therefore
\begin{equation}\label{eq:qco-norm}
    \|x\|_K = \sqrt{\|x^2\|_{K^2}}.
\end{equation}

\textbf{Step 2: Reduction to Gaussian averages.}
By \cite[Lemma 4.5]{ledoux2013probability},
\begin{align*}
    \mathbb{E}\Big\|\sum_{i=1}^m \varepsilon_i x_i\Big\|_K^2
    \le \frac{\pi}{2}\,\mathbb{E}\Big\|\sum_{i=1}^m g_i x_i\Big\|_K^2,
\end{align*}
where $g_i \sim N(0,1)$ are independent Gaussians.

\textbf{Step 3: Coordinatewise representation.}
Let $S := \sum_{i=1}^m g_i x_i \in \mathbb{R}^n$.  
Then for each coordinate $j$,
\begin{align*}
    S_j = \sum_{i=1}^m g_i x_{ij} \sim N(0,\sigma_j^2),
    \quad \sigma_j^2 = \sum_i x_{ij}^2.
\end{align*}
Hence $S_j^2 = \sigma_j^2 \chi_j^2$, where $\chi_j^2 \sim \chi^2_1$, and in vector form
\begin{align*}
    S^2 = \chi^2 \circ \Big(\sum_i x_i^2\Big),
\end{align*}
where ``$\circ$'' denotes coordinatewise multiplication.

\textbf{Step 4: Monotonicity of sign-invariant norms.}
For any unconditional norm,
\begin{align*}
    |u| \le |v| \implies \|u\|_{K^2} \le \|v\|_{K^2}.
\end{align*}
Therefore,
\begin{align*}
    \|S^2\|_{K^2}
    \le \big(\max_j \chi_j^2\big)\, \Big\|\sum_i x_i^2\Big\|_{K^2}.
\end{align*}

\textbf{Step 5: Expectation bound.}
Taking expectations and using subexponential tails of $\chi^2_1$:
\begin{align*}
    \mathbb{E}\|S^2\|_{K^2}
    \le \mathbb{E}\max_j \chi_j^2 \cdot \Big\|\sum_i x_i^2\Big\|_{K^2},
\end{align*}
and
\begin{align*}
    \mathbb{E}\max_{1\le j\le n}\chi_j^2 \lesssim \log n.
\end{align*}
Hence
\begin{equation}\label{eq:exp-bound}
    \mathbb{E}\|S^2\|_{K^2}
    \lesssim (\log n)\Big\|\sum_i x_i^2\Big\|_{K^2}.
\end{equation}

\textbf{Step 6: Subadditivity and conclusion.}
By convexity of the norm,
\begin{align*}
    \Big\|\sum_i x_i^2\Big\|_{K^2} \le \sum_i \|x_i^2\|_{K^2}.
\end{align*}
Using \eqref{eq:qco-norm} and \eqref{eq:exp-bound},
\begin{align*}
    \mathbb{E}\|S\|_K^2
    = \mathbb{E}\|S^2\|_{K^2}
    \lesssim (\log n)\sum_i \|x_i\|_K^2.
\end{align*}
Combining with step 2 yields the desired type-2 inequality.
\end{proof}
\begin{proof}[Proof of Proposition \ref{not:all:type2:sets:are:QCO}]

For $u = (1,0,\dots,0)$ and $v = (0,1,0,\dots,0)$, we have
\begin{align*}
\|u\|_2 = \|v\|_2 = 1, \qquad \|u\|_1 = \|v\|_1 = 1,
\end{align*}
hence
\begin{align*}
\|u\| = \|v\| = \sqrt{1 + \frac{1}{n}} =: c_n > 1.
\end{align*}
Let $\tilde u = u/c_n$ and $\tilde v = v/c_n$.
Then $\|\tilde u\| = \|\tilde v\| = 1$, so $\tilde u, \tilde v \in K$.

Consider the midpoint in squared coordinates for $t = \tfrac{1}{2}$:
\begin{align*}
\frac{1}{2}(\tilde u^2 + \tilde v^2)
   = \frac{1}{2c_n^2}(1,1,0,\dots,0),
\end{align*}
and take the entrywise square root
\begin{align*}
w = \sqrt{\tfrac{1}{2}(\tilde u^2 + \tilde v^2)}
  = \frac{1}{c_n\sqrt{2}}(1,1,0,\dots,0).
\end{align*}
We now compute $\|w\|$:
\begin{align*}
\|w\|_2^2 = \frac{1}{c_n^2}, 
\qquad
\|w\|_1^2 = \Big(\frac{2}{c_n\sqrt{2}}\Big)^2 = \frac{2}{c_n^2}.
\end{align*}
Therefore
\begin{align*}
\|w\|^2
   = \|w\|_2^2 + \frac{1}{n}\|w\|_1^2
   = \frac{1}{c_n^2}\Big(1 + \frac{2}{n}\Big)
   = \frac{1 + \tfrac{2}{n}}{1 + \tfrac{1}{n}}
   = \frac{n+2}{n+1}
   = 1 + \frac{1}{n+1}.
\end{align*}
Hence $\|w\| = \sqrt{1 + \tfrac{1}{n+1}} > 1$, so $w \notin K$.

Thus, although $\tilde u, \tilde v \in K$, the vector
$\sqrt{(\tilde u^2+\tilde v^2)/2}$ leaves the unit ball. The latter shows that $K$ is not quadratically convex.

Finally, since
\begin{align*}
\|x\|_2 \le \|x\| \le \sqrt{2}\|x\|_2,
\end{align*}
the norm is $O(1)$--equivalent to the Euclidean norm and therefore
has type--2 constant bounded by an absolute constant.
\end{proof}

\begin{lemma}\label{K2:is:convex}
Let $K\subseteq\mathbb R^n$ be convex, sign-invariant, and quadratically convex. Define
\begin{align*}
K^2:=\{\operatorname{sign}(x)\circ x^2:x\in K\}.
\end{align*}
Then $K^2$ is convex and sign-invariant.
\end{lemma}

\begin{proof}
Sign-invariance is immediate. Let
\begin{align*}
C:=\{x^2:x\in K\}\subseteq\mathbb R_+^n.
\end{align*}
By quadratic convexity, $C$ is convex. Moreover, since $K$ is convex and sign-invariant, it is coordinatewise solid: if $x\in K$ and $|y|\le |x|$ coordinatewise, then $y\in K$. Hence $C$ is downward closed in $\mathbb R_+^n$.

Now let $a,b\in K^2$ and $t\in[0,1]$. Since $|a|,|b|\in C$,
\begin{align*}
|ta+(1-t)b|
\le t|a|+(1-t)|b|\in C.
\end{align*}
By downward closedness, $|ta+(1-t)b|\in C$. Sign-invariance then implies
\begin{align*}
ta+(1-t)b\in K^2.
\end{align*}
Thus $K^2$ is convex.
\end{proof}

\section{A Counterexample for the Posterior Mean under the Uniform Prior}\label{posterior:mean}

Consider the Gaussian sequence model
\begin{align*}
    Y=\mu+Z,\qquad Z\sim N(0,I_n),
\end{align*}
with parameter constrained to the ellipsoid
\begin{align*}
    K_n
    :=
    \left\{
        \mu=(u,t)\in \RR^{n-1}\times \RR:
        \|u\|_2^2+\frac{t^2}{\sqrt n}\le 1
    \right\}.
\end{align*}
Write
\begin{align*}
    L_n:=n^{1/4},
\end{align*}
so that equivalently
\begin{align*}
    K_n
    =
    \left\{
        (u,t):
        \|u\|_2^2+\frac{t^2}{L_n^2}\le 1
    \right\}.
\end{align*}
Let $\Pi_n=\operatorname{Unif}(K_n)$, and let
\begin{align*}
    \widehat\mu^{\,\mathrm{unif}}(Y)
    :=
    \EE{\Pi_n}[\mu\mid Y]
\end{align*}
denote the corresponding posterior mean. We note that the set $K_n$ is a type-2 body with constant $T_2(K_n) = 1$, which implies our algorithm works optimally on $K_n$.

\begin{proposition}[The uniform-prior posterior mean can be polynomially suboptimal]
There exist universal constants $c>0$ and $n_0<\infty$ such that, for every
$n\ge n_0$,
\begin{align*}
    \sup_{\mu\in K_n}
    \EE\mu
    \bigl\|
        \widehat\mu^{\,\mathrm{unif}}(Y)-\mu
    \bigr\|_2^2
    \ge c\sqrt n.
\end{align*}
On the other hand,
\begin{align*}
    \inf_{\widehat\mu}
    \sup_{\mu\in K_n}
    \EE\mu
    \|\widehat\mu(Y)-\mu\|_2^2
    \le 2.
\end{align*}
Consequently the posterior mean under the uniform prior is not minimax, and
its worst-case risk can exceed the minimax risk by a factor of order
$\sqrt n$.

More concretely, one may take $c=1/32$ for all sufficiently large $n$.
\end{proposition}

\begin{proof}
We prove the lower bound by evaluating the risk at the tip of the long axis,
\begin{align*}
    \mu^\star:=L_n e_n.
\end{align*}

\medskip
\noindent
\textbf{Step 1: Posterior density.}
Since the prior is uniform on $K_n$, the posterior density with respect to
Lebesgue measure on $K_n$ is proportional to
\begin{align*}
    \exp\left(
        -\frac12\|Y-\mu\|_2^2
    \right)
    =
    \exp\left(
        \langle Y,\mu\rangle-\frac12\|\mu\|_2^2
    \right).
\end{align*}
Thus, for fixed $y$, define
\begin{align*}
    q_y(\mu)
    :=
    \langle y,\mu\rangle-\frac12\|\mu\|_2^2.
\end{align*}

Write $\mu=(u,t)$ and $y=(y',y_n)$, where
$u,y'\in\RR^{n-1}$.

\medskip
\noindent
\textbf{Step 2: The uniform volume of the tip is exponentially small.}
For fixed $t\in[-L_n,L_n]$, the cross-section of $K_n$ at height $t$ is the
$(n-1)$-dimensional Euclidean ball
\begin{align*}
    \left\{
        u:
        \|u\|_2
        \le
        \sqrt{1-\frac{t^2}{L_n^2}}
    \right\}.
\end{align*}
Let $v_{n-1}$ denote the volume of the unit ball in $\RR^{n-1}$.

Define the tip region
\begin{align*}
    T_n
    :=
    \left\{
        (u,t)\in K_n:
        |t|\ge \frac{L_n}{2}
    \right\}
\end{align*}
and the central slab
\begin{align*}
    C_n
    :=
    \left\{
        (u,t)\in K_n:
        |t|\le \frac{L_n}{4}
    \right\}.
\end{align*}
For $|t|\ge L_n/2$,
\begin{align*}
    1-\frac{t^2}{L_n^2}
    \le \frac34,
\end{align*}
and therefore
\begin{align*}
    \operatorname{Vol}(T_n)
    \le
    2L_n\,v_{n-1}
    \left(\frac34\right)^{(n-1)/2}.
\end{align*}
Similarly, for $|t|\le L_n/4$,
\begin{align*}
    1-\frac{t^2}{L_n^2}
    \ge \frac{15}{16},
\end{align*}
so
\begin{align*}
    \operatorname{Vol}(C_n)
    \ge
    \frac{L_n}{2}\,v_{n-1}
    \left(\frac{15}{16}\right)^{(n-1)/2}.
\end{align*}
Hence
\begin{equation}\label{eq:volume-ratio}
    \frac{\operatorname{Vol}(T_n)}
         {\operatorname{Vol}(C_n)}
    \le
    4
    \left(\frac45\right)^{(n-1)/2}.
\end{equation}

\medskip
\noindent
\textbf{Step 3: The likelihood cannot compensate for the volume loss.}
Consider the event
\begin{align*}
    \mathcal E_n
    :=
    \left\{
        \|Y'\|_2\le 2\sqrt n,\quad |Y_n|\le 2L_n
    \right\}.
\end{align*}
Suppose that $y\in\mathcal E_n$.

For every $\mu=(u,t)\in K_n$, we have $\|u\|_2\le 1$ and $|t|\le L_n$.
Therefore
\begin{align*}
    q_y(\mu)
    \le
    |\langle y',u\rangle|+|y_n t|
    \le
    2\sqrt n+2L_n^2
    =
    4\sqrt n,
\end{align*}
because $L_n^2=\sqrt n$.

On the other hand, if $\mu=(u,t)\in C_n$, then
$|t|\le L_n/4$, and hence
\begin{align*}
    q_y(\mu)
    &\ge
    -\|y'\|_2\|u\|_2
    -|y_n||t|
    -\frac12\bigl(\|u\|_2^2+t^2\bigr)
    \\
    &\ge
    -2\sqrt n
    -\frac12 L_n^2
    -\frac12\left(1+\frac{L_n^2}{16}\right)
    \\
    &=
    -\frac{81}{32}\sqrt n-\frac12.
\end{align*}

Let $\Pi_n(\,\cdot\mid y)$ denote the posterior distribution. Using the
previous two bounds and \eqref{eq:volume-ratio},
\begin{align}
    \Pi_n(T_n\mid y)
    &=
    \frac{
        \int_{T_n} e^{q_y(\mu)}\,d\mu
    }{
        \int_{K_n} e^{q_y(\mu)}\,d\mu
    }
    \nonumber\\
    &\le
    \frac{
        e^{4\sqrt n}\operatorname{Vol}(T_n)
    }{
        e^{-(81/32)\sqrt n-1/2}\operatorname{Vol}(C_n)
    }
    \nonumber\\
    &\le
    4
    \left(\frac45\right)^{(n-1)/2}
    \exp\left(
        \frac{209}{32}\sqrt n+\frac12
    \right)
    =:\varepsilon_n.
    \label{eq:posterior-tip}
\end{align}
Since
\begin{align*}
    \log \varepsilon_n
    =
    -\frac{n-1}{2}\log\frac54
    +O(\sqrt n),
\end{align*}
we have $\varepsilon_n\to0$. In particular, for all sufficiently large $n$,
\begin{align*}
    \varepsilon_n\le \frac14.
\end{align*}

\medskip
\noindent
\textbf{Step 4: The posterior mean stays far from the tip.}
Let
\begin{align*}
    \widehat t(y)
    :=
    \EE{\Pi_n}[t\mid Y=y]
\end{align*}
be the last coordinate of the posterior mean. Since $|t|\le L_n$ everywhere
on $K_n$, while $|t|<L_n/2$ outside $T_n$, we obtain, on $\mathcal E_n$,
\begin{align*}
    |\widehat t(y)|
    &\le
    \EE{\Pi_n}[|t|\mid y]
    \\
    &\le
    \frac{L_n}{2}
    +
    L_n\,\Pi_n(T_n\mid y)
    \\
    &\le
    \frac{3L_n}{4}
\end{align*}
for all sufficiently large $n$. Therefore, at the true parameter
$\mu^\star=L_ne_n$,
\begin{align*}
    \left|
        \widehat t(Y)-L_n
    \right|
    \ge
    L_n-|\widehat t(Y)|
    \ge
    \frac{L_n}{4}
\end{align*}
on $\mathcal E_n$. Consequently,
\begin{equation}\label{eq:risk-event}
    \EE{\mu^\star}
    \left[
        \bigl(\widehat t(Y)-L_n\bigr)^2
    \right]
    \ge
    \frac{L_n^2}{16}
    \PP_{\mu^\star}(\mathcal E_n).
\end{equation}

Under $\mu^\star$,
\begin{align*}
    Y'=Z',\qquad
    Y_n=L_n+Z_n,
\end{align*}
where $Z'\sim N(0,I_{n-1})$ and $Z_n\sim N(0,1)$. Hence
\begin{align*}
    \PP\bigl(\|Z'\|_2\le 2\sqrt n\bigr)\to1
\end{align*}
and, since $L_n\to\infty$,
\begin{align*}
    \PP\bigl(|L_n+Z_n|\le 2L_n\bigr)
    =
    \PP(-3L_n\le Z_n\le L_n)
    \to1.
\end{align*}
Thus
\begin{align*}
    \PP_{\mu^\star}(\mathcal E_n)\to1.
\end{align*}
In particular, for all sufficiently large $n$,
\begin{align*}
    \PP_{\mu^\star}(\mathcal E_n)\ge\frac12.
\end{align*}
Combining this with \eqref{eq:risk-event} and $L_n^2=\sqrt n$ gives
\begin{align*}
    \EE{\mu^\star}
    \left[
        \bigl(\widehat t(Y)-L_n\bigr)^2
    \right]
    \ge
    \frac{L_n^2}{32}
    =
    \frac{\sqrt n}{32}.
\end{align*}
Since the full squared error is at least the squared error in the last
coordinate,
\begin{align*}
    \sup_{\mu\in K_n}
    \EE\mu
    \bigl\|
        \widehat\mu^{\,\mathrm{unif}}(Y)-\mu
    \bigr\|_2^2
    \ge
    \frac{\sqrt n}{32}.
\end{align*}

\medskip
\noindent
\textbf{Step 5: A uniformly bounded-risk estimator exists.}
Consider
\begin{align*}
    \widetilde\mu(Y)
    :=
    (0,\ldots,0,Y_n).
\end{align*}
For any $\mu=(u,t)\in K_n$,
\begin{align*}
    \EE\mu
    \|\widetilde\mu(Y)-\mu\|_2^2
    =
    \|u\|_2^2+\EE Z_n^2
    =
    \|u\|_2^2+1
    \le 2,
\end{align*}
because membership in $K_n$ implies $\|u\|_2^2\le1$. Hence
\begin{align*}
    \inf_{\widehat\mu}
    \sup_{\mu\in K_n}
    \EE\mu
    \|\widehat\mu(Y)-\mu\|_2^2
    \le2.
\end{align*}
This proves the proposition.
\end{proof}

\begin{remark}[Geometric reason for the failure]
The last semiaxis of $K_n$ has length $L_n=n^{1/4}$, but under the uniform
measure on $K_n$ the marginal density of the last coordinate $t$ is
proportional to
\begin{align*}
    \left(
        1-\frac{t^2}{L_n^2}
    \right)^{(n-1)/2},
    \qquad |t|\le L_n.
\end{align*}
Thus the uniform prior places exponentially little mass near the two tips
$\pm L_ne_n$.  At the same time, these tips are statistically easy to
estimate by retaining the last observation $Y_n$.  The posterior associated
with the uniform prior therefore overweights the enormous central volume of
the ellipsoid and shrinks strongly away from the extremal points.  This
volume effect is precisely what makes the uniform-prior posterior mean
non-minimax.
\end{remark}

\section{Construction and conditional computational lower bound}

We describe a centrally symmetric convex body whose gauge is computable by a
polynomial-size semidefinite program, but for which minimax-rate estimation in the
Gaussian sequence model would imply recovery in the planted-submatrix model in the
computationally hard regime identified by Schramm and Wein
\cite{schrammwein2022}.  The resulting computational lower bound is conditional on
the usual low-degree heuristic/conjecture: Schramm and Wein prove a rigorous
low-degree lower bound, rather than an unconditional lower bound against all
polynomial-time algorithms.

Fix a sparsity parameter $\rho=\rho_n\in(0,1/2)$ and write
\begin{align*}
    k=\rho n.
\end{align*}
For simplicity we ignore integer-rounding issues.  Let $\delta_n\downarrow 0$ satisfy
\begin{align*}
    \delta_n^2 k\longrightarrow\infty,
\end{align*}
for example $\delta_n=k^{-1/4}$, and define
\begin{align*}
    k_-=(1-\delta_n)k,
    \qquad
    k_+=(1+\delta_n)k.
\end{align*}

Define the convex cone
\begin{align*}
\cC_{n,k}
=
\left\{
X\in\Smat^n:
\begin{array}{l}
X\succeq 0,\qquad X_{ij}\ge 0,\\[1mm]
X_{ij}\le X_{ii},\qquad X_{ij}\le X_{jj}
\quad\text{for all }i,j,\\[1mm]
k_-X_{ii}\le \displaystyle\sum_{j=1}^nX_{ij}
\le k_+X_{ii}
\quad\text{for every }i
\end{array}
\right\}.
\end{align*}
Let
\begin{align*}
    K_{n,k}^+
    =
    \left\{
        X\in\cC_{n,k}:\tr(X)\le 1
    \right\},
\end{align*}
and take the centrally symmetric convex body
\begin{equation}
    K_{n,k}
    =
    \conv\bigl(K_{n,k}^+\cup(-K_{n,k}^+)\bigr).
    \label{eq:def-K}
\end{equation}

\subsection{The gauge is efficiently computable}

\begin{lemma}[SDP representation of the gauge]
For every $X\in\Smat^n$,
\begin{equation}
    \|X\|_{K_{n,k}}
    =
    \inf_{\substack{X=X_+-X_-\\ X_+,X_-\in\cC_{n,k}}}
    \left\{
        \tr(X_+)+\tr(X_-)
    \right\}.
    \label{eq:gauge}
\end{equation}
In particular, the gauge of $K_{n,k}$ is computable by a polynomial-size
semidefinite program.
\end{lemma}

\begin{proof}
Since $\cC_{n,k}$ is a convex cone, the absolute convexification in
\eqref{eq:def-K} can equivalently be written as
\begin{align*}
K_{n,k}
=
\left\{
X_+-X_-:
X_\pm\in\cC_{n,k},\;
\tr(X_+)+\tr(X_-)\le 1
\right\}.
\end{align*}
Indeed, if
$X=tX_+-(1-t)X_-$ with $X_\pm\in K_{n,k}^+$ and $t\in[0,1]$,
then $tX_+,(1-t)X_-\in\cC_{n,k}$ and the sum of their traces is at
most one.  The converse follows by normalizing the two summands by
their total trace, i.e. take $X \in \left\{
X_+-X_-:
X_\pm\in\cC_{n,k},\;
\tr(X_+)+\tr(X_-)\le 1
\right\}$ and observe that 
\begin{align*}
X = s\bigg(\frac{\tr(X_+)}{s}\frac{X_{+}}{\tr(X_+)} - \frac{\tr(X_-)}{s}\frac{X_{-}}{\tr(X_-)}\bigg).    
\end{align*}
where $s := \tr(X_+) + \tr(X_-) \leq 1$. Since $s \leq 1$ it follows $X \in K_{n,k}$.
Taking the Minkowski functional gives
\eqref{eq:gauge}.  All constraints defining $\cC_{n,k}$ are linear
matrix inequalities or linear inequalities, so \eqref{eq:gauge} is
an SDP of polynomial size.
\end{proof}

\subsection{Embedding the Schramm--Wein planted submatrix}

Schramm and Wein consider the planted-submatrix model
\begin{equation}
    Y=\lambda vv^\top+W,
    \qquad
    v_i\stackrel{\mathrm{iid}}{\sim}\mathrm{Bernoulli}(\rho),
    \label{eq:SW-model}
\end{equation}
with Gaussian noise $W$; see \cite{schrammwein2022}.  Their theorem is
stated in the parameterization
\begin{align*}
    \rho=n^{-b},
    \qquad
    \lambda=n^{-a},
\end{align*}
with $a>0$ and $0<b<1/2$.

Let
\begin{align*}
    m=\|v\|_0.
\end{align*}
By the Chernoff bound and the assumption $\delta_n^2k\to\infty$,
\begin{equation}
    \Pp\{k_-\le m\le k_+\}=1-o(1).
    \label{eq:support-concentration}
\end{equation}
On the event in \eqref{eq:support-concentration}, set
\begin{align*}
    A_v=\frac{vv^\top}{k_+}.
\end{align*}
Then
\begin{align*}
    \tr(A_v)=\frac{m}{k_+}\le 1.
\end{align*}
Moreover, if $v_i=1$, then
\begin{align*}
    \frac{\sum_j(A_v)_{ij}}{(A_v)_{ii}}
    =
    m\in[k_-,k_+],
\end{align*}
while if $v_i=0$, the entire $i$th row and column are zero.  Also
$A_v\succeq0$, $A_v\ge0$ entrywise, and
\begin{align*}
    (A_v)_{ij}\le (A_v)_{ii},
    \qquad
    (A_v)_{ij}\le (A_v)_{jj}.
\end{align*}
Thus
\begin{equation}
    A_v\in K_{n,k}^+\subset K_{n,k}
    \qquad\text{with probability }1-o(1).
    \label{eq:Av-in-K}
\end{equation}
Consequently, if
\begin{equation}
    \theta=\lambda k_+,
    \label{eq:theta-def}
\end{equation}
then
\begin{equation}
    \lambda vv^\top
    =
    \theta A_v
    \in
    \theta K_{n,k}
    \qquad\text{with probability }1-o(1).
    \label{eq:plant-in-body}
\end{equation}

The Gaussian matrix observation can be identified with an ordinary Gaussian sequence
model on the upper-triangular coordinates.  Any difference between the standard
normalization used for a symmetric Gaussian matrix and iid $N(0,1)$ sequence noise
amounts only to a fixed linear rescaling of the diagonal coordinates and therefore
changes all Euclidean/Frobenius quantities by universal constants.

\subsection{Metric entropy of the convex body}

The particular linear inequalities in the definition of $\cC_{n,k}$ ensure that the
relaxation remains statistically small.

\begin{lemma}[Large-scale entropy bound]
There exists a universal constant $C$ such that, for
$0<\varepsilon<1/2$,
\begin{equation}
\log N\!\left(K_{n,k},\varepsilon B_F\right)
\le
C\left[
    \frac{k}{\varepsilon^2}\log\frac{en}{k}
    +
    \frac{k}{\varepsilon^4}\log\frac{C}{\varepsilon}
\right],
\label{eq:entropy}
\end{equation}
where $B_F$ is the unit ball for the Frobenius norm.
\end{lemma}

\begin{proof}
We first consider $X\in K_{n,k}^+$.  Put $d_i=X_{ii}$.  Since
$0\le X_{ij}\le X_{ii}=d_i$ and
$\sum_jX_{ij}\le k_+d_i$,
\begin{equation}
    \sum_{j=1}^nX_{ij}^2
    \le
    d_i\sum_{j=1}^nX_{ij}
    \le
    k_+d_i^2.
    \label{eq:row-bound}
\end{equation}
Order the diagonal entries so that
$d_1\ge d_2\ge\cdots\ge0$.  Because
$\sum_i d_i=\tr(X)\le1$, we have
$d_{s+1}\le 1/s$ and hence
\begin{align*}
    \sum_{i>s}d_i^2
    \le
    d_{s+1}\sum_{i>s}d_i
    \le
    \frac1s.
\end{align*}
If $T$ is the set of the $s$ largest diagonal coordinates, then
\begin{align}
    \|X-X_{TT}\|_F^2
    &\le
    2\sum_{i\notin T}\sum_{j=1}^nX_{ij}^2 \notag\\
    &\le
    2k_+\sum_{i\notin T}d_i^2
    \le
    \frac{2k_+}{s}.
    \label{eq:support-compress}
\end{align}
Thus, choosing
\begin{align*}
    s\asymp k\varepsilon^{-2}
\end{align*}
makes $X$ $O(\varepsilon)$-close to a matrix supported on an
$s\times s$ principal block.

Next, $X_{TT}\succeq0$ and $\tr(X_{TT})\le1$.  If
$\lambda_1\ge\lambda_2\ge\cdots\ge0$ are its eigenvalues, then
$\sum_j\lambda_j\le1$, and therefore
\begin{align*}
    \sum_{j>r}\lambda_j^2
    \le
    \lambda_{r+1}\sum_{j>r}\lambda_j
    \le
    \frac1r.
\end{align*}
Hence the best rank-$r$ approximation has Frobenius error at most
$r^{-1/2}$.  Taking $r\asymp\varepsilon^{-2}$ reduces the covering
problem to rank-$O(\varepsilon^{-2})$ positive semidefinite matrices
on $O(k\varepsilon^{-2})$ coordinates.

There are at most
\begin{align*}
    \binom{n}{s}
    \le
    \left(\frac{en}{s}\right)^s
\end{align*}
choices of $T$.  For each fixed $T$, a standard covering bound for
rank-$r$ matrices \cite[see Lemma 3.1][together with a trivial discretization of the Frobenius norm]{candes2011tight} in the Frobenius unit ball gives logarithmic
covering number at most
\begin{align*}
    Csr\log\frac{C}{\varepsilon}.
\end{align*}
Substituting $s\asymp k\varepsilon^{-2}$ and
$r\asymp\varepsilon^{-2}$ yields
\begin{align*}
    \log N(K_{n,k}^+,\varepsilon B_F)
    \lesssim
    \frac{k}{\varepsilon^2}\log\frac{en}{k}
    +
    \frac{k}{\varepsilon^4}\log\frac{C}{\varepsilon},
\end{align*}
after harmless changes of constants.

Finally, every $X\in K_{n,k}$ can be written as
$X=X_+-X_-$ with $X_\pm\in K_{n,k}^+$.  Covering the two summands
separately only squares the covering number and changes constants.
This proves \eqref{eq:entropy}.
\end{proof}

\subsection{The statistical minimax rate is much smaller than the planted signal}

Consider the Gaussian sequence model
\begin{equation}
    Z=\mu+\xi,
    \qquad
    \xi\sim N(0,\sigma^2I),
    \qquad
    \mu\in\theta K_{n,k}.
    \label{eq:gsm}
\end{equation}
Let
\begin{align*}
    R^*(\theta K_{n,k})
    =
    \inf_{\widehat\mu}
    \sup_{\mu\in\theta K_{n,k}}
    \E_\mu\|\widehat\mu-\mu\|_F^2
\end{align*}
denote the minimax squared-error risk, up to the harmless fixed linear
identification between symmetric matrices and Gaussian sequence coordinates.

A minimum-distance estimator on an $\varepsilon\theta$-net gives
\begin{equation}
    R^*(\theta K_{n,k})
    \lesssim
    \theta^2\varepsilon^2
    +
    \sigma^2
    \log N(K_{n,k},\varepsilon B_F).
    \label{eq:net-risk}
\end{equation}
This is a folklore result that can be traced back to \cite{yang1999information}; for completeness, we provide a short proof.
\begin{lemma}[Risk bound for a minimum-distance estimator on a net]
Let $T\subseteq \mathbb{R}^d$ and suppose that
\begin{align*}
    Y=\mu+\sigma Z,\qquad Z\sim N(0,I_d),
\end{align*}
where $\mu\in T$. Let $\mathcal N_\delta\subseteq T$ be a finite
$\delta$-net of $T$ in Euclidean norm, and define
\begin{align*}
    \widehat\mu\in\arg\min_{\nu\in\mathcal N_\delta}\|Y-\nu\|_2.
\end{align*}
Then there exists a universal constant $C>0$ such that
\begin{align*}
    \sup_{\mu\in T}\mathbb E_\mu\|\widehat\mu-\mu\|_2^2
    \le
    C\left(
        \delta^2+\sigma^2\log |\mathcal N_\delta|
    \right).
\end{align*}
Consequently,
\begin{align*}
    R^*(T)
    \lesssim
    \delta^2+\sigma^2\log N(T,\delta B_2).
\end{align*}
In particular, for $T=\theta K_{n,k}$ and
$\delta=\theta\varepsilon$,
\begin{align*}
    R^*(\theta K_{n,k})
    \lesssim
    \theta^2\varepsilon^2
    +
    \sigma^2\log N(K_{n,k},\varepsilon B_F).
\end{align*}
\end{lemma}

\begin{proof}
Fix $\mu\in T$, and choose $\nu_0\in\mathcal N_\delta$ such that
\begin{align*}
    \|\mu-\nu_0\|_2\le\delta.
\end{align*}
By the definition of $\widehat\mu$,
\begin{align*}
    \|Y-\widehat\mu\|_2^2
    \le
    \|Y-\nu_0\|_2^2.
\end{align*}
Substituting $Y=\mu+\sigma Z$ and expanding gives
\begin{align*}
    \|\widehat\mu-\nu_0\|_2^2
    \le
    2\langle \mu-\nu_0,\widehat\mu-\nu_0\rangle
    +
    2\sigma\langle Z,\widehat\mu-\nu_0\rangle.
\end{align*}
Define
\begin{align*}
    M
    :=
    \max_{\substack{\nu\in\mathcal N_\delta\\ \nu\neq\nu_0}}
    \frac{\langle Z,\nu-\nu_0\rangle}
         {\|\nu-\nu_0\|_2},
\end{align*}
with the convention $M=0$ if $\mathcal N_\delta=\{\nu_0\}$.
Since $\widehat\mu\in\mathcal N_\delta$, Cauchy--Schwarz and the
definition of $M$ imply
\begin{align*}
    \|\widehat\mu-\nu_0\|_2^2
    \le
    2\delta\|\widehat\mu-\nu_0\|_2
    +
    2\sigma M\|\widehat\mu-\nu_0\|_2.
\end{align*}
Therefore, either $\widehat\mu=\nu_0$, or, after dividing by
$\|\widehat\mu-\nu_0\|_2$,
\begin{align*}
    \|\widehat\mu-\nu_0\|_2
    \le
    2\delta+2\sigma M.
\end{align*}
Hence
\begin{align*}
    \|\widehat\mu-\nu_0\|_2^2
    \lesssim
    \delta^2+\sigma^2 M_+^2,
\end{align*}
where $M_+=\max\{M,0\}$.

For every fixed $\nu\neq\nu_0$,
\begin{align*}
    \frac{\langle Z,\nu-\nu_0\rangle}
         {\|\nu-\nu_0\|_2}
    \sim N(0,1).
\end{align*}
Thus, by the standard Gaussian maximal inequality,
\begin{align*}
    \mathbb E M_+^2
    \lesssim
    \log |\mathcal N_\delta|
\end{align*}
whenever $|\mathcal N_\delta|\ge 2$; the case
$|\mathcal N_\delta|=1$ is immediate. Consequently,
\begin{align*}
    \mathbb E_\mu
    \|\widehat\mu-\nu_0\|_2^2
    \lesssim
    \delta^2
    +
    \sigma^2\log |\mathcal N_\delta|.
\end{align*}
Finally,
\begin{align*}
    \|\widehat\mu-\mu\|_2^2
    \le
    2\|\widehat\mu-\nu_0\|_2^2
    +
    2\|\nu_0-\mu\|_2^2,
\end{align*}
and therefore
\begin{align*}
    \sup_{\mu\in T}
    \mathbb E_\mu\|\widehat\mu-\mu\|_2^2
    \lesssim
    \delta^2
    +
    \sigma^2\log |\mathcal N_\delta|.
\end{align*}
Taking $\mathcal N_\delta$ to be a minimal $\delta$-net proves the
entropy bound. Applying it with $T=\theta K_{n,k}$ and
$\delta=\theta\varepsilon$, and using
\begin{align*}
    N(\theta K_{n,k},\theta\varepsilon B_F)
    =
    N(K_{n,k},\varepsilon B_F),
\end{align*}
gives the claimed specialization.
\end{proof}
Combining \eqref{eq:net-risk} with \eqref{eq:entropy}, and choosing
$\varepsilon=\varepsilon_n\downarrow0$ sufficiently slowly, yields
\begin{equation}
    R^*(\theta K_{n,k})
    =
    o(\theta^2)
    \qquad\text{whenever}\qquad
    \frac{\theta^2}
    {\sigma^2 k\log(en/k)}
    \longrightarrow\infty.
    \label{eq:small-relative-risk}
\end{equation}
For the planted signal in \eqref{eq:plant-in-body},
$\theta\asymp\lambda k$, so the condition becomes
\begin{equation}
    \lambda^2 k
    \gg
    \sigma^2\log\frac{en}{k}.
    \label{eq:stat-condition}
\end{equation}

\subsection{The Schramm--Wein computational regime}\label{computational:appendix}

Schramm and Wein \cite[Theorem~1.2]{schrammwein2022} consider the planted
submatrix model in the parametrization
\begin{align*}
    \rho=n^{-b},
    \qquad
    \lambda=n^{-a},
    \qquad
    0<b<\frac12.
\end{align*}
Their lower bound is already an \emph{estimation} lower bound.  More precisely,
letting $\mathsf{MMSE}_{\leq D}$ denote the minimum mean-squared error among
degree-$D$ polynomial estimators of a planted coordinate $v_1$, they prove that
if
\begin{equation}
    a>\frac12-b,
    \label{eq:SW-hard}
\end{equation}
then, for some $\epsilon=\epsilon(a,b)>0$,
\begin{equation}
    \mathsf{MMSE}_{\leq n^\epsilon}
    =
    \rho-(1+o(1))\rho^2.
    \label{eq:SW-mmse}
\end{equation}
Thus no polynomial of degree $n^\epsilon$ can achieve mean-squared error
$o(\rho)$ for estimating a planted coordinate; in particular, it cannot
estimate the full planted vector $v$ with total squared error
$o(\rho n)=o(k)$.  Conversely, if $a<1/2-b$, Schramm and Wein show that a
constant-degree polynomial achieves mean-squared error $o(\rho)$.  Hence
$a=1/2-b$ is the sharp low-degree \emph{estimation} threshold in their model.

On the other hand, our statistical condition \eqref{eq:stat-condition} with
$\sigma\asymp1$ reads
\begin{align*}
    n^{-2a}n^{1-b}\gg\log n,
\end{align*}
which holds whenever
\begin{equation}
    a<\frac{1-b}{2}
    =
    \frac12-\frac b2.
    \label{eq:stat-range}
\end{equation}
Therefore the interval
\begin{equation}
    \frac12-b
    <
    a
    <
    \frac12-\frac b2
    \label{eq:gap-range}
\end{equation}
is nonempty for every $0<b<1/2$.  Throughout this interval, the full convex
body $\theta K_{n,k}$ is statistically estimable to vanishing relative error,
while the embedded planted-submatrix instances lie in the Schramm--Wein
low-degree hard regime for estimation.

For example, take
\begin{align*}
    b=\frac13,
    \qquad
    a=\frac14.
\end{align*}
Then
\begin{align*}
    \rho=n^{-1/3},
    \qquad
    k=n^{2/3},
    \qquad
    \lambda=n^{-1/4}.
\end{align*}
We have
\begin{align*}
    \lambda^2k=n^{1/6}\gg\log n,
\end{align*}
so \eqref{eq:small-relative-risk} holds, while
\begin{align*}
    \frac14>\frac16=\frac12-b,
\end{align*}
placing the planted instance strictly in the Schramm--Wein low-degree hard
regime.

\subsection{Reduction from minimax-rate estimation to planted-submatrix estimation}

\begin{proposition}[Conditional computational-statistical separation]
Fix constants $a,b$ satisfying \eqref{eq:gap-range}, and let
\begin{align*}
    \rho=n^{-b},
    \qquad
    k=\rho n,
    \qquad
    \lambda=n^{-a},
    \qquad
    \theta=\lambda k_+.
\end{align*}
Let $L_n$ be any polylogarithmic sequence, i.e.,
\begin{align*}
    L_n\leq (\log n)^A
\end{align*}
for some fixed $A<\infty$. Suppose there exists a polynomial-time estimator
$\widehat\mu$ for the Gaussian sequence model over $\theta K_{n,k}$ satisfying
\begin{equation}
    \sup_{\mu\in\theta K_{n,k}}
    \E_\mu\|\widehat\mu-\mu\|_F^2
    \leq
    L_n R^*(\theta K_{n,k}).
    \label{eq:minimax-opt}
\end{equation}
Then, from $\widehat\mu$, one can construct in polynomial time an estimator
$\widetilde v=\widetilde v(Y)\in[0,1]^n$ of the planted vector in the
Schramm--Wein model such that
\begin{equation}
    \E\|\widetilde v-v\|_2^2
    =
    o(k)
    =
    o(\rho n).
    \label{eq:vector-small-mse}
\end{equation}
Consequently, after a permutation symmetrization, one obtains a polynomial-time
estimator of a planted coordinate with mean-squared error $o(\rho)$.

Therefore, under the standard low-degree conjecture for planted-submatrix
\emph{estimation}---namely, that the degree-$n^{\Omega(1)}$ estimation barrier
of \cite[Theorem~1.2]{schrammwein2022} reflects a corresponding barrier for
polynomial-time algorithms---there is no polynomial-time estimator satisfying
\eqref{eq:minimax-opt} for any polylogarithmic approximation factor
$L_n=(\log n)^{O(1)}$.
\end{proposition}

\begin{proof}
We first strengthen the statistical estimate
\eqref{eq:small-relative-risk}.  Define
\begin{align*}
    \Delta
    :=
    1-b-2a.
\end{align*}
Since $a<1/2-b/2$ by \eqref{eq:gap-range}, we have $\Delta>0$. Moreover,
\begin{align*}
    \lambda^2 k
    =
    n^{-2a}n^{1-b}
    =
    n^\Delta.
\end{align*}

Recall from \eqref{eq:net-risk} and \eqref{eq:entropy} that, for every
$0<\varepsilon<1/2$,
\begin{align}
    R^*(\theta K_{n,k})
    &\lesssim
    \theta^2\varepsilon^2
    +
    \sigma^2 k
    \left[
        \varepsilon^{-2}\log\frac{en}{k}
        +
        \varepsilon^{-4}\log\frac{C}{\varepsilon}
    \right].
    \label{eq:polylog-risk-start}
\end{align}
For the planted-submatrix model under consideration, $\sigma\asymp1$, and
$\theta=\lambda k_+\asymp\lambda k$. Hence
\begin{align}
    \frac{R^*(\theta K_{n,k})}{\theta^2}
    \lesssim
    \varepsilon^2
    +
    \frac{1}{\lambda^2 k}
    \left[
        \varepsilon^{-2}\log\frac{en}{k}
        +
        \varepsilon^{-4}\log\frac{C}{\varepsilon}
    \right].
    \label{eq:polylog-relative-risk}
\end{align}

Choose
\begin{align*}
    \varepsilon
    :=
    n^{-\Delta/8}.
\end{align*}
For all sufficiently large $n$, $\varepsilon\in(0,1/2)$. Since
\begin{align*}
    \log\frac{en}{k}=O(\log n),
    \qquad
    \log\frac{C}{\varepsilon}=O(\log n),
\end{align*}
and $\lambda^2k=n^\Delta$, \eqref{eq:polylog-relative-risk} yields
\begin{align}
    \frac{R^*(\theta K_{n,k})}{\theta^2}
    &\lesssim
    n^{-\Delta/4}
    +
    n^{-3\Delta/4}\log n
    +
    n^{-\Delta/2}\log n
    \nonumber\\
    &\lesssim
    n^{-\Delta/4}.
    \label{eq:polynomial-relative-risk}
\end{align}
Thus
\begin{equation}
    R^*(\theta K_{n,k})
    \lesssim
    \theta^2 n^{-\Delta/4}.
    \label{eq:strong-small-relative-risk}
\end{equation}

This polynomial slack absorbs every fixed polylogarithmic approximation factor.
Indeed, if $L_n\leq(\log n)^A$ for some fixed $A<\infty$, then
\begin{align*}
    L_n n^{-\Delta/4}
    \longrightarrow
    0.
\end{align*}
Consequently, \eqref{eq:minimax-opt} and
\eqref{eq:strong-small-relative-risk} imply
\begin{equation}
    \sup_{\mu\in\theta K_{n,k}}
    \E_\mu\|\widehat\mu-\mu\|_F^2
    =
    o(\theta^2).
    \label{eq:polylog-small-relative-risk}
\end{equation}

Now let
\begin{align*}
    \mathcal G
    :=
    \{k_-\leq m\leq k_+\},
    \qquad
    m=\|v\|_0.
\end{align*}
By \eqref{eq:support-concentration},
$\Pp(\mathcal G)=1-o(1)$. On $\mathcal G$, the embedding
\eqref{eq:plant-in-body} gives
\begin{align*}
    \mu=\lambda vv^\top\in\theta K_{n,k}.
\end{align*}
Therefore, uniformly over planted vectors satisfying $\mathcal G$,
\begin{equation}
    \E_W
    \bigl\|
        \widehat\mu-\lambda vv^\top
    \bigr\|_F^2
    =
    o(\theta^2)
    =
    o(\lambda^2k^2).
    \label{eq:matrix-small-mse}
\end{equation}
Here the expectation is over the Gaussian noise conditional on $v$.

Since the target matrix is symmetric, we may replace $\widehat\mu$ by its
symmetrization
\begin{align*}
    \widehat M
    :=
    \frac{\widehat\mu+\widehat\mu^\top}{2},
\end{align*}
which is still computable in polynomial time and satisfies
\begin{align*}
    \|\widehat M-\lambda vv^\top\|_F
    \leq
    \|\widehat\mu-\lambda vv^\top\|_F.
\end{align*}
Let $\widehat u$ be a leading unit eigenvector of $\widehat M$, and write
\begin{align*}
    u:=\frac{v}{\sqrt m}.
\end{align*}
The matrix $\lambda vv^\top$ has a single nonzero eigenvalue $\lambda m$.
Combining Davis--Kahan with the trivial bound on the sine of the angle gives,
for an absolute constant $C_0$,
\begin{equation}
    1-|\langle\widehat u,u\rangle|^2
    \leq
    C_0
    \frac{
        \|\widehat M-\lambda vv^\top\|_F^2
    }{
        \lambda^2m^2
    }.
    \label{eq:DK-mean}
\end{equation}
Since $m\asymp k$ on $\mathcal G$, taking conditional expectations and using
\eqref{eq:matrix-small-mse} yields
\begin{equation}
    \E_W
    \left[
        1-|\langle\widehat u,u\rangle|^2
    \right]
    =
    o(1)
    \qquad
    \text{uniformly on }\mathcal G.
    \label{eq:angle-small}
\end{equation}

We now convert this normalized eigenvector estimate directly into an estimator
of the planted vector.  Let $|\widehat u|$ denote coordinatewise absolute
value, and define
\begin{equation}
    \widetilde v
    :=
    \Pi_{[0,1]^n}
    \bigl(
        \sqrt{k}\,|\widehat u|
    \bigr),
    \label{eq:vtilde-def}
\end{equation}
where the projection onto $[0,1]^n$ is coordinatewise clipping.  This
construction removes the irrelevant global sign of $\widehat u$ and is
polynomial time.  Since $u$ is entrywise nonnegative,
\begin{align*}
    \bigl\||\widehat u|-u\bigr\|_2
    \leq
    \min_{s\in\{-1,1\}}
    \|s\widehat u-u\|_2.
\end{align*}
For unit vectors,
\begin{align*}
    \min_{s\in\{-1,1\}}
    \|s\widehat u-u\|_2^2
    =
    2\bigl(1-|\langle\widehat u,u\rangle|\bigr)
    \leq
    2\bigl(1-|\langle\widehat u,u\rangle|^2\bigr).
\end{align*}
Also, projection onto $[0,1]^n$ cannot increase the distance to
$v\in\{0,1\}^n$.  Hence, on $\mathcal G$,
\begin{align}
    \|\widetilde v-v\|_2^2
    &\leq
    \bigl\|
        \sqrt{k}\,|\widehat u|-\sqrt m\,u
    \bigr\|_2^2
    \nonumber\\
    &\leq
    2k\bigl\||\widehat u|-u\bigr\|_2^2
    +
    2(\sqrt{k}-\sqrt m)^2.
    \label{eq:vtilde-error}
\end{align}
By \eqref{eq:angle-small}, the conditional expectation of the first term is
$o(k)$.  Furthermore, on $\mathcal G$ we have
$m/k=1+O(\delta_n)$, and therefore
\begin{align*}
    (\sqrt{k}-\sqrt m)^2
    =
    O(\delta_n^2k)
    =
    o(k).
\end{align*}
It follows that
\begin{equation}
    \E\left[
        \|\widetilde v-v\|_2^2
        \mathbf 1_{\mathcal G}
    \right]
    =
    o(k).
    \label{eq:good-event-v-risk}
\end{equation}

It remains only to control the negligible event $\mathcal G^c$.  Since
coordinatewise clipping can only decrease the Euclidean norm,
\begin{align*}
    \|\widetilde v\|_2^2\leq k,
\end{align*}
and hence
\begin{align*}
    \|\widetilde v-v\|_2^2
    \leq
    2k+2m.
\end{align*}
Thus
\begin{align*}
    \E\left[
        \|\widetilde v-v\|_2^2
        \mathbf 1_{\mathcal G^c}
    \right]
    &\leq
    2k\,\Pp(\mathcal G^c)
    +
    2\E\left[m\mathbf 1_{\mathcal G^c}\right]\\
    &\leq
    2k\,\Pp(\mathcal G^c)
    +
    2(\E m^2)^{1/2}\Pp(\mathcal G^c)^{1/2}
    =
    o(k),
\end{align*}
where we used $\E m^2\asymp k^2$ and
$\Pp(\mathcal G^c)=o(1)$.  Combining this with
\eqref{eq:good-event-v-risk} proves \eqref{eq:vector-small-mse}.

We now show a lemma about symmmetrizing the estimator.

\begin{lemma}[Deterministic permutation averaging]
Consider the planted-submatrix model
\begin{align*}
    Y=\lambda vv^\top+W,
    \qquad
    v_i\stackrel{\mathrm{iid}}{\sim}\operatorname{Bernoulli}(\rho),
\end{align*}
where the distribution of $W$ is invariant under simultaneous permutations of
its rows and columns. Let $A(Y)\in\RR^n$ be any estimator of $v$.

For each $i\in[n]$, choose a fixed permutation $\pi_i\in S_n$ satisfying
\begin{align*}
    \pi_i(1)=i,
\end{align*}
and let $P_{\pi_i}$ denote its permutation matrix. Define
\begin{align*}
    A^{\pi_i}(Y)
    :=
    P_{\pi_i}^\top
    A\!\left(P_{\pi_i}YP_{\pi_i}^\top\right),
\end{align*}
and define the scalar estimator
\begin{equation}
    \overline A_1(Y)
    :=
    \frac1n
    \sum_{i=1}^n
    A^{\pi_i}_1(Y).
    \label{eq:deterministic-permutation-average}
\end{equation}
Then
\begin{equation}
    \EE
    \left[
        \bigl(\overline A_1(Y)-v_1\bigr)^2
    \right]
    \leq
    \frac1n
    \EE
    \|A(Y)-v\|_2^2.
    \label{eq:deterministic-coordinate-risk}
\end{equation}
In particular, if
\begin{align*}
    \EE\|A(Y)-v\|_2^2=o(\rho n),
\end{align*}
then
\begin{align*}
    \EE
    \left[
        \bigl(\overline A_1(Y)-v_1\bigr)^2
    \right]
    =
    o(\rho).
\end{align*}
Moreover, if $A$ is computable in polynomial time, then so is
$\overline A_1$.
\end{lemma}

\begin{proof}
For each $i\in[n]$, define
\begin{align*}
    v^{\pi_i}:=P_{\pi_i}v,
    \qquad
    Y^{\pi_i}:=P_{\pi_i}YP_{\pi_i}^\top.
\end{align*}
Since the coordinates of $v$ are i.i.d.\ Bernoulli$(\rho)$, and since the
noise distribution is invariant under simultaneous permutations of the rows
and columns,
\begin{equation}
    (Y^{\pi_i},v^{\pi_i})
    \stackrel{d}{=}
    (Y,v).
    \label{eq:deterministic-permutation-invariance}
\end{equation}
By the definition of $A^{\pi_i}$,
\begin{align*}
    A^{\pi_i}_1(Y)
    =
    \bigl[
        P_{\pi_i}^\top A(Y^{\pi_i})
    \bigr]_1
    =
    A_{\pi_i(1)}(Y^{\pi_i})
    =
    A_i(Y^{\pi_i}),
\end{align*}
while
\begin{align*}
    (v^{\pi_i})_i
    =
    v_1.
\end{align*}
Therefore, using \eqref{eq:deterministic-permutation-invariance},
\begin{align}
    \EE
    \left[
        \bigl(A^{\pi_i}_1(Y)-v_1\bigr)^2
    \right]
    &=
    \EE
    \left[
        \bigl(
            A_i(Y^{\pi_i})-(v^{\pi_i})_i
        \bigr)^2
    \right]
    \nonumber\\
    &=
    \EE
    \left[
        \bigl(
            A_i(Y)-v_i
        \bigr)^2
    \right].
    \label{eq:permuted-coordinate-risk}
\end{align}

Now, by convexity of $x\mapsto x^2$,
\begin{align*}
    \EE
    \left[
        \bigl(\overline A_1(Y)-v_1\bigr)^2
    \right]
    &=
    \EE
    \left[
        \left(
            \frac1n
            \sum_{i=1}^n
            \bigl(A^{\pi_i}_1(Y)-v_1\bigr)
        \right)^2
    \right]\\
    &\leq
    \frac1n
    \sum_{i=1}^n
    \EE
    \left[
        \bigl(A^{\pi_i}_1(Y)-v_1\bigr)^2
    \right]\\
    &=
    \frac1n
    \sum_{i=1}^n
    \EE
    \left[
        \bigl(A_i(Y)-v_i\bigr)^2
    \right]\\
    &=
    \frac1n
    \EE
    \|A(Y)-v\|_2^2,
\end{align*}
where the penultimate equality follows from
\eqref{eq:permuted-coordinate-risk}. This proves
\eqref{eq:deterministic-coordinate-risk}.

Finally, $\overline A_1$ requires only $n$ evaluations of $A$, together
with permutations of the rows and columns of $Y$. Hence if $A$ is
polynomial-time computable, then $\overline A_1$ is also polynomial-time
computable.
\end{proof}

\paragraph{Application to the planted vector.}
Apply the lemma to the polynomial-time estimator $\widetilde v(Y)\in[0,1]^n$
constructed above from the leading eigenvector of the matrix estimator.
We have already shown that
\begin{align*}
    \EE
    \|\widetilde v-v\|_2^2
    =
    o(k)
    =
    o(\rho n).
\end{align*}
Therefore, defining
\begin{align*}
    \overline v_1(Y)
    :=
    \frac1n
    \sum_{i=1}^n
    \left[
        P_{\pi_i}^\top
        \widetilde v\!\left(
            P_{\pi_i}YP_{\pi_i}^\top
        \right)
    \right]_1,
\end{align*}
the lemma gives
\begin{align*}
    \EE
    \left[
        \bigl(\overline v_1(Y)-v_1\bigr)^2
    \right]
    =
    o(\rho).
\end{align*}
Moreover, $\overline v_1$ is computable in polynomial time, since it requires
only $n$ calls to the polynomial-time procedure producing $\widetilde v$.
Thus a polynomial-time minimax-rate estimator over $\theta K_{n,k}$ would
yield a polynomial-time estimator of the planted coordinate $v_1$ with
mean-squared error $o(\rho)$. Under the standard low-degree conjecture
interpreting the Schramm--Wein degree-$n^{\Omega(1)}$ estimation lower bound as a polynomial-time barrier, this gives the desired contradiction.
\end{proof}

\begin{remark}[Scope of the Schramm--Wein result]
The final step is deliberately stated conditionally.  The theorem of Schramm
and Wein is a rigorous lower bound for degree-$n^{\Omega(1)}$ polynomial
estimators; it does not by itself rule out every polynomial-time algorithm.
Thus the complexity-theoretic conclusion for the convex-body example is
conditional on the standard low-degree conjecture for planted-submatrix
\emph{estimation}.
\end{remark}

\begin{remark}[Why the construction is useful]
The point of $K_{n,k}$ is that the computational difficulty is not hidden in
the representation of the constraint set.  Its gauge is the SDP
\eqref{eq:gauge}, while the metric-entropy bound \eqref{eq:entropy} shows that
the entire convex relaxation remains statistically estimable at a scale well
inside the conjectured computationally hard regime for planted-submatrix
estimation. Hence the example separates efficient access to the convex body
from efficient minimax-rate estimation over that body.
\end{remark}

\end{document}